\documentclass[a4paper,12pt]{amsbook}
\usepackage[english]{babel}

\usepackage[latin1]{inputenc}
\usepackage{mathrsfs}
\usepackage{amsmath}
\usepackage{amsthm}
\usepackage{amssymb}
\usepackage{indentfirst}
\usepackage{type1cm}

\usepackage{ifthen}
\newcommand{\cit}[1]{{\rm \textbf{#1}}}

\newcommand{\Ref}[2]{\cit{%
\ifthenelse{\equal{#1}{thm}}{Theorem}{}%
\ifthenelse{\equal{#1}{cap}}{Chapter}{}%
\ifthenelse{\equal{#1}{prop}}{Proposition}{}%
\ifthenelse{\equal{#1}{lem}}{Lemma}{}%
\ifthenelse{\equal{#1}{cor}}{Corollary}{}%
\ifthenelse{\equal{#1}{defn}}{Definition}{}%
\ifthenelse{\equal{#1}{oss}}{Remark}{}%
\ifthenelse{\equal{#1}{sec}}{Section}{}%
\ifthenelse{\equal{#1}{ex}}{Example}{}%
\ifthenelse{\equal{#1}{conj}}{Conjecture}{}%
\ifthenelse{\equal{#1}{ssec}}{Subsection}{}%
\ifthenelse{\equal{#1}{tab}}{Table}{}%
\ifthenelse{\equal{#1}{cla}}{Claim}{}%
\ifthenelse{\equal{#1}{ques}}{Question}{}%
\  \ref{#1:#2}%
}}

\newcommand{\hk}{Hyperk\"{a}hler }
\newcommand{\kahl}{K\"{a}hler }
\newcommand{\ktipo}{$K3^{[2]}$-type }
\newcommand{\kntipo}{$K3^{[n]}$-type }
\newcommand{\ie}{i.~e.~}

\theoremstyle{plain} 
\newtheorem{prop}{Proposition}[section]
\newtheorem{thm}[prop]{Theorem}
\newtheorem{lem}[prop]{Lemma} 
\newtheorem{cor}[prop]{Corollary}
\newtheorem{conj}[prop]{Conjecture}
\newtheorem{ques}[prop]{Question}

\theoremstyle{remark}
\newtheorem{oss}[prop]{Remark}
\newtheorem{ex}[prop]{Example}

\theoremstyle{definition}
\newtheorem{defn}[prop]{Definition}
\title{Automorphisms of \hk manifolds}
\address{Ph.D thesis}
\author{Giovanni Mongardi}
\numberwithin{equation}{chapter}
\numberwithin{section}{chapter}

\makeindex
\begin{document}

\frontmatter
\maketitle
\tableofcontents
\numberwithin{prop}{chapter}
\chapter*{Introduction}\label{cap:intro}



This thesis is devoted to study certain interesting properties of \hk manifolds and of their automorphism groups. \hk manifolds have been studied mainly due to their appearance in the famous Bogomolov's decomposition theorem. This theorem states that any manifold with a Ricci flat metric is, up to a finite cover, a direct product of complex tori, Calabi-Yau manifolds and \hk manifolds. By Yau's proof of Calabi's conjecture having a Ricci flat metric is equivalent to having trivial first Chern class.\\
The field of \hk geometry is quite recent, although a classical example consists in K3 surfaces\footnote{All \hk manifolds of dimension 2 are $K3$ surfaces.}. The first higher dimensional \footnote{\ie of dimension greater than 2.} examples where found by Fujiki \cite{fuj} and by Beauville \cite{beau3}. They consist of the Hilbert scheme\footnote{See \Ref{ex}{kntipo}.} of length $n$ subschemes on a $K3$ surface $S$, denoted $S^{[n]}$, and of generalized Kummer manifolds\footnote{See \Ref{ex}{gen_kum_tipo}.}. We remark that the generic deformation of $S^{[n]}$ for any $K3$ surface $S$ is not the Hilbert scheme on another $K3$ surface. We will call elements of this deformation class Manifolds of \kntipo.\\ Fujiki, Beauville and Bogomolov developed much of the theory concerning the second cohomology of \hk  manifolds, proving the existence of what is commonly known as Beauville-Bogomolov form (or also Fujiki-Beauville-Bogomolov in Japanese literature). An interesting feature of \hk manifolds is that any family of \hk manifolds has a dense subset consisting of \emph{projective} \hk manifolds, therefore in their study it is possible to apply both analytical and geometric methods. Some results, who are apparently deeply algebraic in nature, such as \Ref{prop}{proj_dense}, have a complex analytical proof. While other results, such as \Ref{thm}{moving_cone}, deeply use recent progress in minimal model program.\\
For quite some time all known examples consisted of manifolds deformation equivalent to those found by Beauville \cite{beau3}, with interesting projective examples as Fano variety of lines on cubic fourfolds (Beauville and Donagi, \cite{bedon}), Double covers of certain special sextic fourfolds (O'Grady, \cite{ogr}), variety of sums of powers of cubic fourfolds (Iliev and Ranestad, \cite{ir1}) and subspaces of certain grassmanians (Debarre and Voisin, \cite{dv}).\\
A notable impulse to the research in this field is due to the discovery made by Mukai \cite{muk_f} of a symplectic form on moduli spaces of certain sheaves on symplectic surfaces. This fact led to the hope that new \hk manifolds could be found with these construction and a good theory was developed by various mathematicians; for a complete set of references the interested reader can consult \cite{hl_moduli}. However it has been proved that all nonsingular \hk manifolds obtained in this way were a deformation of known examples and the singular ones had a resolution of singularities which is \hk just in two cases, namely in O'Grady's six dimensional manifold \cite{ogr3} and in O'Grady's ten dimensional manifold \cite{ogr4}.\\
Recently a long standing question of \hk geometry has been partially resolved, namely the Torelli problem. Verbitsky \cite{ver}, Markman \cite{mar1} and Huybrechts \cite{huy_tor} have proven theorems explaining to which extent a \hk manifold can be recovered by its integral Hodge structure on its second cohomology. This result will be instrumental in our work and will allow, under some hypothesis, to construct a group of birational transformation on a \hk manifolds from a group of isometries on a lattice.\\
 In recent years there have been several works concerning automorphisms of Hyperk\"{a}hler manifolds, starting from the foundational work of Nikulin \cite{nik1}, Mukai \cite{muk} and Kondo \cite{kon} and an explicit example of Morrison \cite{mor} in the case of K3 surfaces. Then isolated examples of automorphisms of higher dimensional \hk manifolds were given by Namikawa \cite{nami}, by Beauville \cite{beau2} and later by Kawatani \cite{kaw} and Amerik \cite{amerik}. Some further work was done by Boissi\`{e}re, Nieper-Wi\ss kirchen and Sarti \cite{boiniesar}, which also paved the way for a generalization of the notion of Enriques surface, independently developed also by Oguiso and Schr\"{o}er in \cite{ogu1} and \cite{ogu2}. Some general work on the automorphisms and birationalities was done by Oguiso \cite{ogu3} and by Boissi\`{e}re \cite{boi} and yet again recently by Boissi\`{e}re and Sarti \cite{bs} while order 2 automorphisms were fully analyzed by Beauville if they are antisymplectic \cite{beau} and partially analyzed by Camere (\cite{cam}, where also an exhaustive list of examples can be found). Before those works on involutions came the work of O'Grady (\cite{ogr} and \cite{ogr2}) on Double-EPW sextics which are naturally endowed with an antisymplectic involution and form a family of the maximal dimension for such involutions.\\
\newpage
\section*{Overview of the results}
In order to present our results we introduce some definitions:
\begin{defn}\label{defn:def_equiv}
Let $X,Y$ be two deformation equivalent \hk manifolds and let $G\subset Aut(X)$, $H\subset Aut(Y)$. Then $(X,G)$ is deformation equivalent to $(Y,H)$ if $G\cong G'\cong H$ and there exists a flat family $\mathcal{X}\rightarrow B$ and two maps $\{a\}\rightarrow B$, $\{b\}\rightarrow B$ such that $\mathcal{X}_a\cong X$ and $\mathcal{X}_b\cong Y$. Moreover we require that there exists a faithful action of the group $G'$ on $\mathcal{X}$ inducing fibrewise faithful actions of $G'$ such that its restriction to $\mathcal{X}_a$ and $\mathcal{X}_b$ coincides with $G$ and $H$.
\end{defn}
\begin{defn}\label{defn:natural}
Let $S$ be a K3 surface and let $G\subset Aut(S)$. Then $G$ induces automorphisms on $S^{[n]}$ which are called \emph{natural}. We will also call $(S^{[n]},G)$ a natural couple.
\end{defn}
\begin{defn}\label{defn:standard_intro}
Let $X$ be a manifold of \kntipo and let $G\subset Aut(X)$. Then the couple $(X,G)$ is standard if it is deformation equivalent to $(S^{[n]},H)$, where $S$ is a K3 surface and $H$ is a group of natural automorphisms. We call $G$ exotic otherwise.
\end{defn}
\begin{defn}
Let $X$ be a \hk manifold, we define \\$Aut_s(X)\subset Aut(X)$ \index{Group, symplectic automorphisms of $X$, $Aut_s(X)$} as the subgroup of the automorphism group that preserves the symplectic form on $X$. We call its elements symplectic automorphisms and we call it the group of symplectic automorphisms.
\end{defn}
 
The results contained in the present thesis can be grouped in 3 deeply connected areas: \emph{new exotic symplectic automorphisms}, \emph{standardness of known automorphisms} and a \emph{classification of prime order symplectic automorphisms} on a wide class of \hk manifolds.\\
New exotic symplectic automorphisms can be found in \Ref{cap}{exoexa}. We remark that a symplectic automorphism on a $K3$ surface has at most order $8$, we have written down two examples of order 11 automorphisms on manifolds of \ktipo. Namely one is defined on the Fano scheme of lines of a cubic fourfold and the other on a double EPW-sextic\footnote{See \Ref{ssec}{epw}}. Moreover we have given also an example of an order 15 symplectic automorphism again on the Fano scheme of lines of a cubic fourfold. Mukai \cite{muk} proved that a group of symplectic automorphisms on a K3 surface has order at most 960, however the situation is very different on manifolds of \ktipo as we have given an example with a group of symplectic automorphisms of order 2520 and one of order 29160.\\
For what concerns deformational behaviour of symplectic automorphisms we have proven that any couple $(X,\varphi)$ consisting of a manifold of \ktipo and a symplectic automorphism of order 2 or 5 is standard. Moreover we proved that points corresponding to natural couples are dense in the moduli space of manifolds of \ktipo having a symplectic automorphism of order 2 and 5. The same result also holds for manifolds of \ktipo with a symplectic automorphism of order 3 if this automorphism satisfies a condition on the fixed locus $X^\varphi$.\\
In the last part of the thesis we take a lattice theoretic approach, much in the spirit of what Nikulin \cite{nik1}, Mukai \cite{muk} and Kondo \cite{kon} did in the case of $K3$ surfaces. We work towards classification results. For a general \hk manifold of the known types we can only prove some limitations on the possible prime orders of a symplectic automorphism (or, more precisely, on the order of the induced Hodge isometry on $H^2$). If we specialize to manifolds of \kntipo we can give a full classification of all possible prime order symplectic automorphisms and also of their Co-invariant lattices\footnote{In the sense of \Ref{defn}{inv_locus}}. In the case of $K3$ surfaces we have that the maximal prime order is 7, while in this case it is 11. Moreover we prove also a theorem stating sufficient conditions to give a group of symplectic birational transformations of a manifold of \kntipo from a group of isometries on a Niemeier lattice. Finally, restricting even more to manifolds of \ktipo, we improve our classification result giving all prime order symplectic automorphisms together with their fixed locus and their co-invariant lattice. We wish to stress that this does not give a result on the number of deformation classes of couples $(X,\varphi)$, where $X$ is of \ktipo and $\varphi$ is symplectic of prime order.\\
Now let us briefly talk about the interplay between these three kind of results: our classification of prime order symplectic automorphisms is made more precise thanks to our results on the standardness of automorphisms, namely because it allows to prove that in these cases there is only one deformation equivalence class. On the other hand our examples can be improved by our classification and by the lattice theoretic approach, since it allows us to compute the Picard lattices of these examples as we did in \Ref{sec}{exa2}. Moreover it allows also to determine, in the case of order 11 automorphisms, the number of deformation equivalence classes in dimension 4. There are 2 such deformation classes and we give a projective element of both of them among our examples. Finally, using our classification, the standardness question proven in \Ref{thm}{standard_morph} can be reformulated with less conditions as we did in \Ref{cor}{prime_stand}.

\section*{Structure of the thesis}
\Ref{cap}{hk} provides a survey on several well known results on \hk manifolds, almost all of the material contained here is present in the literature, apart for \Ref{prop}{vgir_mail} whose proof was communicated to me by Prof. B. Van Geemen \cite{vgir} and is due to him, A. Iliev and K. Ranestad.\\
\Ref{cap}{lattice} contains various results on lattices which are instrumental in the analysis of symplectic automorphisms on \hk manifolds. Some of the material contained here is classical, like Niemeier's list of unimodular lattices and Nikulin's results on discriminant forms. However a portion of this chapter contains some results I could not find in the literature, like \Ref{sec}{prime_nieme}.\\
\Ref{cap}{k3_case} gives a brief overview on symplectic automorphisms of K3 surfaces, illustrating the classical results of Nikulin and Mukai. We give explicit proofs of most of the material contained here. Whenever possible these proofs are adapted to exploit only the \hk structure of $K3$ surfaces, providing an easier generalization.\\
\Ref{cap}{exoexa} gathers a series of examples of symplectic automorphisms on \hk manifolds, it focuses on manifolds of \ktipo but some more examples are given. This chapter contains also some examples of exotic automorphisms, \ie automorphisms which are not obtained as deformations of the automorphism group induced on the Hilbert scheme of points of a $K3$ surface by the underlying automorphism group of the surface itself.\\
\Ref{cap}{deformations} is devoted to establish whether a symplectic automorphism on a manifold of \ktipo can be obtained through the above cited process of deformation of an induced automorphism group. 
In the literature automorphisms of $S^{[2]}$ induced by those on $S$ are called \emph{Natural}\footnote{See \cite{boi} and \Ref{defn}{natural}}. Often the same terminology is used for automorphisms obtained as deformations of $(S^{[2]},\varphi^{[2]})$, however we prefer to use the term \emph{Standard}\footnote{See \Ref{defn}{standard_intro}} to denote automorphisms obtained by deformation. This chapter proves that the deformation-theoretical question on the standardness of a symplectic automorphism of order 2,3 or 5 on a manifold of \ktipo is only a cohomological condition: such an automorphism is standard as soon as its action on the second cohomology is the same of that of a natural symplectic automorphism.\\
\Ref{cap}{fixed} consists of a series of computations providing a generalization of \Ref{sec}{fix_k3} in the case of manifolds of \ktipo. Let us stress that some of the techniques used on $K3$ surfaces do not work in the higher dimensional case: in dimension 2 the fixed locus of a symplectic automorphism $\varphi$ consists of isolated points, therefore there is a crepant resolution of $X/\varphi$ which is again a $K3$ surface. This allows to compute the possible orders of $\varphi$ and the number of fixed points. However in higher dimensions usually there is no crepant resolution of the quotient, it only exists if the fixed locus has pure codimension 2. Furthermore in this case the resolution is still a \hk manifold which need not be deformation equivalent to the one we started with. At present it is not possible to provide a full generalization of \Ref{sec}{fix_k3} mainly because the computations are quite hard and, as the case of order 3 shows in \Ref{thm}{fixedp3}, these computations provide only a series of \emph{possible} fixed loci for a symplectic automorphism without giving any hint on the existence of such a morphism. Therefore the use of the methods contained in the following chapter are preferrable.\\
\Ref{cap}{groups_sporadic} contains the main general results of the present work: we are able to establish a connection between finite groups of symplectic automorphisms on \hk manifolds and isometries of certain well known lattices. In the particular case of \hk manifolds of \kntipo we can embed finite groups of symplectic automorphisms in the sporadic simple group\footnote{See \Ref{defn}{co1} or \cite{con2}.} $Co_1$ and, to some extent, we provide also a converse. Moreover we give also a classification result for prime order symplectic automorphisms on manifolds of \kntipo. In the particular case of manifolds of \ktipo our classification provides indeed all known examples, however for higher dimensions there are no known examples for 3 of the possible cases appearing in \Ref{tab}{prime_autom_k3n_tab}. 

\section*{Notations}
In this section we gather several definitions that will be used throughout the rest of the paper, most of our notation is standard, apart for \Ref{defn}{leech_24dim} where we define the Leech lattice as a negative definite lattice instead of a positive definite one. Also \Ref{defn}{trascend} is nonstandard but coincides with the standard definition in the projective case.
\begin{defn}
Let $R$ be a free $\mathbb{Z}$-module and let $(\,,\,)_R\,:\,R\times R\,\rightarrow\,\mathbb{Z}$ be a bilinear pairing. We call the couple $(R,(\,,\,)_R)$ \index{Lattice, $(R,(\,,\,)_R)$} a lattice and we denote it by $R$ whenever the pairing is understood.
\end{defn}
If the pairing takes values in $2\mathbb{Z}$ we will say that the lattice is even, odd otherwise.
\begin{defn}
A pair $(R,(\,,\,)_R)$ is called generalized lattice whenever $R$ is as above and the pairing takes values over $\mathbb{Q}$.\index{Lattice, Generalized} 
\end{defn}
Moreover we denote $R(n)$ \index{Lattice, Multiple of a, $(R,n(\,,\,)_R)$, $R(n)$,} the lattice $R$ with pairing multiplied by $n$ and we denote $(n)$ \index{Lattice, rank 1, $(n)$} the lattice $(\mathbb{Z},q)$ with $q(1)=n$. We also call a $n$ vector an element of square $n$ inside a lattice. We will say that a lattice $R\subset L$ is primitive if the quotient $L/R$ is torsion free. If on the other hand the quotient is a finite group we say that $L$ is an overlattice of $R$. We will say that a lattice $R$ represents an integer $n$ if there exists a (primitive) element of $R$ with square $n$.\\  
We will show that a \hk manifold $X$ has a lattice structure on its second cohomology and we will often denote $(\,,\,)_X$ the pairing and $q_X$ \index{Quadratic form on a \hk manifold, $q_X$} the induced quadratic form. If on the other hand the pairing is understood we will denote $q_X(e)=e^2$ for all $e\in H^2(X,\mathbb{Z})$.
\begin{defn}
Let $R$ be a lattice and let $e\in R$. We denote $div(e)_R=n$ \index{Lattice, divisibility of an element, $div(e)_R,\,div(e)$} if $(e,R)=n\mathbb{Z}$ and we say that $e$ is $n$-divisible. If the lattice $R$ is understood we just denote it $div(e)$.
\end{defn}
\begin{defn}\label{defn:leech_24dim}
The Leech lattice $\Lambda$ \index{Lattice, Leech, $\Lambda$} is the unique negative definite unimodular lattice of rank 24 that does not contain any element of square $-2$. We will also provide explicit definitions in \Ref{ex}{leech_24} and \Ref{ex}{leech_26}.
\end{defn}
\begin{defn}\label{defn:co1}
Let $Co_0=Aut(\Lambda)$ \index{Group, Conway's dotto, $Co_0$} be the automorphism group of the Leech lattice and let $Co_1=Co_0/(\pm Id)$ \index{Group, Conway's first simple, $Co_1$} be its quotient by its center. It is a well known fact (see \cite{atlas} and \cite{con2}) that $Co_1$, usually called Conway's first sporadic group, is a simple group.
\end{defn}

\begin{defn}\label{defn:trascend}
Let $X$ be a symplectic manifold. Then we define the transcendental part $T(X)$ \index{Transcendental part, $T(X)$}as the smallest integral Hodge structure containing the symplectic form $\sigma_X$. If $X$ is \hk there is a quadratic form on $H^2(X,\mathbb{Z})$ and we will denote $S(X)=T(X)^\perp$ \index{Algebraic part, $S(X)$}.
\end{defn}

\numberwithin{prop}{section}
\mainmatter
\setcounter{prop}{0}
\chapter{\hk  manifolds}\label{cap:hk}
This chapter gathers several known results on \hk manifolds and provides an introductory guide to such manifolds. Many of these results are taken from the survey of Huybrechts \cite{huy}.\\ Obviously we start with the following:
\begin{defn}
Let $X$ be a \kahl manifold, it is called a irreducible holomorphic symplectic manifold if the following hold:
\begin{itemize}
\item $X$ is compact.
\item $X$ is simply connected.
\item $H^{2,0}(X)=\mathbb{C}\sigma_X$, where $\sigma_X$ \index{Symplectic form, $\sigma_X$} is an everywhere nondegenerate symplectic 2-form.
\end{itemize}
\end{defn}
\begin{defn}
Let $X$ be a \kahl manifold and let $i,j,k$ be three complex \kahl metrics such that $Re(i)=Re(j)=Re(k)=g$. Then the Riemannian metric $g$ is called \hk if the three complex structures $I,J,K$ induced by $i,j,k$ respectively on $TX$ satisfy $IJ=K$.
\end{defn}

\begin{defn}
Let $X$ be a \kahl manifold, it is called a \hk manifold if the following hold:
\begin{itemize}
\item $X$ is compact.
\item $X$ is simply connected.
\item There exists a \hk metric $g$ on $X$.
\end{itemize}
\end{defn}
Often in the literature \hk manifolds denotes just manifolds with a \hk metric, without requiring compactness and simple connectedness.
At first sight a \hk manifold seems quite different from an irreducible holomorphic symplectic manifold but Yau's proof of Calabi's Conjecture \cite{yau} can be used to associate a \hk metric to any \kahl class on an irreducible symplectic holomorphic manifold, therefore we will not distinguish between the two definitions.\\
There are not many known examples of \hk manifolds, for a long time the only known \hk manifolds were $K3$ surfaces (which are the only example in dimension 2), but two families of examples were given by Beauville \cite{beau3}:
\begin{ex}\label{ex:kntipo}
Let $S$ be a $K3$ surface and let $S^{(n)}$ be its $n$-th symmetric product. There exists a minimal resolution of singularities 
\begin{equation}\nonumber
S^{[n]}\,\stackrel{HC}{\rightarrow}\,S^{(n)},
\end{equation}
where $S^{[n]}$ \index{Douady space or Hilbert scheme of points, $S^{[n]}$} is the Douady space parametrizing zero dimensional analytic subsets of $S$ of length $n$. Furthermore this resolution of singularities endows $S^{[n]}$ with a symplectic form induced by the symplectic form on $S$. Moreover if $n\geq2$ we have $b_2(S^{[n]})=23$.\\ The case $n=2$ was first studied by Fujiki \cite{fuj}, notice that in this case the resolution of singularities is simply the blow-up along the diagonal.\\ Whenever $X$ is a \hk manifold deformation of one of these manifolds we will call $X$ of $K3^{[n]}$-type.
\end{ex}
\begin{ex}\label{ex:gen_kum_tipo}
Let $T$ be a complex torus and let
\begin{equation}\nonumber
T^{[n+1]}\,\stackrel{HC}{\rightarrow}\,T^{(n+1)}
\end{equation}
be the minimal resolution of singularities of the symmetric product. As in \Ref{ex}{kntipo} Beauville proved that the symplectic form on $T$ induces a symplectic form on $T^{[n+1]}$. However this manifold is not \hk since it is not simply connected, but if we consider
\begin{align}
T^{(n+1)} \stackrel{\Sigma}{\rightarrow}& T\\\nonumber
(t_1,\dots,t_{n+1}) \rightarrow & \sum_i t_i.
\end{align}
And we set $K_n(T)=HC^{-1}\circ \Sigma^{-1}(0)$ \index{Generalized Kummer manifold, $K_n(T)$} we obtain a new \hk manifold called Generalized Kummer manifold of $T$. If $n=1$ then $K_n(T)$ is just the usual Kummer surface, otherwise it has $b_2=7$.\\
Whenever $X$ is a \hk manifold deformation of one of these manifolds we will call $X$ of Kummer $n$-type.
\end{ex}
Two more examples of \hk manifolds are known and they were both discovered by O'Grady (see \cite{ogr3} and \cite{ogr4}), we do not give a precise definition but we will call $Og_6$ \index{O'Grady's six dimensional manifold, $Og_6$} the 6-dimensional example and $Og_{10}$ \index{O'Grady's ten dimensional manifold, $Og_{10}$} the 10 dimensional example. It is known that $b_2(Og_6)=8$ and $b_2(Og_{10})=24$.

\begin{ex}
Let $X$ be a \hk manifold of dimension $2n>2$ and let $\mathbb{P}^n\cong P\subset X$. Then there exists a birational map, called Mukai flop \cite{muk_f}, defined as follows: let $Z$ be the blowup of $X$ along $P$ and let $D$ be the exceptional divisor. The projection $D\rightarrow P$ is isomorphic to the projective bundle $\mathbb{P}(\mathcal{N}_{P|X})\cong\mathbb{P}(\Omega_P)\rightarrow \mathbb{P}^n$. It is a well known fact that this projective bundle has a second projection $D\,\rightarrow\,(\mathbb{P}^{n})^\vee$ and this gives a blowdown $Z\,\rightarrow X'$ on a smooth manifold $X'$ such that $D$ is the exceptional divisor. Moreover if $X'$ is \kahl it is also \hk.
\end{ex}
\section{Cohomology of \hk manifolds}
This section is devoted to illustrate the peculiar nature of the cohomology of a \hk manifold. The first interesting fact is the following:
\begin{thm}
Let $X$ be a \hk manifold of dimension $2n$. Then there exists a canonically defined pairing $(\,,\,)_X$ \index{Beauville-Bogomolov form, $(\,,\,)_X$} on $H^2(X,\mathbb{C})$, the Beauville-Bogomolov pairing, and a constant $c_X$ \index{Fujiki constant, $c_X$} (the Fujiki constant) such that the following holds:
\begin{equation}
\small{(\alpha,\alpha)_X=c_X\left(n/2\int_X\alpha^2(\sigma\overline{\sigma})^{n-1}+(1-n)(\int_X \alpha\sigma^{n-1}\overline{\sigma}^n)(\int_X \alpha\sigma^{n}\overline{\sigma}^{n-1})\right).}
\end{equation}
Here $\sigma$ is a symplectic form such that $\int_X(\sigma\overline{\sigma})^n=1$.
Moreover $c_X$ and $(\,,\,)_X$ are deformation and birational invariants.
\end{thm}
This fact is quite striking and unexpected for $n\geq 2$, furthermore the Beauville-Bogomolov pairing can be used to define a lattice on $H^2(X,\mathbb{Z})$ of signature $(3,b_2(X)-3)$. These lattices have been studied by Beauville \cite{beau3} for manifolds of $K3^{[n]}$-type or of Kummer $n$-type. The cases of O'Grady's examples were studied by Rapagnetta (\cite{rap} and \cite{rap2}).
\begin{ex}
Let $X$ be a \hk manifold of $K3^{[n]}$-type. Then $H^2(X,\mathbb{Z})$ endowed with its Beauville-Bogomolov pairing is isomorphic to the lattice \index{Lattice on a manifold of \kntipo, $L_n$}
\begin{equation}\label{latticeK3n}
L_n:=U\oplus U\oplus U\oplus E_8(-1)\oplus E_8(-1)\oplus (2-2n).
\end{equation}
Where $U$ is the hyperbolic lattice, $E_8(-1)$ is the unique unimodular even negative definite lattice of rank $8$, $(2-2n)$ is $(\mathbb{Z},q)$ with $q(1)=2-2n$ and $\oplus$ denotes orthogonal direct sum. In the following chapters we will often denote $L=L_2$.
\end{ex}
\begin{ex}\label{ex:chi_k3n}
The Beauville-Bogomolov form on the second cohomology allows also an easy computation of Euler characteristic of a divisor. In the case $X$ is a manifold of $K3^{[n]}$-type and $D$ a divisor on it we have
\begin{equation}
\chi(D)=\binom{(D,D)_X/2+n+1}{n}.
\end{equation}
See \cite[Example 23.19]{huy} for a proof.
\end{ex}
\begin{ex}
Let $X$ be a \hk manifold of Kummer $n$-type. Then $H^2(X,\mathbb{Z})$ endowed with its Beauville-Bogomolov pairing is isomorphic to the lattice \index{Lattice on a manifold of Kummer $n$-type, $L_{K.n}$}
\begin{equation}\label{latticeGKum}
L_{K.n}=U\oplus U\oplus U\oplus (-2-2n).
\end{equation}
\end{ex}
\begin{ex}
Let $X$ be a \hk manifold deformation equivalent to O'Grady's 6 dimensional example. Then $H^2(X,\mathbb{Z})$ endowed with its Beauville-Bogomolov pairing is isomorphic to the lattice \index{Lattice on $Og_6$, $L_{O.6}$}
\begin{equation}\label{latticeOG6}
L_{O.6}=U\oplus U\oplus U\oplus (-2)\oplus (-2).
\end{equation}

\end{ex}
\begin{ex}
Let $X$ be a \hk manifold deformation equivalent to O'Grady's 10 dimensional example. Then $H^2(X,\mathbb{Z})$ endowed with its Beauville-Bogomolov pairing is isomorphic to the lattice \index{Lattice on $Og_{10}$, $L_{O.10}$}
\begin{equation}\label{latticeOG10}
L_{O.10}=U\oplus U\oplus U\oplus E_8(-1)\oplus E_8(-1)\oplus A_2(-1).
\end{equation}
Here $A_2(-1)$ is a Dynkin lattice defined in \Ref{ex}{an_form}.
\end{ex}
The Beauville-Bogomolov form also allows a useful projectivity criterion:
\begin{prop}\label{prop:proj_dense}
Let $X$ be a \hk manifold. Then $X$ is projective if and only if there exists $v\in H^{1,1}(X,\mathbb{Z})$ such that $q_X(v)>0$.
\end{prop}


Another interesting result relating $H^2(X,\mathbb{C})$ to higher cohomologies has been given by Verbitsky \cite{ver2}:
\begin{thm}
Let $X$ be a \hk manifold of dimension $2n$ and Beauville-Bogomolov form $q_X$. Let $Sym\,H^2(X,\mathbb{C})$ be the subalgebra generated by $H^2(X,\mathbb{C})$. Then
\begin{equation}
Sym\,H^2(X,\mathbb{C})=S^*H^2(X,\mathbb{C})/<\alpha^{n+1}|q_X(\alpha)=0>.
\end{equation}
\end{thm}

\section{Moduli of \hk manifolds and the Torelli problem}
Let us start with two well known deformations of a \hk manifold $X$: the universal deformation $Def(X)$ and the twistor space $TW(X)$ \index{Twistor space, $TW(X),\,TW_\omega(X)$}.
\begin{lem}
Let $X$ be a \hk manifold with \kahl class $\omega$ and symplectic form $\sigma_X$. Then there exists a family
\begin{align} 
TW_{\omega}(X) := &X\times\mathbb{P}^1\\\nonumber
 & \downarrow\\\nonumber
\{(a,b,c)\,\in\,\mathbb{R}^3,\,a^2+b^2+c^2=1\}=S^2 \cong&\mathbb{P}^1 
\end{align}
called Twistor space such that $TW_{\omega}(X)_{(a,b,c)}=X$ with complex structure given by the \kahl class $a\omega+b(\sigma_X+\overline{\sigma}_X)+c(\sigma_X-\overline{\sigma}_X)$.
\end{lem}
A proof that the above defined class gives a complex structure can be found in \cite{hit}
\begin{lem}
Let $X$ be a \hk manifold. Then there exists a flat family $\mathcal{X}\,\rightarrow\,Def(X)$ such that $0\,\in\,Def(X)$, $\mathcal{X}_0\,\cong\,X$. Moreover for every flat family $\mathcal{Y}\,\rightarrow\,S$ such that $\mathcal{Y}_0\cong X$ there exists a commutative diagram
\begin{eqnarray}\nonumber
\mathcal{Y} &\rightarrow & \mathcal{X}\\\nonumber
\downarrow  & &\downarrow  \\ \nonumber
S &\rightarrow & Def(X).
\end{eqnarray}
\end{lem}

\begin{defn}
Let $X$ be a \hk manifold and let \\$H^2(X,\mathbb{Z})\,\cong\,N$. An isometry $f\,:\,H^2(X,\mathbb{Z})\,\rightarrow\,N$ is called a marking of $X$. A couple $(X,f)$ is called a marked \hk manifold.
\end{defn}
We can use this universal deformation and the twistor family to define a moduli space of marked \hk manifolds:
\begin{defn}
Let $(X,\phi)$ be a marked \hk manifold and let $H^2(X,\mathbb{Z})\cong N$. Let $\mathcal{M}_N$ \index{Moduli space of Marked \hk manifolds, $\mathcal{M}_N$} be the set $\{(X,\phi)\}/\sim$ of marked \hk manifolds  where $(X,\phi)\sim (X',\phi')$ if and only if there exists an isomorphism $f\,:X\,\rightarrow\,X'$ such that $f^*=\phi^{-1}\circ \phi'$.
\end{defn}
A priori this definiton endows $\mathcal{M}_N$ only with the structure of a set, but we will use \Ref{thm}{local_torelli} to prove that this is indeed a compact non Hausdorff complex space.
\begin{defn}
Let $X$ be a \hk manifold and let $N$ be a lattice such that $H^2(X,\mathbb{Z})\cong N$. Then we define the period domain $\Omega_N$ \index{Period domain, $\Omega_N$} as
\begin{equation}
\Omega_N=\{x\in \mathbb{P}(N\otimes\mathbb{C})\,|\,(x,x)_N=0,\,(x+\overline{x},x+\overline{x})_N>0\}.
\end{equation}
\end{defn}
\begin{defn}
Let $\mathcal{X}\rightarrow S$ be a flat family of deformations of $X$ and let $f$ be a marking of $X$ into the lattice $N$. Let moreover $F$ be a marking of $\mathcal{X}$ compatible with $f$. Then the period map $\mathcal{P}\,:\,S\,\rightarrow\,\Omega_N$ \index{Period map, $\mathcal{P}$} is defined as follows:
\begin{equation}
\mathcal{P}(s)=F_s(H^{2,0}(\mathcal{X}_s)).
\end{equation}
\end{defn}
Whenever we choose $\mathcal{X}\rightarrow Def(X)$ as flat family of deformations of $X$ we call $\mathcal{P}$ the local period map.
\begin{thm}[Local Torelli, Beauville \cite{beau3}]\label{thm:local_torelli}
Let $(X,f)$ and $N$ be as above, let moreover $F$ be a compatible marking of $\mathcal{X}\rightarrow Def(X)$. Then the map $Def(X)\stackrel{\mathcal{P}}{\rightarrow} \Omega_N$ is a local isomorphism.
\end{thm}
Now this local isomorphism allows to glue the various universal deformations into a complex space.
Another well known fact about the period map is the following:
\begin{thm}[Huybrechts, \cite{huy}]\label{thm:surj_period}
Let $\mathcal{M}_N^0$ be a connected component of $\mathcal{M}_N$. Then the period map $\mathcal{P}\,:\mathcal{M}_N^0\,\rightarrow\,\Omega_N$ is surjective.
\end{thm}
An interesting question is whether we have a global Torelli theorem as in the case of $K3$ surfaces, in general this is false as in the following
\begin{ex}
Let $S$ be a $K3$ surface such that $Pic(S)=\mathbb{Z}\mathcal{O}(C)$ were $C\subset S$ is a smooth rational curve. 
Let $X=S^{[2]}$, it contains $C^{[2]}\cong \mathbb{P}^2$. Let $X'$ be the Mukai flop of $X$ along $C^{[2]}$. Debarre \cite{deb} proved that $X'$ is not isomorphic to $X$ and moreover there exist two markings $f,f'$ such that $\mathcal{P}(X,f)=\mathcal{P}(X',f')$.
\end{ex}
This implies that we cannot hope to have an isomorphism between manifolds with the same period, however the situation is even worse, as the following shows:
\begin{ex}
Let $T$ be a complex torus such that $T^\vee$ is not isomorphic to it and such that $Pic(T)=0$. Let $X=K_2(T)$ and $X'=K_2(T^\vee)$ and let $E$ and $E'$ be respectively the exceptional divisors of  $X\,\rightarrow\,T^{(3)}$ and $X'\,\rightarrow\,(T^\vee)^{(3)}$.\\ It has been shown by Namikawa \cite{nam2} that there exist two markings $f $ and $f'$ such that $\mathcal{P}(X,f)=\mathcal{P}(X',f')$. Suppose there exists a birational map $\psi\,:\,X\,\dashrightarrow\,X'$. This map is regular in codimension two and defines a birational map $E\,\dashrightarrow\,E'$. However the Albanese is a birational invariant hence we obtain $T\,\cong\, Alb(E)\,\cong\,Alb(E')\,\cong\,T^\vee$ which is absurd. 
\end{ex}
However, under some more hypothesis, a weaker Global Torelli theorem holds, see \cite{huy_tor}, \cite{mar1} and \cite{ver}.
\begin{thm}[Global Torelli, Verbitsky, Markman and Huybrechts]\label{thm:global_torelli}
Let $X$ and $Y$ be two \hk manifold of $K3^{[n]}$-type and let $n-1$ be a prime power. Suppose $\psi\,:\,H^2(X,\mathbb{Z})\,\rightarrow\,H^2(Y,\mathbb{Z})$ is an isometry preserving the Hodge structure. Then there exists a birational map $\phi\,:\,X\,\dashrightarrow\,Y$.
\end{thm}
Related to this there is also the following useful theorem, due to Huybrechts \cite[Theorem 3.2]{mar1}:
\begin{thm}\label{thm:graph}
Let $(X,f)$ be a marked Hyperk\"{a}hler manifold and $(X',g)$ be another marked Hyperk\"{a}hler manifold such that $\mathcal{P}(X,f)=\mathcal{P}(X',g)$ and such that the points $(X,f)$ and $(X',g)$ are not separated.  
 Then there exists an effective cycle $\Gamma=Z+\sum_{j}Y_j$ in $X\times X'$ satisfying the following conditions:
\begin{itemize}
\item Z is the graph of a bimeromorphic map from $X$ to $X'$.
\item The codimensions of $\pi_1(Y_j)$ and $\pi_2(Y_j)$ are equal.
\item The composition $g^{-1}\circ f$ is equal to $\Gamma_*\,:\,H^2(X,\mathbb{Z})\,\rightarrow\,H^2(X',\mathbb{Z})$.
\item If $\pi_i(Y_j)$ has codimension 1 then it is supported by an effective uniruled divisor.
\end{itemize}
\end{thm}

\subsection{Moduli spaces of polarized \hk manifolds}\label{ssec:moduli_polar}

In this subsection we will analyze the behaviour of families of \hk manifolds with some conditions imposed on their Picard lattice. Let us start with the following:
\begin{defn}\label{defn:hk_polar}
Let $(X,f)$ be a marked \hk manifold with $H^2(X,\mathbb{Z})\cong N$. Let $h\in N$ be a primitive vector such that $h^2\geq 0$. We call $X$ a $h$-polarized \hk manifold if $f^{-1}(h)$ is represented by an ample divisor on $X$.
\end{defn}

If $X$ is a $h$-polarized manifold for some marking $f$ it is clear that in a projective family $\mathcal{X}$ of \hk manifolds containing $X$ the generic element is still $h$-polarized for a marking $F$ of the family compatible with $f$. Moreover on those elements $\mathcal{X}_t$ such that $F_t^{-1}(h)$ is not ample we still have a closed subset where $F_t^{-1}(h)$ is Nef. Thus we will weaken \Ref{defn}{hk_polar} by imposing only the Nefness of the divisor. The interesting fact is that there is a honest moduli space of such manifolds, see \cite{ghs2}. 
In \Ref{cap}{deformations} we will be interested in a more general case, namely in \hk manifolds such that a given lattice is primitively contained in the Picard lattice:
\begin{defn}\label{defn:hk_gen_polar}
Let $(X,f)$ be a marked \hk manifold with $H^2(X,\mathbb{Z})\cong N$. Let $R\subset N$ be a primitive sublattice of signature $(a,b)$. We call $X$ a $R$-polarized \hk manifold if $f^{-1}(R)\subset Pic(X)$ and, if $a>0$, $f^{-1}(h)$ is represented by a Nef divisor on $X$ for some $h\in R,\,h^2>0$.
\end{defn} 

Also in this case we have a moduli space of such manifolds and we denote it $\mathcal{M}_{R,N}$ or $\mathcal{M}_{R}$ whenever $N$ is understood.

\section{\kahl and positive Cones}

This section analyzes the shape of the Ample (\kahl in the general case) Cones of \hk manifolds, and of other related cones of interest in \Ref{cap}{groups_sporadic}. Recall that the Beauville-Bogomolov form allows to define a notion of positivity on divisors as in the case of surfaces.
\begin{defn}
Let $X$ be a Hyperk\"{a}hler manifold and let $\omega$ be a K\"{a}hler class.
Let $\{l\,\in\,H^{1,1}_{\mathbb{R}}(X),\,\,l^2>0\}$ be the set of positive classes in $H^{1,1}_{\mathbb{R}}(X)$ and let the positive cone $\mathcal{C}_X$ \index{Cone, Positive, $\mathcal{C}_X$} be its connected component containing $\omega$.\\
Let the K\"{a}hler cone $\mathcal{K}_X\,\subset\,\mathcal{C}_X$ \index{Cone, \kahl, $\mathcal{K}_X$} be the set of K\"{a}hler classes.\\
The birational K\"{a}hler cone \index{Cone, Birational \kahl, $\mathcal{BK}_X$} is the union
\begin{equation}
\mathcal{BK}_X\,=\,\bigcup_{f\,:\,X\,\dashrightarrow\,X'}\,f^*\mathcal{K}_{X'},
\end{equation}
where $f\,:\,X\,\dashrightarrow\,X'$ runs through all birational maps $X\,\dashrightarrow\,X'$ from $X$ to another Hyperk\"{a}hler manifold $X'$.
\end{defn}
There are several results on the structure of these cones and there are also some conjectures, see \cite[Section 27 and 28]{huy}, \cite[Section 9]{mar2} and \cite{ht}. Let us summarize most of them:
\begin{prop}
Let $X$ be a \hk manifold. The closure $\overline{\mathcal{K}}_X$ of the \kahl cone is the set of all classes $\alpha\in\overline{\mathcal{C}}_X$ such that $\int_C\alpha\geq0$ for all rational curves $C\subset X$.
\end{prop}
\begin{prop}\label{prop:birat_cone_unir}
Let $X$ be a \hk manifold. The closure $\overline{\mathcal{BK}}_X$ of the birational \kahl cone is the set of all classes $\alpha\in\overline{\mathcal{C}}_X$ such that $q_X(\alpha,D)\geq0$ for all uniruled divisors $D\subset X$.
\end{prop}
The latter is often used together with the following numerical criterion:
\begin{prop}
Let $X$ be as before and let $D\subset X$ be an irreducible effective divisor such that $q_X(D,D)<0$. Then $D$ is uniruled.
\end{prop}
If we specialize to the case of $K3^{[n]}$-type manifolds there are more precise results due to Markman (\cite{mar1} and \cite{mar2}). Let $X$ be a manifold of $K3^{[n]}$-type, he proved that, if $n>3$, the quotient $Q(X):=H^4(X,\mathbb{C})/(S^2H^2(X,\mathbb{C}))$ \index{Mukai Lattice, inside a manifold of \kntipo, $Q(X)$} is an integer Hodge structure of weight 2 and there is a bilinear pairing giving its integer part the structure of a lattice isometric to $U^4\oplus E_8(-1)^2$. Moreover he proved the existence of a unique primitive embedding $i:\,H^2(X,\mathbb{Z})\rightarrow\,Q(X)$. We remark that if $n=2$ or $3$ there exists a unique (up to isometry) primitive embedding $H^2(X,\mathbb{Z})\,\rightarrow\,U^4\oplus E_8(-1)^2$.
 Let $e\in H^{2}(X,\mathbb{Z})$. Let $r$ be the divisibility of $e$ in $H^2(X,\mathbb{Z})$. Let $H^2(X,\mathbb{Z})\subset U^4\oplus E_8(-1)^2$ using either $i$ or the unique embedding. Let $v$ be a generator of $H^2(X,\mathbb{Z})^\perp$ in this embedding. Let $\rho$ be the integer such that $\frac{e+v}{\rho}$ is a primitive class in $U^4\oplus E_8(-1)^2$, let $\sigma$ be the integer such that $\frac{e-v}{\sigma}$ is a primitive class in $U^4\oplus E_8(-1)^2$. We set $rs(e)$ to be the unordered set $\{\rho,\sigma\}$ if $n$ is even and $r=n-1$, otherwise we set it to be $\{\rho/2,\sigma/2\}$.
 
\begin{defn}\label{defn:numerical_exc}
Let $X$ be a manifold of $K3^{[n]}$-type and let $h$ be a \kahl class. A primitive class $e\in H^{2}(X,\mathbb{Z})$ is called numerically exceptional if $(h,e)>0$, $(e,e)=-2$ or $(e,e)=2-2n$ and one of the following holds:
\begin{itemize}
\item $div(e)=2n-2$ and $rs(e)=\{1,n-1\}$.
\item $div(e)=2n-2$, $rs(e)=\{2,(n-1)/2\}$ and $n\equiv 3$ mod $4$.
\item $div(e)=n-1$, $n$ is even and $rs(e)=\{1,n-1\}$.
\item $div(e)=n-1$, $n$ is odd and $rs(e)=\{1,(n-1)/2\}$.
\end{itemize}
We denote $\mathcal{NE}xc_X$ \index{Exceptional classes, Numerical, $\mathcal{NE}xc_X$} the set of numerically exceptional classes
\end{defn}
We can use this to define the following:
\begin{defn}
Let $X$ be as before, we define the fundamental exceptional chamber $\mathcal{FE}xc_X$ \index{Exceptional classes, Fundamental exceptional chamber, $\mathcal{FE}xc_X$} to be the set of $\alpha\in H^2(X,\mathbb{Z})$ such that $(\alpha,e)>0$ for all numerically exceptional class $e$. 
\end{defn}
\begin{thm}\cite[Theorem 1.11 and Proposition 1.5]{mar2}\label{thm:num_eff}
Let $X$ be a \hk manifold of $K3^{[n]}$-type. 
 Let $e\in H^2(X,\mathbb{Z})$ be a numerically exceptional class. Then $ke$ is the class of a reduced irreducible effective divisor, where $k$ is as follows:\newline 
If $e^2=2-2n$ then
\begin{itemize}
\item $k=2$ if $div(e)=2n-2$ and $rs(e)=\{1,n-1\}$.
\item $k=1$ if $div(e)=2n-2$ and $rs(e)=\{2,(n-1)/2\}$.
\item $k=1$ if $div(e)=n-1$.
\end{itemize}
If $e^2=-2$ we have
\begin{itemize}
\item $k=2$ if $div(e)=2$ and $n=2$.
\item $k=1$ if $div(e)=2$ and $n>2$.
\item $k=1$ if $div(e)=1$.
\end{itemize}
\end{thm}
This yields the following result:
\begin{thm}\cite[Prop 5.6]{mar1}\label{thm:birat_cone}
Let $X$ be a manifold of \kntipo, then $\overline{\mathcal{BK}}_X=\overline{\mathcal{FE}xc}_X$.
\end{thm}
In the case of $K3$ surfaces we indeed have $\mathcal{FE}xc=\mathcal{BK}=\mathcal{K}$.\\
Let us specialize further to the case of \ktipo manifolds: \Ref{thm}{birat_cone} implies that $\overline{\mathcal{BK}}_X$ is cut out by $(-2)$ divisors. Let moreover $\mathcal{NK}_X$ \index{Cone, Numerical \kahl, $\mathcal{NK}_X$} be the following cone:
\begin{align}
\mathcal{NK}_X =& \{\alpha\,\in\,\mathcal{C}_X\,|\,(\alpha,e)>0\,\forall\,\text{effective}\,e\in\,Pic(X)\\\nonumber
 & \,\text{s.t.}\,\,\,e^2=-2\,\text{or}\, 
 e^2=-10\,\text{and}\,div(e)=2\}.
\end{align}
Then there is the following conjecture made by Hassett and Tschinkel:
\begin{conj}\cite{ht}\label{conj:ht}
Let $X$ be a fourfold of $K3^{[2]}$-type, then $\mathcal{K}_X=\mathcal{NK}_X$ and moreover for all $e\in Pic(X)$ such that $e^2=-10$ and $div(e)=2$ either $e$ or $-e$ is represented by an effective divisor.
\end{conj}
Some evidence for this conjecture is given by the following:
\begin{oss}
Let $(X,g)$ and $(X',g')$ be two marked projective manifolds of \ktipo and let $f\,:\,X\,\dashrightarrow\,X'$ be a Mukai flop. Then the induced map $g\circ f^*\circ g'^{-1}$ on $L_2$ is the reflection along an element $e$ such that $e^2=-10$ and $div(e)=2$. 
\end{oss}
and by this result:
\begin{thm}[Hassett and Tschinkel, \cite{ht2}]\label{thm:moving_cone}
Let $X$ be a projective manifold of \ktipo and let $w$ be a \kahl class. Let $h$ be the class of a divisor such that $(h,w)>0$ and $(e,h)>0$ for all $e\in\mathcal{NE}xc_X$ and for all $e$ such that $e^2=-10,$ $div(e)=2$ and $(e,w)>0$. Then $h$ is ample.
\end{thm}
There are higher dimensional analogues for this behaviour, see \cite{hht}. 
\begin{ex}
Let $S$ be a $K3$ surface such that $Pic(S)=\mathbb{Z}h$, $h^2=14$. Let $X=S^{[2]}$. We wish to use \Ref{thm}{moving_cone} and \Ref{thm}{birat_cone} to compute the \kahl cone of $X$. Let $2\delta$ be the class of the exceptional divisor on $X$, then $Pic(X)=<h,\delta>$. Let now $C$ be a curve on $S$ in the same numerical class of $h$ and let $C_p$ be a curve in $X$ given by analytical subsets of $S$ consisting in a point of $C$ and a point $p\notin C$. $C_p$ is an effective curve dual to $h$, therefore the positive cone $\mathcal{C}_X$ consists of all elements $ah+b\delta$ with $a>0$ and $b\leq \sqrt{7}a$ or $-b\leq \sqrt{7}a$. Let us remark that all \kahl classes $\omega=ah+b\delta$ must satisfy $b<0$. By \Ref{thm}{birat_cone} the birational \kahl cone $\mathcal{BK}_X$ is cut out by $(-2)$ effective divisors. A direct computation shows that these divisors are the integer solutions of $7a^2-b^2=-1$ with $a\geq 0$. Therefore we have $\delta$ and two series of divisors $\{a_nh +b_n\delta\}_n$ and $\{a_nh-b_n\delta\}_n$ where $a_n/b_n$ tends to $\frac{\sqrt{7}}{7}$ from below.\\ Therefore $\mathcal{BK}_X$ is the set of elements orthogonal to these divisors. It is easy to see that $\mathcal{BK_X}=\mathcal{C}_X\cap \{b<0\}$. Finally we apply \Ref{thm}{moving_cone}: there are no elements of $Pic(X)$ with square $-10$, therefore $\mathcal{K}_X=\mathcal{BK}_X$.  
\end{ex}
\section{Projective families of manifolds of $K3^{[n]}$-type}
In this section we gather several examples of projective \hk manifolds, mainly in the case of manifolds of \ktipo. Notice that most of these examples share the property of being a locally complete family, \ie the image of the local period map has the maximal dimension. In some cases compactifications of these families have been studied, see for example \cite{loo} and \cite{ogr6}.
\subsection{Double EPW Sextics}\label{ssec:epw}
Double EPW sextics were first introduced by O'Grady in \cite{ogr}, they are in many ways a higher dimensional analogous to $K3$ surfaces obtained as the double cover of $\mathbb{P}^2$ ramified along a sextic curve.\\
Let $V\cong\mathbb{C}^6$ be a six dimensional vector space with basis given by $\{e_0,e_1,e_2,e_3,e_4,e_5\}$ and let 
\begin{equation}\nonumber
vol(e_0\wedge e_1\wedge e_2\wedge e_3\wedge e_4\wedge e_5)=1
\end{equation}
be a volume form, giving a symplectic form $\sigma$ on $\Lambda^3V$ defined by
\begin{equation}\nonumber
\sigma(\alpha,\beta)=vol(\alpha\,\wedge\,\beta).
\end{equation}
Let $\mathbb{LG}(\Lambda^3V)$ be the set of lagrangian subspaces of $\Lambda^3V$ with respect to $\sigma$. 
Furthermore let $F$ \index{EPW-sextics, vector bundle $F$} be the vector bundle on $\mathbb{P}(V)$ with fibre 
\begin{equation}\nonumber
F_{v}=\{\alpha\,\in\,\Lambda^3V\,,\,\alpha\,\wedge\,v=0\}.
\end{equation}
Let $A\subset\mathbb{LG}(\Lambda^3V)$ and let $\lambda_A(v)$ \index{EPW-sextics, degeneracy map, $\lambda_A$} be the following composition
\begin{equation}\label{EPW_eq}
F_{v}\,\rightarrow\,\Lambda^3 V\rightarrow (\Lambda^3 V)/A,
\end{equation}
where the first map is injection of $F_{v}$ as a subspace of $\Lambda^3 V$ and the second is the projection to the quotient with respect to $A$.\\
Therefore we define $Y_A[i]$ \index{EPW sextics, degeneracy loci, $Y_A[i]$} as the following locus:
\begin{equation}
Y_A[i]=\{[v]\,\in\,\mathbb{P}(V)\,,\,\dim(A\,\cap\,F_v)\geq i\}.
\end{equation}
Here $Y_A[1]=Y_A$ is the EPW-sextic associated to $A$ and coincides with the degeneracy locus of $\lambda_A$ if $A$ is general.\\
\begin{defn}
Let $\mathbb{LG}(\Lambda^3 V)^0$ \index{EPW-sextics, Open subset of lagrangian subspaces, $\mathbb{LG}(\Lambda^3 V)^0$} be the open subset of lagrangian subspaces $A$ such that the following hold
\begin{itemize}
\item $Y_A[3]=\emptyset$.
\item $Gr(3,6)\cap \mathbb{P}(A)=\emptyset$, where $Gr(3,6)\,\subset\,\mathbb{P}(\Lambda^3 \mathbb{C}^6)$ via the Pl\"{u}cker embedding.
\end{itemize}
\end{defn}
Let us remark that $\mathbb{LG}(\Lambda^3V)^0$ contains the general lagrangian subset.
\begin{thm}\cite[Theorem 1.1]{ogr}\label{thm:kieran_epw}
For $A\in\,\mathbb{LG}(\Lambda^3V)^0$ there exists a double cover $X_A\,\rightarrow\,Y_A$  ramified along $Y_A[2]$ such that $X_A$ is a hyperk\"{a}hler manifold of \ktipo.
\end{thm}
A polarization $h$ of a Double EPW sextic $X_A$ is given by the pullback of the hyperplane section $\mathcal{O}_{Y_A}(1)$ of $Y_A$, a direct computation yields $h^2=2$.
\begin{lem}\label{lem:dual_epw}
Let $A$ be a generic lagrangian subspace and let $Y_A^\vee\,\subset\,\mathbb{P}(V^\vee)$ be the dual hypersurface. then $Y_A^\vee=Y_{\delta(A)}$ \index{EPW-sextics, dual, $Y_{\delta(A)}$} where $\delta$ is the following map
\begin{align}\label{dual_lagr}
\delta\,:\,\mathbb{LG}(\Lambda^3V)  \rightarrow & \mathbb{LG}(\Lambda^3V^\vee)\\\nonumber
A \rightarrow & \{\alpha\,\in\,\Lambda^3V^\vee\,,\,s.t.\,<\alpha,A>=0\}.
\end{align}
Here $<\,,\,>$ is the standard pairing.
\end{lem}
\begin{oss}
It is a well known fact that for $A\,\notin\,\mathbb{LG}(\Lambda^3V)^0$ the situation can be dire indeed: there are degenerate examples where $Y_A=\mathbb{P}(V)$ or when $X_A$ has very bad singularities. 
\end{oss}
Let us look a little into what can happen if the lagrangian $A$ contains some decomposable tensors, first of all we have a result of O'Grady:
\begin{prop}\label{prop:kieran_planes}\cite[Proposition 4.8]{ogr5}
Let $A$ be as before and let $I$ be the set of decomposable tensors contained in $A$. If $I$ is finite then its cardinality is at most 20.
\end{prop}
If on the other hand we have an infinite set of decomposable tensors inside $A$ we obtain an infinite set of planes intersecting each other in a point and we have the following classical result of Morin \cite{mori}.
\begin{thm}\label{thm:morin_planes}
Let $\mathcal{W}$ be an infinite complete (\ie contains all planes meeting all elements of $\mathcal{W}$ in a point) set of planes in $\mathbb{P}^5$ meeting each other in a point. Then $\mathcal{W}$ satisfies one of the following
\begin{itemize}
\item There exists $v\in\mathbb{P}^5$ s.t. $v\in\cap_{W\in\mathcal{W}}W$.
\item There exists a plane $W'$ such that $W\cap W'$ is a line for all $W\in\mathcal{W}$.
\item The set of vectors contained inside some element of $\mathcal{W}$ spans a hyperplane.
\item All elements of $\mathcal{W}$ are the planes contained in a quadric $Q$.
\item All elements of $\mathcal{W}$ are tangent to a Veronese surface.
\item All elements of $\mathcal{W}$ meet a Veronese surface in a conic.
\end{itemize}
\end{thm}

Since we will be interested in automorphisms of EPW-sextics we will need the following:
\begin{prop}\label{prop:EPW_autom}
Let $G\subset PGL(6,\mathbb{C})$ be a simple group of automorphisms of $\mathbb{P}^5$ leaving a lagrangian subspace $A\in\mathbb{LG}(\Lambda^3\mathbb{C}^6)^0$ invariant. Let $X_A\,\rightarrow Y_A$ be the \hk cover of $Y_A$. Then $G$ extends to a group of automorphisms of $Y_A$ and of $X_A$.  Suppose that $G$ acts faithfully on $Y_A$ and trivially on a section of $K_{Y_A}$. Then $G$ acts on $X_A$ as a group of symplectic automorphisms acting trivially on its natural polarization.
\begin{proof}
Let $y\in Y_A[i]$, since $G$ preserves $A$ we have $g(y)\in Y_A[i]$ for all $g\in G$. Therefore $G$ induces automorphisms on $Y_A$ and, since it preserves also $Y_A[i]$, it extends also to its cover $X_A\rightarrow Y_A$ ramified along $Y_A[2]$. Notice that we obtain a (usually nontrivial) extension of $G$ with the covering involution $\tau_A$ of $X_A$. If we suppose moreover that $G$ acts trivially on a section of the Canonical divisor of $Y_A$ we have that all its elements act trivially also on sections of $K_{X_A}$, such as $\sigma_{X_A}^2$. Therefore $g\in G$ acts as $\pm Id$ on the symplectic form $\sigma_{X_A}$. If $g(\sigma_{X_A})=-\sigma_{X_A}$ we can use $\tau_A$ to obtain nonetheless a symplectic automorphism $g\tau_A$ of $X_A$.  
\end{proof}
\end{prop}
\subsection{Fano schemes of lines on cubic fourfolds}
Fano schemes of lines on cubic fourfolds were first studied By Beauville and Donagi \cite{bedon}, where the authors proved that they are \hk manifolds of \ktipo.\\
Let $X\,\subset\,\mathbb{P}^5$ be a smooth cubic fourfold and let $F(X)$ \index{Fano scheme of lines, $F(X)$} be the scheme parametrizing lines contained in $X$.
\begin{thm}\label{thm:fano_map}
Keep notation as above, then the following hold:
\begin{itemize}
\item $F(X)$ is a \hk manifold.
\item $F(X)$ is deformation equivalent to $K3^{[2]}$.
\item the Abel-Jacobi map
\begin{equation}\label{abel} \alpha\,:\,H^4(X,\mathbb{C})\,\rightarrow\,H^2(F(X),\mathbb{C})
\end{equation}
is an isomorphism of rational Hodge structures.
\end{itemize}
\end{thm} 
 
Let us remark that the proof of \Ref{thm}{fano_map} gives also a polarization $h$ of $F(X)$ which is the restriction of the Pl\"{u}cker polarization of $Gr(2,6)$ to $F(X)$. Moreover we have $h^2=6$ and $div(h)=2$.
\begin{oss}\label{oss:residue_fano}
Since $X$ is a hypersurface it is possible to give generators of its cohomology in terms of the residues of its defining equation $f$, as proved classically by Griffiths \cite{grif}. By \cite[Théorème 18.1]{voi} we have that $H^{3,1}(X)=\mathbb{C}\,\frac{Res(\Omega)}{f^2}$, where $\Omega=\sum (-1)^ix_i x_0\wedge \dots \widehat{x_i}\dots \wedge x_5$. This formula is particularly useful in determining whether an automorphism induced on $F(X)$ by one on $X$ is symplectic or not. 
\end{oss}
\begin{oss}\label{oss:cover_fano}
Let $f=x_0^3+f'$ be a nonsingular cubic polynomial, where $f'$ is a polynomial in $x_1,\dots,x_5$ and let $Y=V(f)$ be a cubic fourfold. The natural projection of $Y$ from $e_0$ to the hyperplane $e_0^\vee$ is a $3:1$ cover ramified along a cubic threefold. Obviously the covering morphism induces an order 3 morphism on $F(Y)$. A direct computation with \Ref{oss}{residue_fano} shows that this automorphism is nonsymplectic. 
\end{oss}
\subsection{Moduli spaces of sheaves on $K3$ surfaces}
Some very interesting examples of projective \hk manifold are given by moduli spaces of sheaves on polarized K3 surfaces. These examples have been studied by several people, we will refer to \cite{hl_moduli} for a complete list of references. First of all we have to define what are the Mukai vector of a sheaf and the Mukai pairing.
\begin{defn}
Let $X$ be a smooth manifold and let $E,F$ be two coherent sheaves. Then $v(E)=ch(E)\sqrt{td(X)}$ \index{Mukai vector, $v(E)$} is the Mukai vector of $E$ and $(v(E),v(F))_M:=-\chi(E,F)$ \index{Mukai pairing, $(v(E),v(F))_M$} is their Mukai pairing.
\end{defn}
We will denote $\mathcal{M}_v(S,H)$ \index{Moduli space of sheaves, $\mathcal{M}_v(S,H)$} to be the moduli spaces of stable sheaves on a K3 surface $S$ with Mukai vector $v$ with respect to the polarization $H$ of $S$. To state the fundamental result we will need the concept of $v$-generic polarization, we will not state this condition precisely but it is sufficient to know that this condition is indeed generic.
\begin{thm}
Let $S$ be a $K3$ surface, let $v\in H^{*}(S,\mathbb{Z})$ be a primitive Mukai vector such that $(v,v)_M\geq0$ and $rk(v)>0$. Let $H$ be a $v$-generic polarization. Then $\mathcal{M}_v(S,H)$ is a \hk manifold of dimension $2+(v,v)_M$.
\end{thm}
In the realm of \hk manifolds the most interesting case is when $(v,v)_M=n\geq 2$, in this case $\mathcal{M}_v(S,H)$ is a \hk manifold of $K3^{[n]}$-type and the Hodge structure of its second cohomology is given by the weight 2 Hodge structure on $v^\perp\subset H^{*}(S,\mathbb{Z})$ with pairing given by the Mukai pairing.\\
\subsection{Varieties of sums of powers}
In \cite{ir1} and \cite{ir2} Iliev and Ranestad introduce another maximal projective family of manifolds of \ktipo, namely the variety of sums of powers of a cubic fourfold.
\begin{defn}
Let $f$ be a homogeneous polynomial of degree $d$ in $n+1$ variables, defining the hypersurface $X\subset\mathbb{P}^n$. Let $VSP(f,s)$ \index{Variety of sums of powers, $VSP(f,s)$} be the closure of
\begin{equation}
\{\{<l_1>,\dots,<l_s>\}\in\,Hilb_s((\mathbb{P}^n)^\vee)\,,\,\exists\,\lambda_i\in\mathbb{C}\,:\,f=\lambda_1l_1^d+\dots+\lambda_sl_s^d\}.
\end{equation}
\end{defn}

\begin{thm}\cite{ir1}
Let $f$ be a general cubic polynomial in 6 variables, then $VSP(f,10)$ is a \hk manifold of \ktipo with a polarization given by an embedding $VSP(f,10)\rightarrow\, Gr(4,\Lambda^2 \mathbb{C}^6)$. 
\end{thm}
Iliev and Ranestad also analyze the natural correspondence between the Fano scheme of lines of the cubic hypersurface $X=V(f)$ and $VSP(f,10)$ and prove that the two families are distinct:
\begin{lem}\cite{ir1}
Let $f=f(x_0,x_1,\dots,x_5)$ be a general cubic polynomial, and let $F=F(V(f))$ be the associated fano scheme of lines. Let $VSP(f,10)$ be the fourfold given by the variety of sums of powers. Then the two families 
\begin{eqnarray*}
\mathcal{F}(U) &  & \mathcal{VSP}(U,10)\\
\searrow & & \swarrow\\
 &U=\{\text{general cubic polynomials}\}.&
\end{eqnarray*}
Intersect transversally along a locus given by Hilbert schemes of two points on a $K3$ surface of degree $14$.
\end{lem}
The following was first computed by A. Iliev, K. Ranestad and B. Van Geemen:
\begin{prop}\cite{vgir}\label{prop:vgir_mail}
Keep notation as above, then the natural polarization on $VSP(f,10)\subset Gr(4,\Lambda^2 \mathbb{C}^6)$ has square $38$ and divisibility $2$.
\begin{proof}
Without loss of generality we can compute everything in the codimension 1 locus of Hilbert schemes of 2 points on a $K3$ surface of degree 14.
Let $S$ be a generic $K3$ surface of degree $14$ such that $X=S^{[2]}\cong F(V(f))\cong VSP(f',10)$ for some cubic polynomials $f$ and $f'$. Since $S$ is generic we know $Pic(S)=\mathbb{Z}h$, where $h^2=14.$ Moreover this implies $Pic(S^{[2]})=<h,\delta>$ where $2\delta$ is the class of the exceptional fibre, $(h,\delta)=0$, $(\delta,\delta)=-2$ and $\delta$ has divisibility 2. We also have the polarizations induced by those on the Fano scheme of lines and on the variety of sum of powers,  namely two classes $l_{Fano}$ and $l_{vsp}$. It is known that $l_{Fano}=2h-5\delta$ by \cite{bedon}, it has square 6 and divisibility $2$.
Moreover let $C$ be a divisor of $S$ representing the polarization $h$, let $p\in S$ and let $C_p\subset X$ be a curve given by length 2 subschemes of $S$ containing $p$ and a point of $C$. Let $D_p$ be the rational curve parametrizing length 2 subschemes of $S$ supported on $p$. Notice that $h$ is represented by $C_S$, \ie subschemes of length 2 of $S$ supported on a point of $S$ and one of $C$. Moreover $2\delta$ is represented by $\{D_p\}_{p\in S}$.\\
$C_p$ and $D_p$ are nothing else than a basis of $H^{3,3}(X,\mathbb{Q})$ dual to $h$ and $\delta$ in the following sense:
\begin{align}\label{hcp_vsp}
(h,C_p) =& 14,\\
\label{hdp_vsp} (h,D_p)=&0,\\
\label{dcp_vsp} (\delta,C_p)=& 0, \\
\label{ddp_vsp} (\delta,D_p)=& -1.
\end{align}
Here \eqref{hcp_vsp} is obtained by setting $h=C'_S,$ for $C'\subset S$ a curve linearly equivalent to $C$, \eqref{hdp_vsp} and \eqref{dcp_vsp} are a consequence of $(h,\delta)=0$ and \eqref{ddp_vsp} is obtained by linearity from the fact that $(l_{Fano},D_p)=5$ (see \cite{ir2}).
In the same paper it is proved that $(l_{vsp},D_p)=3$, therefore $l_{vsp}=ah-3\delta$. Now we only need to evaluate $a$.
First of all $a\geq 2$ otherwise $l_{vsp}$ would have negative square. The polarization $l_{vsp}$ is the restriction of the Pl\"{u}cker embedding $Gr(4,\Lambda^2 \mathbb{C}^6)\subset\mathbb{P}(\Lambda^4\mathbb{C}^6)$ (see \cite[Proof of Lemma 3.6]{ir1}) and the dimension of the Pl\"{u}cker embedding of this grassmanian is $1364$, therefore $h^0(l_{vsp})\leq 1365$. It is possible to compute $\chi(l_{vsp})$ in terms of its Beauville-Bogomolov form (see \Ref{ex}{chi_k3n}). Moreover by Kodaira's vanishing one can conclude that $h^0(l_{vsp})=\chi(l_{vsp})\leq 1365.$ However $\chi(l_{vsp})=\binom{(l_{vsp})^2/2+3}{2}$, which is greater than $1365$ as soon as $a\geq 3$. Therefore $l_{vsp}=2h-3\delta$, it has square 38 and divisibility 2.   
\end{proof}

\end{prop}
We must remark that the hypothesis on the generality of $f$ is indeed necessary, as the following example with a nonsingular polynomial shows:
\begin{oss}
Let $f=x_0^3+x_1^3+f'$, where $f'$ is a general cubic polynomial on $\{x_1,\dots,x_5\}$. This is a 10-dimensional subset of cubic polynomials and we have the following inclusion:
\begin{equation}
VSP(f',8)\subset VSP(f,10).
\end{equation}
Here $VSP(f',8)$ is obtained by points of the form $(x_0,x_1,l_1,\dots,l_8)$ inside $Hilb_2((\mathbb{P}^2)^\vee)\,\times\,Hilb_8((\mathbb{P}^4)^\vee)$.
However it was proven in \cite{ras} that $VSP(f',8)$ has dimension 5, therefore $VSP(f,10)$ has dimension greater than $4$.
\end{oss} 

\subsection{Subspaces of Grassmannians}
This last example was introduced by Debarre and Voisin \cite{dv} and deals with a certain subspace of a Grassmannian. Let $V$ be a 10 dimensional vector space and let $\sigma\in\Lambda^3V^\vee$ be a generic 3-form on $V$. Let moreover $Y_\sigma\subset G(6,V)$ be the set of six dimensional subspaces where $\sigma$ vanishes identically and let $F_{\sigma}\subset G(3,V)$ be the set of 3 dimensional subspaces where $\sigma$ vanishes. Notice that $F_\sigma$ is a hypersurface in $G(3,V)$
\begin{thm}\cite{dv}
Let $V,\,\sigma,\,F_\sigma $ and $Y_\sigma$ be as before. Then $Y_\sigma$ is a \hk fourfold of \ktipo. Moreover there is an isomorphism of weight 2 rational Hodge structures
\begin{equation}
H^{20}(F_\sigma,\mathbb{C})_{van}\cong H^2(Y_\sigma,\mathbb{C})_{van}.
\end{equation}
\end{thm}
Debarre and Voisin also prove that $Y_\sigma$ has a polarization of square 22.
\setcounter{prop}{0}

\chapter{Lattice theory}\label{cap:lattice}
This chapter is devoted to gather all necessary results about lattices and quadratic forms, in particular we make extensive use of discriminant forms and groups. The interested reader can consult \cite{nik2} for what concerns discriminant forms, \cite{con} for what concerns most of the lattices treated in this chapter and also \cite{atlas} for some information on the groups we treat often here.
\section{Discriminant forms and applications}
First of all let us start with the basic notions of discriminant groups and forms: given an even lattice $N$ with quadratic form $q$ we can consider the group $A_N=N^{\vee}/N$ \index{Lattice, Discriminant group, $A_N$} which is called discriminant group and whose elements are denoted $[x]$ for $x\in N^\vee$. We denote with $l(A_N)$ \index{Length of a group $G$, $l(G)$} the least number of generators of $A_N$. On $A_N$ there is a well defined quadratic form $q_{A_N}$ \index{Lattice, Discriminant form on $A_N$, $q_{A_N}$} taking values inside $\mathbb{Q}/2\mathbb{Z}$ which is called discriminant form; moreover we call $(n_+,n_-)$ \index{Lattice, Signature, $(n_+,n_-)$} the signature of $q$ and therefore of $N$ as a lattice.
It is possible to define the signature $sign(q)$ \index{Lattice, Signature of a Discriminant form, $sign(q)$} of a discriminant form $q$ (modulo 8) as the signature modulo 8 of a lattice having that discriminant form. This notion is well defined since 2 lattices $N,N'$ such that $q_{A_N}=q_{A_{N'}}$ are stably equivalent, \ie there exist 2 unimodular lattices $T,T'$ such that $N\oplus T\cong N'\oplus T'$.\\
One more definition we will need is that of the genus of a lattice: two lattices $N$ and $N'$ are said to have the same genus if $N\otimes\mathbb{Z}_p\cong N'\otimes\mathbb{Z}_p$ for all primes $p$. Notice that there might be several isometry classes in the same genus.\\


\begin{lem}\cite[Corollary 1.13.5]{nik2}\label{lem:nik_spezza}
Let $S$ be an even lattice of signature $(t_+,t_-)$. Then the following hold:
\begin{itemize}
\item If $t_+>0$, $t_->0$ and $t_++t_->2+l(A_S)$ then $S\cong U\oplus T$ for some lattice $T$.
\item If $t_+>0$, $t_->7$ and $t_++t_->8+l(A_S)$ then $S\cong E_8(-1)\oplus T$ for some lattice $T$.
\end{itemize}
\end{lem}

\begin{lem}\cite[Proposition 1.4.1]{nik2}\label{lem:nik_overlattice}
Let $S$ be an even lattice. There exists a bijection $S'\rightarrow H_{S'}$ between even overlattices of finite index of $S$ and isotropic subfactors of $A_S$, moreover the following hold:
\begin{enumerate}
\item $A_{S'}=(H_{S'}^\perp)/H_{S'}\subset A_S$.
\item $q_{A_{S'}}=q_{A_S|A_{S'}}$.
\end{enumerate}
\begin{proof}
Let us briefly give an idea of the proof: suppose $v\in S'-S$, then it defines an element of $A_S$ and its square is $v^2$ modulo $2\mathbb{Z}$, \ie it is $0$. Let $H_{S'}$ be the image of all elements in $S'-S$. Clearly in the natural inclusion $A_{S'}\subset A_S$ all elements of $A_{S'}$ are orthogonal to $H_S$ and the intersection is 0. Conversely, if we have an isotropic subgroup $H_S'\subset A_S$ where its elements are of the form $v/n$, $v\in S$ and $v^2$ is a multiple of $2n^2$ we define an overlattice $S'$ by adding the vectors $v/n$.   
\end{proof}
\end{lem}
\begin{oss}\label{oss:nik_overlattice2}
\Ref{lem}{nik_overlattice} is particularly useful in the particular case of overlattices $R$ of $T\oplus S$ where both $S$ and $T$ are primitive in $R$. In this case our isotropic subgroup $H_{T\oplus S}$ is of the form $(a,\phi(a))$, where $a\in B_T\subset A_T$ is not isotropic and $\phi$ is an isometry between $B_T$ and its image in $A_S(-1)$.
\end{oss}
We will often need to analyze primitive embeddings of an even lattice into another one, let us make some useful remarks whose proofs can also be found in \cite{nik2}:
\begin{oss}\label{oss:pre_prim}
A primitive embedding of an even lattice $S$ into an even lattice $N$ is equivalent to giving $N$ as an overlattice of $S\oplus S^{\perp_N}$ corresponding to an isotropic subgroup $H_S$ of $A_S\oplus A_{S^{\perp_N}}$. Moreover there exists an isometry $\gamma\,:\,p_S(H_S)\,\rightarrow\,p_{S^{\perp_N}}(H_S)$ between $q_S$ and $q_{S^{\perp_N}}$ ($p_S$ denotes the natural projection $A_S\oplus A_{S^{\perp_N}}\rightarrow A_S$). Note moreover that this implies $H_S=\Gamma_\gamma(p_S(H_S))$ where $\Gamma_\gamma$ is the pushout of $\gamma$ in $A_S\oplus A_{S^{\perp_N}}$.
\end{oss}
\begin{oss}
Suppose we have a lattice $S$ with signature $(s_+,s_-)$ and discriminant form $q(A_S)$ primitively embedded into a lattice $N$ with signature $(n_+,n_-)$ and discriminant form $q(A_N)$ and let $K$ be a lattice, unique in its genus and such that $O(K)\rightarrow O(q_{A_K})$ is surjective, with signature $(k_+,k_-)$ and discriminant form $-q(A_N)$.\\ It follows from \cite{nik2} that primitive embeddings of $S$ into $N$ are equivalent to primitive embeddings of $S\oplus K$ into an unimodular lattice $T$ of signature $(n_++k_+,n_-+k_-)$ such that both $S$ and $K$ are primitively embedded in $T$. By \Ref{oss}{pre_prim} an embedding of $S\oplus K$ into a finite overlattice $V$ such that both $S$ and $K$ are primitively embedded into it is equivalent to giving subgroups $H_S$ of $A_S$ and $H_N$ of $A_N$ and an isometry $\gamma\,:\,q_{A_S|H_S}\,\rightarrow\,-q_{A_N|H_N}$. Finally a primitive embedding of $V$ into $T$ is given by the existence of a lattice with signature $(v_-,v_+)$ and discriminant form $-q_V$. 
\end{oss}
Keeping the same notation as before we give a converse to these remarks:
\begin{lem}\cite[Proposition 1.15.1]{nik2}\label{lem:nik_immerge}
Primitive embeddings of S into an even lattice $N$ are determined by the sets $(H_S,H_N,\gamma,K,\gamma_K)$ \index{Lattice, Quintuple defining a primitive embedding, $(H_S,H_N,\gamma,K,\gamma_K)$} where $K$ is an even lattice with signature  $(n_+-s_+,n_--s_-)$ and discriminant form $-\delta$ where $\delta\,\cong\,(q_{A_S}\oplus -q_{A_N})_{|\Gamma_\gamma^\perp/\Gamma_\gamma}$ 
and $\gamma_K\,:\,q_K\,\rightarrow\,(-\delta)$ is an isometry.\\
Moreover two such sets $(H_S,H_N,\gamma,K,\gamma_K)$ and $(H'_S,H'_N,\gamma',K',\gamma'_K)$ determine isomorphic sublattices if and only if
\begin{itemize}
\item $H_S=\lambda H'_S$, $\lambda\in O(q_S)$,
\item $\exists\,\epsilon\,\in\,O(q_{A_N})$ and $\psi\,\in\,Isom(K,K')$ such that $\gamma'=\epsilon\circ\gamma$ and $\overline{\epsilon}\circ\gamma_K=\gamma'_K\circ\overline{\psi}$, where $\overline{\epsilon}$ and $\overline{\psi}$ are the isometries induced among discriminant groups.
\end{itemize} 
\end{lem}
For many purposes we will use only the following simplified version of \Ref{lem}{nik_immerge}:
\begin{lem}\label{lem:nik_immerge1}
Let $S$ be an even lattice of signature $(s_+,s_-)$. The existence of a primitive embedding of $S$ into some unimodular lattice $L$ of signature $(l_+,l_-)$ is equivalent to the existence of a lattice $M$ of signature $(m_+,m_-)$ and discriminant form $q_{A_M}$ such that the following are satisfied: 
\begin{itemize}
\item $s_++m_+=l_+$ and $s_-+m_-=l_-$.
\item $A_M\cong A_S$ and $q_{A_M}=-q_{A_S}$.
\end{itemize}
\end{lem}
We will also use a result on the existence of lattices, the following is a simplified version of \cite[Theorem 1.10.1]{nik2}
\begin{lem}\label{lem:nik_esiste}
Suppose the following are satisfied:
\begin{itemize}
\item $sign(q_T)\equiv t_+-t_-$ $mod\,8$.
\item $t_+\geq0$, $t_-\geq0$ and $t_++t_-\geq l(A_T)$.
\item There exists a lattice $T'$ of rank $t_++t_-$ and discriminant form $q_T$ over the group $A_T$.
\end{itemize}
Then there exists an even lattice $T$ of signature $(t_+,t_-)$, discriminant group $A_T$ and form $q_{A_T}$.
\end{lem}
\begin{oss}\label{oss:overl_group}
Let $M$ and $M'$ be lattices and let $N$ be an overlattice of $M\,\oplus\,M'$. Then $l(A_N)\leq\,l(A_M)+l(A_{M'})$.
\end{oss}

Let us give a few examples on the computation of discriminant forms and groups:
\begin{ex}\label{ex:K32_lat}
Let $L$ be as in \eqref{latticeK3n} for $n=2$, then
\begin{equation}
A_{L}=\mathbb{Z}_{/(2)},\,\,\,\,\,\,\,\,q_{A_L}(1)=-\frac{1}{2}.
\end{equation}
\end{ex}
\begin{ex}\label{ex:an_form}
Let $A_n$ \index{Lattice, Dynkin of type $A_n$, $A_n$} be the Dynkin lattice given by $\{\,v=\sum a_ie_i\in\mathbb{Z}^{n+1},\,\sum a_i=0\,\}$ with the bilinear form induced by the euclidean bilinear form, then it has discriminant group $\mathbb{Z}_{/(n+1)}$ generated by an element of the form $(\frac{1}{n+1},\dots,\frac{1}{n+1},-\frac{n}{n+1})$.
\end{ex}
\begin{ex}
Let $\{(x_1,\dots,x_n)\in\mathbb{Z}^{n},\,\,\sum x_i\,\text{is even}\}$ be the positive definite Dynkin lattice of type $D_n$ \index{Lattice, Dynkin of type $D_n$, $D_n$}. Then its discriminant group is $\mathbb{Z}_{/(2)}^2$ if $n$ is even and $\mathbb{Z}_{4}$ otherwise. In any case its 4 elements are the modulo $D_n$ classes of $(0,\dots,0),(\frac{1}{2},\dots,\frac{1}{2}),(0,\dots,0,1)$ and $(\frac{1}{2},\dots,\frac{1}{2},-\frac{1}{2})$.
\end{ex}
\begin{ex}
Let $n=4k$ and let $D_n^+\subset\mathbb{Z}^n$ \index{Lattice, $D_n^+$} be the lattice generated by $D_n$ and $(\frac{1}{2},\dots,\frac{1}{2})$. It is an unimodular lattice and it is even if $k$ is even. Moreover $D_8^+$ is usually defined as the Dynkin lattice $E_8$ \index{Lattice, Dynkin of type $E_8$, $E_8$}. If $k$ is even this gives an easy example to \Ref{lem}{nik_overlattice} where the isotropic subgroup of $A_{D_n}$ is generated by the class of $(\frac{1}{2},\dots,\frac{1}{2})$.
\end{ex}
\begin{ex}\label{ex:e67_form}
Let $v,w\in E_8$ be two elements of square $2$ such that $<v,w>\cong A_2$. Then $v^{\perp}=E_7$ \index{Lattice, Dynkin of type $E_7$, $E_7$} and $<v,w>^{\perp}=E_6$ \index{Lattice, Dynkin of type $E_6$, $E_6$}. By \Ref{lem}{nik_immerge1} we have $A_{E_6}=A_{A_2}$ and $q_{E_6}=-q_{A_2}$. Analogously $A_{E_7}=A_{A_1}$ and $q_{E_7}=-q_{A_1}$.
\end{ex}
\begin{ex}\label{ex:discrinvol}
The lattice $E_8(-2)$ has discriminant group $(\mathbb{Z}_{/(2)})^8$ and discriminant form $q_{E_8(-2)}$ given by the following matrix:
\begin{equation}\nonumber
\left(\begin{array}{cccccccc}
1&0&0&\frac{1}{2}&0&0&0&0\\
0&1&\frac{1}{2}&0&0&0&0&0\\
0&\frac{1}{2}&1&\frac{1}{2}&0&0&0&0\\
\frac{1}{2}&0&\frac{1}{2}&1&\frac{1}{2}&0&0&0\\
0&0&0&\frac{1}{2}&1&\frac{1}{2}&0&0\\
0&0&0&0&\frac{1}{2}&1&\frac{1}{2}&0\\
0&0&0&0&0&\frac{1}{2}&1&\frac{1}{2}\\
0&0&0&0&0&0&\frac{1}{2}&1
\end{array}\right).
\end{equation}
\end{ex}
\begin{ex}\label{ex:latticeLam}
Let
\begin{align}\label{latticeLam} \index{Mukai lattice, abstract, $L'$}
L' &=U^4\oplus E_8(-1)^2,\\\label{latticeM} 
M_2 &=E_8(-2) \oplus U^3\oplus\,(-2). \index{Lattice, $M_2$}
\end{align}
Since $L'$ is unimodular $A_{L'}=\{0\}$.\\
The lattice $(-2)$ has discriminant group $\mathbb{Z}_{/(2)}$ and discriminant form $q'$ with $q'(1)=q_{A_{1}(-1)}(1)=-\frac{1}{2}$ as in \Ref{ex}{K32_lat}.
Therefore the lattice $M_2$ has discriminant form $q_{E_8(-2)}\oplus q'$ over the group $(\mathbb{Z}_{/(2)})^9$.
\end{ex}
We wish to remark that often $L'$ is called the Mukai lattice because it is isometric to the lattice given by $H^*(S,\mathbb{Z})$, where $S$ is a $K3$ surface and the pairing is the Mukai pairing.
\begin{ex}\label{ex:un_form}
Let $R=U(n)$. Then its discriminant group is $\mathbb{Z}_{/(n)}^2$ with discriminant form $\left(\begin{array}{cc} 0 & \frac{1}{n}\\ \frac{1}{n}& 0\end{array}\right)$.
\end{ex}
\begin{ex}\label{ex:discr_autom3}
Let $M_3=U\oplus U(3)^2\oplus A_2(-1)^2\oplus (-2)$. \index{Lattice, $M_3$}
Then $A_{M_3}=\mathbb{Z}_{/(3)}^6\times\mathbb{Z}_{/(2)}$ and its discriminant form is obtained as the direct sum of those of its addends as detailed in \Ref{ex}{un_form}, \Ref{ex}{K32_lat} and \Ref{ex}{an_form} (with the appropriate sign changes). 
\end{ex}
\begin{ex}\label{ex:discr_autom5}
Let $M_5=U\oplus U(5)^2\oplus (-2)$. \index{Lattice, $M_5$} Then $A_{M_5}=\mathbb{Z}_{/(5)}^4\times\mathbb{Z}_{/(2)}$ with discriminant form obtained from \Ref{ex}{un_form} and \Ref{ex}{K32_lat}.
\end{ex}
\begin{ex}\label{ex:discr_autom7}
Let $M_7=U(7)\oplus \left(\begin{array}{cc} 4 & 1\\ 1&2 \end{array}\right)\oplus (-2)\cong U\oplus U(7)\oplus (14)$. \index{Lattice, $M_7$} Then $M_7$ has discriminant group $\mathbb{Z}_{/(7)}^2\times\mathbb{Z}_{/(14)}$ with discriminant form 
\begin{equation}
\left(\begin{array}{ccc} 
0 & \frac{1}{7}& 0\\
\frac{1}{7} & 0 & 0\\
0 & 0& \frac{1}{14}
\end{array}\right).
\end{equation}
\end{ex}

To conclude this section we analyze the behaviour of (-2) vectors inside $L$, $M_2,M_3$ and $M_7$, since they will play a fundamental role in \Ref{lem}{denselemma}. Hence we will need the following:
\begin{lem}\label{lem:evenembed}
Let $(-2)$ be $A_1(-1)$ and let $e$ be one of its generators. Let $L$,$M_2,M_3$ and $M_5$ be as before.
Then the following hold:
\begin{itemize}
\item Up to isometry there is only one primitive embedding $(-2)\hookrightarrow M_2$ such that $(e,M_2)=2\mathbb{Z}$ (\ie $e$ is 2-divisible). Moreover $e\oplus e^\perp=M_2$.
\item Up to isometry there is only one primitive embedding $(-2)\hookrightarrow L$ such that $(e,L)=2\mathbb{Z}$. Moreover $e\oplus e^\perp=L$.
\item Up to isometry there is only one primitive embedding $(-2)\hookrightarrow M_3$ such that $(e,M_3)=2\mathbb{Z}$. Moreover $e\oplus e^\perp=M_3$.
\item Up to isometry there is only one primitive embedding $(-2)\hookrightarrow M_5$ such that $(e,M_5)=2\mathbb{Z}$. Moreover $e\oplus e^\perp=M_5$.

\end{itemize}

%
Furthermore all other primitive embeddings into $M_2$ given by\\ $(H_e,H_{M_2},\gamma,K,\gamma_K)$ satisfy the following:
\begin{equation}\label{halfness_eq}
\exists s\,\in\, A_K,\,\,q_{A_K}(s,s)=\pm\frac{1}{2}.
\end{equation}

\begin{proof}



By \Ref{lem}{nik_immerge} we know that the quintuple \\$(H_e,H_{M_2},\gamma,K,\gamma_K)$ determines primitive embeddings of $e$ inside $M_2$ and the quintuple $(H_e,H_L,\gamma,K,\gamma_K)$ provides those into $L$.\\
A direct computation shows that primitive embeddings of $e$ into $L$ are 2-divisible only for the quintuple $(\mathbb{Z}_{/(2)},A_L,Id,U^3\oplus E_8(-1)^2,Id)$.\\ Now let us move on to the case of $M_2$:\\
 If $H_e=Id$ then we have $K\cong U^2\oplus E_8(-2)\oplus (2)\oplus (-2)$, obviously $e$ is not 2-divisible in this case and this satisfies \eqref{halfness_eq}. If $H_e=\mathbb{Z}_{/(2)}$ and $(H_{M_2},A_{M_2}^{\perp_{A_{E_8}}})\neq 0$ we obtain nonetheless condition \eqref{halfness_eq} and again $e$ is not 2-divisible in this embedding since $e\oplus e^{\perp_{M_2}}$ is properly contained in $M_2$ with index a multiple of 2. Therefore $(\mathbb{Z}_{/(2)},A_{M_2}^{\perp_{A_{E_8}}},Id,U^3\oplus E_8(-2),Id)$ is the only possible case. The proof goes the same for $M_3$ and $M_5$.

\end{proof}
\end{lem}

\section{Lattices over cyclotomic fields}\label{sec:cyclotomic_lattices}
In this section we define a special class of lattices, namely lattices defined over rings different from $\mathbb{Z}$. Most of the lattices of \Ref{cap}{groups_sporadic} can be better understood in this context. Throughout this section all integer lattices will be \emph{definite}, either positive or negative.
\begin{defn}
Let $\omega_n$ be an $n$-th primitive root of unity and let $\mathbb{D}_n:=\mathbb{Z}[\omega_n]$ \index{Ring of cyclotomic integers, $\mathbb{D}_n$} be the ring of cyclotomic integers. A free  $\mathbb{D}_n$-module $C$ is a $\mathbb{D}_n$-lattice if it is endowed with a nondegenerate hermitian pairing \index{Lattice, Nondegenerate hermitian pairing over $\mathbb{D}_n$, $(\,,\,)_{\mathbb{D}_n}$}
\begin{equation}
(\,,\,)_{\mathbb{D}_n}\,:\,C\times C\,\rightarrow \mathbb{D}_n\otimes\mathbb{Q}.
\end{equation}

\end{defn}
Notice that we allow the bilinear pairing to take non-integer values, the reason will become apparent in the following examples.
\begin{oss}\label{oss:cyclo_to_int}
Let $R$ be a $\mathbb{D}_n$-lattice generated by $e_1\,\dots\,e_l$ and let $\omega_1\,\dots,\omega_{\phi(n)}$ be the set of primitive $n$-th roots of unity. Then $R$ has the structure of a free $\mathbb{Z}$-module with generators $e_i\omega_j$ and rank $\phi(n)rank_{\mathbb{D}_n}(R)$. Moreover it is a generalized lattice when endowed with the following pairing:
\begin{equation}
(a,b)_{\mathbb{Z}}=\frac{1}{\phi(n)}\sum_{\rho\in\Gamma_n}\rho(a,b)_{\mathbb{D}_n}.
\end{equation}
Here $\Gamma_n=Gal(\mathbb{D}_n\otimes\mathbb{Q},\mathbb{Q})$. Notice moreover that multiplication by $\omega_n$ defines an isometry of the integer lattice.
\end{oss}

\begin{oss}\label{oss:int_to_cyclo}
Let $R$ be a definite lattice of rank $m$ and let $\varphi\subset O(R)$ be a free isometry of order $n\geq3$, \ie $\varphi^i(v)=v$ if and only if $ v=0\,\text{or}\,i\equiv 0\,mod\,n$. 
Suppose moreover there is an isomorphism $\eta$ of $\mathbb{D}_n$ modules between $R$ and $\mathbb{D}_n^{m/\phi(n)}$, where $\eta\varphi(v)=\omega_n\eta(v)$.
 Then $R$ is a $\mathbb{D}_n$ lattice of rank $\frac{m}{\phi(n)}$.
\end{oss}
If a lattice $R$ can be given both structures we denote as $R$ the $\mathbb{Z}$ lattice and as $R_{\mathbb{D}_n}$ \index{Lattice, cyclotomic, $R_{\mathbb{D}_n}$} the $\mathbb{D}_n$ lattice.

This section cries out for examples, so let us give quite a few:
\begin{ex}\label{an_cyclo}
Let $n\geq 3$ and $A_n\subset\mathbb{Z}^{n+1}$ be the Dynkin lattice as defined in \Ref{ex}{an_form} and let $\varphi$ be the automorphism defined by the permutation $(1\,2\,\dots\,n+1)$ on the standard basis of $\mathbb{Z}^{n+1}$.  If $n+1$ is prime then $A_n$ is a $\mathbb{D}_{n+1}$ lattice of rank $1$.
Let us see two particular cases, $n=2$ and $n=5$:
$A_2$ has basis $\{e_1,e_2\}=\{(1,-1,0),(0,1,-1)\}$ and $\varphi(1,-1,0)=(0,1,-1)$ hence as a $\mathbb{D}_3$ lattice we have $\omega_3e_1=e_2$ and $(e_1,e_1)_{\mathbb{D}_3}=2$, $(e_1,e_2)_{\mathbb{D}_3}=-2\omega_3$.\\
$A_5$ has basis $\{e_1,e_2,e_3,e_4,e_5\}=\{(1,-1,0,0,0,0),\dots,(-1,0,0,0,0,1)\}$ as before. As a $\mathbb{D}_6$ lattice we should have $\varphi^3(e_i)=\omega_6^3(e_i)=-e_i$ but this is not the case.
\end{ex}
\begin{ex}
Let $E_6$ be the Dynkin lattice as in \Ref{ex}{e67_form}. $E_6$ is isometric to the rank 3 $\mathbb{D}_3$ lattice with matrix
\begin{equation}\nonumber
\left(\begin{array}{ccc} 
2 & 0& \frac{2i\sqrt{3}}{3}\\
0& 2 & \frac{2i\sqrt{3}}{3}\\
-\frac{2i\sqrt{3}}{3} & -\frac{2i\sqrt{3}}{3} & 2
\end{array}
\right). 
\end{equation}
\end{ex}
\begin{ex}\label{ex:cox-todd}
Let $K_{12}\subset \mathbb{D}_3^6$ \index{Lattice, Coxeter-Todd, $K_{12}$} be the $\mathbb{D}_3$ lattice generated by
\begin{equation}\nonumber
\frac{1}{\sqrt{2}}(\pm i\sqrt{3},\pm 1,\pm 1,\pm 1,\pm 1,\pm 1),
\end{equation}
where $i\sqrt{3}$ can be in any position and there are an even number of minus signs. This is the Coxeter-Todd lattice. Applying \Ref{oss}{cyclo_to_int} we obtain a lattice with discriminant group $\mathbb{Z}_{/(3)}^6$, which we still call $K_{12}$. The minimal norm of its elements is 4 and its hermitian form over $\mathbb{D}_3$ is the following: 
\begin{equation}
\left(
\begin{array}{cccccc}
4 & 0& 0 & 2 & 2\overline{\omega_3} & 2\overline{\omega_3}\\
0& 4& 0 & 2\overline{\omega_3} & 2 & 2\overline{\omega_3}\\
0& 0& 4 & 2\overline{\omega_3} & 2\overline{\omega_3} & 2\\
2 & 2\omega_3 & 2\omega_3 & 4 & 2 & 2\\
2\omega_3 & 2 & 2\omega_3 & 2 & 4 & 2\\
2\omega_3 & 2\omega_3 & 2& 2 & 2 & 4
\end{array}
\right).
\end{equation}
In its integer form the isometry induced by multiplication by $\omega_3$ acts trivially on the discriminant group.
\end{ex}
The following two lattices are taken from \cite{gs}, where to my knowledge they were explicitly computed for the first time. They correspond to the Co-invariant lattice (cfr. \Ref{defn}{inv_coinv}) of a symplectic automorphism of order 5 and 7 of a $K3$ surface.
\begin{ex}\label{ex:S_5K3}
Let $S_{5.K3}$ \index{Lattice, $S_{5.K3}$} be the lattice associated with the following bilinear form on $\mathbb{Z}^{16}$:
\begin{equation}
\left(
\tiny{\begin{array}{cccccccccccccccc}
 -4 & 2 & 0 & 0 & 0 & -1 & 0 & 0 & 0 & -1 & 0 & 0 & -1 & 1 & -1 & 0 \\
 2 & -4 & 2 & 0 & 5 & 2 & -1 & 0 & 0 & 2 & -1 & 0 & 1 & -1 & 1 & 1 \\
 0 & 2 & -4 & 2 & -5 & -1 & 2 & -1 & 0 & -1 & 2 & -1 & 1 & -1 & 0 & -1 \\
 0 & 0 & 2 & -4 & 0 & 0 & -1 & 2 & 0 & 0 & -1 & 2 & -1 & 1 & 1 & -1 \\
 0 & 5 & -5 & 0 & -50 & 0 & 0 & 0 & 0 & 0 & 0 & 0 & 0 & 0 & 5 & -15 \\
 -1 & 2 & -1 & 0 & 0 & -6 & 4 & -1 & -3 & 0 & 0 & 0 & 0 & 0 & 0 & 0 \\
 0 & -1 & 2 & -1 & 0 & 4 & -6 & 4 & 1 & 0 & 0 & 0 & 0 & 0 & 0 & 0 \\
 0 & 0 & -1 & 2 & 0 & -1 & 4 & -6 & 0 & 0 & 0 & 0 & 0 & 0 & 0 & 0 \\
 0 & 0 & 0 & 0 & 0 & -3 & 1 & 0 & -4 & 3 & -1 & 0 & 2 & 0 & 0 & 0 \\
 -1 & 2 & -1 & 0 & 0 & 0 & 0 & 0 & 3 & -6 & 4 & -1 & -3 & 0 & 0 & 0 \\
 0 & -1 & 2 & -1 & 0 & 0 & 0 & 0 & -1 & 4 & -6 & 4 & 1 & 0 & 0 & 0 \\
 0 & 0 & -1 & 2 & 0 & 0 & 0 & 0 & 0 & -1 & 4 & -6 & 0 & 0 & 0 & 0 \\
 -1 & 1 & 1 & -1 & 0 & 0 & 0 & 0 & 2 & -3 & 1 & 0 & -4 & 3 & -1 & 0 \\
 1 & -1 & -1 & 1 & 0 & 0 & 0 & 0 & 0 & 0 & 0 & 0 & 3 & -6 & 4 & -1 \\
 -1 & 1 & 0 & 1 & 5 & 0 & 0 & 0 & 0 & 0 & 0 & 0 & -1 & 4 & -6 & 4 \\
 0 & 1 & -1 & -1 & -15 & 0 & 0 & 0 & 0 & 0 & 0 & 0 & 0 & -1 & 4 & -6
\end{array}}
\right).
\end{equation}
This lattice has the structure of a $\mathbb{D}_5$ lattice (see \cite{gs}) and is isometric to the following:
\begin{align}\nonumber
\{\,(x_1,x_2,x_3,x_4)\in\mathbb{D}_5^4, \,x_1\equiv x_2\equiv 2x_3\equiv 2x_4\,\,mod\,(1-\omega_5)\,\text{and}\,\\\nonumber
 (3-\omega_5)(x_1+x_2)+x_3+x_4\equiv 0\,mod\,(1-\omega_5)^2\}.
\end{align}
With the following hermitian form:
\begin{equation}\nonumber
(x,y)=x_1\overline{y_1}+x_2\overline{y_2}+fx_3\overline{fy_3}+fx_4\overline{fy_4},\,\text{where}\,\,\, f=1-(\omega_5^2+\omega_5^3).
\end{equation}
\end{ex}
\begin{ex}\label{ex:S_7K3}
Let $S_{7.K3}$ \index{Lattice, $S_{7.K3}$} be the lattice associated with the following bilinear form on $\mathbb{Z}^{18}$:
\begin{equation}
\left(
\tiny{\begin{array}{cccccccccccccccccc}
 -4 & 2 & 0 & 0 & 0 & 0 & 0 & -1 & 0 & 0 & 0 & 0 & 0 & 1 & -1 & 0 & 0 & 0 \\
 2 & -4 & 2 & 0 & 0 & 0 & 0 & 2 & -1 & 0 & 0 & 0 & 0 & -1 & 1 & 1 & -1 & 0 \\
 0 & 2 & -4 & 2 & 0 & 0 & 7 & -1 & 2 & -1 & 0 & 0 & 0 & 0 & 0 & -1 & 1 & 1 \\
 0 & 0 & 2 & -4 & 2 & 0 & -7 & 0 & -1 & 2 & -1 & 0 & 1 & -1 & 0 & 0 & 0 & -1 \\
 0 & 0 & 0 & 2 & -4 & 2 & 0 & 0 & 0 & -1 & 2 & -1 & -1 & 1 & 1 & -1 & 0 & 0 \\
 0 & 0 & 0 & 0 & 2 & -4 & 0 & 0 & 0 & 0 & -1 & 2 & -1 & 0 & -1 & 1 & 1 & -1 \\
 0 & 0 & 7 & -7 & 0 & 0 & -98 & 0 & 0 & 0 & 0 & 0 & 0 & 0 & 0 & 0 & 7 & -21 \\
 -1 & 2 & -1 & 0 & 0 & 0 & 0 & -6 & 4 & -1 & 0 & 0 & 0 & 0 & 0 & 0 & 0 & 0 \\
 0 & -1 & 2 & -1 & 0 & 0 & 0 & 4 & -6 & 4 & -1 & 0 & 0 & 0 & 0 & 0 & 0 & 0 \\
 0 & 0 & -1 & 2 & -1 & 0 & 0 & -1 & 4 & -6 & 4 & -1 & 0 & 0 & 0 & 0 & 0 & 0 \\
 0 & 0 & 0 & -1 & 2 & -1 & 0 & 0 & -1 & 4 & -6 & 4 & -1 & 0 & 0 & 0 & 0 & 0 \\
 0 & 0 & 0 & 0 & -1 & 2 & 0 & 0 & 0 & -1 & 4 & -6 & 3 & 0 & 0 & 0 & 0 & 0 \\
 0 & 0 & 0 & 1 & -1 & -1 & 0 & 0 & 0 & 0 & -1 & 3 & -4 & 3 & -1 & 0 & 0 & 0 \\
 1 & -1 & 0 & -1 & 1 & 0 & 0 & 0 & 0 & 0 & 0 & 0 & 3 & -6 & 4 & -1 & 0 & 0 \\
 -1 & 1 & 0 & 0 & 1 & -1 & 0 & 0 & 0 & 0 & 0 & 0 & -1 & 4 & -6 & 4 & -1 & 0 \\
 0 & 1 & -1 & 0 & -1 & 1 & 0 & 0 & 0 & 0 & 0 & 0 & 0 & -1 & 4 & -6 & 4 & -1 \\
 0 & -1 & 1 & 0 & 0 & 1 & 7 & 0 & 0 & 0 & 0 & 0 & 0 & 0 & -1 & 4 & -6 & 4 \\
 0 & 0 & 1 & -1 & 0 & -1 & -21 & 0 & 0 & 0 & 0 & 0 & 0 & 0 & 0 & -1 & 4 & -6
\end{array}}
\right).
\end{equation}
This lattice has the structure of a $\mathbb{D}_7$ lattice (see \cite{gs}) and is isometric to the following:
\begin{align}\nonumber
\{\,(x_1,x_2,x_3)\in\mathbb{D}_7^3, \,x_1\equiv x_2\equiv 6x_3\,\,\,mod\,(1-\omega_7)\,\text{and}\,\\\nonumber
 (1+5\omega_7)x_1+3x_2+2x_3\equiv 0\,mod\,(1-\omega_7)^2\}.
\end{align}
With the following hermitian form:
\begin{equation}\nonumber
(x,y)=x_1\overline{y_1}+f_1x_2\overline{f_1y_2}+f_2x_3\overline{f_2y_3},
\end{equation}
where $f_1=3+2(\omega_7+\overline{\omega}_7)+(\omega_7^2+\overline{\omega}_7^2)$ and $f_2=2+(\omega_7+\overline{\omega}_7)$.
\end{ex}

\begin{ex}\label{ex:craig_leech}
Let $\Lambda$ be the Leech lattice, Craig \cite{craig} proved that it has the structure of a $\mathbb{D}_{39}$ lattice, however its proof is fairly complicated. Notice that this implies that $\Lambda$ has also the structure of $\mathbb{D}_3$ and $\mathbb{D}_{13}$ lattice, a fact that will lead to \Ref{lem}{max_p_order}. 
\end{ex}
\begin{ex}\label{ex:cyclo_orderp}
Let us consider a lattice over $\mathbb{D}_p$ of rank 1, where $p$ is a prime. It is generated by an element $v$ of square $a\in\mathbb{D}_p$. Let us look at its integer form: it has rank $p-1$ and its basis is $v,\omega_pv,\dots,\omega_p^{p-2}v$. Let us suppose $a\in\mathbb{Z}$. In this case we have \\$(v,\omega_pv)_{\mathbb{Z}}=-a/2.$ Thus this lattice is nothing else than $A_{p-1}(a/2)$ in a different form, see \cite[Section 4.6]{con}. This implies moreover that the discriminant group of a rank 1 cyclotomic lattice has (in its integer form) $p-1$ generators as soon as $|a|>2$.
\end{ex}
We need to consider one more lattice, which was first introduced by Wall \cite{wal} and studied also by Nebe and Plesken \cite[page 65]{neb}:
\begin{ex}\label{ex:wall}
Let $W$ \index{Lattice, Wall's, $W$} be the lattice associated to the following bilinear form on $\mathbb{Z}^{18}$:
\begin{equation}\label{wall_lattice}
\left(
\tiny{\begin{array}{cccccccccccccccccc}
 4 & 2 & -2 & -2 & -2 & 0 & 0 & 0 & -2 & 1 & -2 & -2 & -2 & -2 & 2 & -2 & 2 & 2 \\
 2 & 4 & -2 & 0 & -2 & -1 & -1 & -1 & -2 & -1 & 0 & 0 & 0 & -1 & 2 & 0 & 1 & 2 \\
 -2 & -2 & 4 & 1 & 2 & 1 & -1 & 1 & 1 & -1 & 0 & 2 & 2 & 2 & -2 & 0 & 0 & -2 \\
 -2 & 0 & 1 & 4 & 2 & 1 & 0 & 1 & 1 & -1 & 0 & 2 & 2 & 0 & -1 & 1 & -2 & -2 \\
 -2 & -2 & 2 & 2 & 4 & 2 & 1 & 1 & 1 & -1 & 0 & 2 & 2 & 0 & -2 & 0 & -2 & -2 \\
 0 & -1 & 1 & 1 & 2 & 4 & 0 & 0 & 0 & 0 & 0 & 0 & 0 & -1 & -2 & -2 & -1 & 0 \\
 0 & -1 & -1 & 0 & 1 & 0 & 4 & 1 & 1 & 0 & 0 & 0 & -1 & -2 & 0 & 0 & 0 & -1 \\
 0 & -1 & 1 & 1 & 1 & 0 & 1 & 4 & 0 & 1 & -2 & 0 & 1 & -1 & 1 & -1 & 1 & -1 \\
 -2 & -2 & 1 & 1 & 1 & 0 & 1 & 0 & 4 & 1 & 0 & 0 & 1 & 1 & -2 & 2 & -1 & -1 \\
 1 & -1 & -1 & -1 & -1 & 0 & 0 & 1 & 1 & 4 & -2 & -2 & -1 & 0 & 1 & 0 & 1 & 1 \\
 -2 & 0 & 0 & 0 & 0 & 0 & 0 & -2 & 0 & -2 & 4 & 1 & 0 & 1 & -1 & 1 & -1 & 0 \\
 -2 & 0 & 2 & 2 & 2 & 0 & 0 & 0 & 0 & -2 & 1 & 4 & 2 & 1 & -1 & 1 & -1 & -2 \\
 -2 & 0 & 2 & 2 & 2 & 0 & -1 & 1 & 1 & -1 & 0 & 2 & 4 & 1 & -1 & 1 & -1 & -1 \\
 -2 & -1 & 2 & 0 & 0 & -1 & -2 & -1 & 1 & 0 & 1 & 1 & 1 & 4 & -1 & 2 & 0 & -1 \\
 2 & 2 & -2 & -1 & -2 & -2 & 0 & 1 & -2 & 1 & -1 & -1 & -1 & -1 & 4 & 0 & 2 & 1 \\
 -2 & 0 & 0 & 1 & 0 & -2 & 0 & -1 & 2 & 0 & 1 & 1 & 1 & 2 & 0 & 4 & -1 & -1 \\
 2 & 1 & 0 & -2 & -2 & -1 & 0 & 1 & -1 & 1 & -1 & -1 & -1 & 0 & 2 & -1 & 4 & 1 \\
 2 & 2 & -2 & -2 & -2 & 0 & -1 & -1 & -1 & 1 & 0 & -2 & -1 & -1 & 1 & -1 & 1 & 4
\end{array}}
\right).
\end{equation}
This lattice has $Aut(W)=3^{1+4}:Sp_4(\mathbb{Z}_{/(3)}).2$ (in the notation of \Ref{sec}{niemeier}). We also have $A_W=(\mathbb{Z}_{/(3)})^5$, therefore it can be embedded in a positive definite unimodular lattice of rank $24$. 
\end{ex}

\section{Isometries, Invariant and Co-invariant Lattices}
In this section we analyze two kind of lattices linked to an isometry, namely the co-invariant and invariant lattices. We will give also some proofs related to some lattices useful in \Ref{cap}{deformations} and we will present an easy construction of fundamental relevance in the proof of \Ref{thm}{sporadic}.

\begin{defn}\label{defn:inv_coinv}
Let $R$ be a lattice and let $G\subset O(R)$. Then we define $T_G(R)=R^G$ \index{Lattice, Invariant, $T_G(R)$} as the invariant lattice of $G$ and \\$S_G(R)=T_G(R)^\perp$ \index{Lattice, Co-invariant, $S_G(R)$} as the co-invariant lattice. 
\end{defn}
\begin{oss}\label{oss:G_tors}
Let $R$ be a lattice, and let $G\subset O(R)$. Then the following hold:
\begin{itemize}
\item $T_G(R)$ contains $\sum_{g\in G}gv$ for all $v\in R$.
\item $S_G(R)$ contains $v-gv$ for all $v\in R$ and all $g\in G$.
\item If $R$ is definite then $T_G(R)$ and $S_G(R)$ are nondegenerate.
\item $R/(T_G(R)\oplus S_G(R))$ is of $|G|$-torsion.
\item Suppose $G$ is of prime order $p$ and $R$ is definite, then $S_G(R)$ is a $\mathbb{Z}[\omega_p]$ lattice.
\end{itemize}
\begin{proof}
It is obvious that $\sum_{g\in G}gv$ is $G$-invariant for all $v\in R$. Let $w\in T_G(R)$, since $g$ is an isometry we have $(w,v)=(gw,gv)=(w,gv)$ for all $v\in R$ and all $g\in G$. Therefore $v-gv$ is orthogonal to all $G$-invariant vectors, hence it lies in $S_G(R)$.
Obviously whenever $R$ is definite all of its sublattices are nondegenerate. 
Let $t\in R$, we can write $|G|t=\sum_{g\in G} g(t) + \sum_{g\in G}(t-g(t))$, where the first term lies in $T_G(R)$ and the second in $S_G(R)$.
Finally if $|G|=p$ we let $g$ be one of its generators. $g$ acts freely on $S_G(R)$, therefore by \Ref{oss}{int_to_cyclo} it can be defined as a $\mathbb{Z}[\omega_p]$ lattice.
\end{proof}
\end{oss}

\begin{lem}\label{lem:biggen_g}
Let $L$ be as in \Ref{ex}{K32_lat} and let $L'$ be the Mukai lattice. Let $g\in O(L)$, then there exists an embedding $L\subset L'$ and an isometry $\overline{g}\,\in\,O(L')$ such that $\overline{g}_{|L}=g$.
\begin{proof}
The isometry $g$ induces an automorphism of the discriminant group $A_L$. Since $A_L=\mathbb{Z}_{/(2)}$ this automorphism is the identity. Let $[v/2]$ be a generator of $A_L$ such that $v^2=-2$. We then have $g([v/2])=[v/2]$ \ie $g(v)=v+2w$. Consider now a lattice of rank 1 generated by an element $x$ of square 2, its discriminant group is still $\mathbb{Z}_{/(2)}$ and is generated by $[x/2]$ with discriminant form given by $q(x/2)=1/2$.\\ Notice that $L\oplus\,\mathbb{Z}x$ has an overlattice isometric to $L'$ which is generated by $L$ and $\frac{x+v}{2}$.\\ We now extend $g$ on $L\oplus x$ by imposing $g(x)=x$ and we thus obtain an extension $\overline{g}$ of $g$ to $L'$ defined as follows:
\begin{align}\nonumber
\overline{g}(e) &= g(e)\,\,\forall\,e\,\in\,L,\\ \nonumber
\overline{g}(x)&= x,\\ \nonumber
\overline{g}(\frac{x+v}{2})&= \frac{x+g(v)}{2}.
\end{align}
\end{proof}
\end{lem}

\begin{oss}\label{oss:s_in_24}
Let $L$ be as in \Ref{ex}{K32_lat} and let $G\subset O(L)$. Then there exists a primitive embedding $L\rightarrow L'\cong U^4\oplus E_8(-1)^2$ such that $G$ extends to a group of isometries of $L'$ and $S_G(L)=S_G(L')$.
\begin{proof}
Let $x$ be a vector of square $2$ and $v\in L$ a vector of square $-2$ such that $(v,L)=2\mathbb{Z}$. Let $L'$ be the overlattice of $L\oplus\mathbb{Z}x$ generated by $L$ and $\frac{x+v}{2}$ and let us extend the action of $G$ to $L'$ as in \Ref{lem}{biggen_g}. A direct computation shows $S_G(L)=S_G(L')$.
\end{proof}
\end{oss}

\subsection{$\mathcal{S}$-lattices}
In this subsection we analyze briefly a few sublattices of the Leech lattice which arise as $T_G(\Lambda)$ for some interesting groups $G$. Let us start with the basics:
\begin{defn}\label{defn:slat}
Let $M\subset \Lambda$. Then $M$ is a $\mathcal{S}$-lattice \index{Lattice, Special sublattice of $\Lambda$, $\mathcal{S}$-lattice} if all elements of $M$ are congruent modulo $2\Lambda$ to an element of $M$ of norm $0,-4$ or $-6$.
\end{defn}
There are not many examples of $\mathcal{S}$-lattices and they were classified by Curtis:
\begin{lem}\cite{cur}
Up to isomorphisms there are 12 $\mathcal{S}$-lattices inside $\Lambda$.
\end{lem}
He classified also their stabilizers and their automorphism groups inside $Co_0$, a full table can be found in \cite[page 180]{atlas}.
For our purpose it is better to give an explicit presentation of the Leech Lattice $\Lambda$:
\begin{ex}\label{ex:leech_24}
Let us consider the vector space $\mathbb{R}^{W}$, where $W=\mathbb{P}^1(\mathbb{Z}_{/(23)})$ is a set with 24 elements and let us endow it with  a quadratic form defined as the opposite of the euclidean form. Let $Q\subset W$ be the set whose elements are quadratic residues modulo $23$ and $0$, and let $a=8^{-1/2}$. Then $\Lambda\subset \mathbb{R}^W$ is spanned by the following elements:
\begin{align*}
a(2,\dots,2,0,\dots,0), & \text{    where the twelve non zero elements are}\\
& \text{supported on a translate of $Q$ by an element of $W$},\\
a(-3,1,\dots,1), & \\
a(\pm4,\pm4,0,\dots,0). & 
\end{align*}  
\end{ex} 
Let us introduce a piece of notation: a $\mathcal{S}$-lattice $M$ is denoted $2^i3^j$ \index{Lattice, Notation for $\mathcal{S}$-lattices, $2^i3^j$} if (up to sign) it contains $i$ vectors of norm $-4$ and $j$ vectors of norm $-6$.
\begin{ex}
The easiest example possible is that of a lattice $M=(-4)=2^1$ in the above notation. The condition of \Ref{defn}{slat} is trivially satisfied and $Aut(M)=\pm Id$, $Stab(M)=Co_2$ \ie $M=T_{Co_2}(\Lambda)$.
\end{ex} 

\begin{ex}\label{ex:slat_5rk4}
Let us consider the $\mathcal{S}$-lattice $M=2^53^{10}$, it has rank 4 and it is $T_G(\Lambda)$ where $G$ is an extension of $(\mathbb{Z}_{/(5)})^3$ with $\mathbb{Z}_{/(4)}$. We wish to remark that $G$ contains an element of the conjugacy class $5C$ in the notation of \cite{atlas}. We will later denote $M^\perp$ as the lattice $S_{5.exo}$ \index{Lattice, $S_{5.exo}$}, which is shown in \Ref{ex}{A46} to be just $S_{5C}(\Lambda)$. Moreover $M$ is isometric to the following
\begin{equation*}
\left(
\begin{array}{cccc}
-4 & -1 & -1 & 1\\
-1 & -4 & 1 & -1\\
-1 & 1 & -4 & -1\\
1 & -1 & -1 & -4
\end{array}
\right).
\end{equation*} 
From Nipp's \cite{nip} list of definite quadratic forms we have that $M$ is the unique lattice in its genus, moreover $M\oplus M$ has $E_8(-1)$ as overlattice. 
\end{ex}

\begin{ex}\label{ex:slat_3rk4}
The $\mathcal{S}$-lattice $M=2^93^6$ is a lattice of rank 4 and it is the stabilizer of a group $G\subset Co_0$ which is a nontrivial extension of $(\mathbb{Z}_{/(3)})^4$ with $A_6$. If we consider $\Lambda$ as the lattice defined in \Ref{ex}{leech_24} then $M$ is spanned by the following 9 elements:
\begin{align*}
a(0,0,0,0,0,0,0,0,-4,0,0,0,0,0,0,0,4,0,0,0,0,0,0,0), & \\
a(4,0,0,0,0,0,0,0,0,0,0,0,0,0,0,0,-4,0,0,0,0,0,0,0), & \\
a(-4,0,0,0,0,0,0,0,4,0,0,0,0,0,0,0,0,0,0,0,0,0,0,0), & \\
a(0,0,0,0,0,0,0,0,2,2,2,2,0,0,0,0,-2,-2,-2,-2,0,0,0,0), & \\
a(-2,-2,-2,-2,0,0,0,0,0,0,0,0,0,0,0,0,2,2,2,2,0,0,0,0), &\\
a(2,2,2,2,0,0,0,0,-2,-2,-2,-2,0,0,0,0,0,0,0,0,0,0,0,0), &\\
a(0,0,0,0,0,0,0,0,2,-2,-2,-2,0,0,0,0,-2,2,2,2,0,0,0,0), & \\
a(-2,2,2,2,0,0,0,0,0,0,0,0,0,0,0,0,2,-2,-2,-2,0,0,0,0), & \\
a(2,-2,-2,-2,0,0,0,0,-2,2,2,2,0,0,0,0,0,0,0,0,0,0,0,0), &
\end{align*}
where $a=8^{-1/2}$. A direct computation shows that it is isometric to the lattice 
\begin{equation*}
\left(\begin{array}{cccc} 
-4 & 2 & -2 & 1 \\
2 & -4 & 1 & -2 \\
-2 & 1 & -4 & 2 \\
1 & -2 & 2 & -4 \\
\end{array}\right).
\end{equation*}
A look at Nipp's table \cite{nip} shows again that it is unique in its genus. This time however $M\oplus M$ has not an unimodular overlattice, however $M\oplus A_2(-1)\oplus A_2(-3)$ does. Notice moreover that $A_2(-1)\oplus A_2(-3)$ is again unique in its genus by \cite{nip}.
\end{ex}

\begin{ex}\label{ex:slat_rk6}
The $\mathcal{S}$-lattice $M=2^{27}3^{36}$ is a lattice of rank 6 and discriminant group $(\mathbb{Z}_{/(3)})^5$ and it is the stabilizer of a group $G\subset Co_0$, where $G$ is a nontrivial extension of $(\mathbb{Z}_{/(3)})^5$ with $\mathbb{Z}_{/(2)}$. Its orthogonal inside $\Lambda$ is a lattice which contains the group $O(E_6)$ in its automorphism group. We do not give a direct proof that its orthogonal is isometric to the lattice $W(-1)$ of \Ref{ex}{wall}, however it is implied by the following facts:
\begin{enumerate}
\item $G$ contains an element of conjugacy class $3C$ in the notation of \cite{atlas}.
\item $W(-1)$ embeds into a negative definite unimodular lattice $N$ of rank $24$ such that $S_\varphi(N)=W(-1)$, where $\varphi$ is induced by an isometry of order 3 of $W(-1)$ acting trivially on $W(-1)^{\perp_N}$.
\item There are $24$ lattices in the genus of $N$, see \Ref{sec}{niemeier} for details, and only $\Lambda$ has an element $\varphi$ of order 3 (of conjugacy class $3C$) such that $S_\varphi(N)$ has rank 18.
\end{enumerate}  
\end{ex}

\section{Eichler transvections}
In this section we make good use of a certain class of Isometries, known as Eichler's transvections \cite{eich}. Our exposition follows very much \cite{ghs1}, where all the proofs we omit can be found.
\begin{defn}
Let $R$ be a lattice and let $e\in R$ be an isotropic vector. Let $a\in\,e^\perp$. The map
\begin{equation}
t'(e,a)\,:\,v\,\rightarrow\,v-(a,v)e
\end{equation}
defines an isometry of $e^\perp$. 
\end{defn}
\begin{lem}
$t'(e,a)$ extends to a unique isometry $t(e,a)$ \index{Lattice, Eichler's transvection, $t(e,a)$} of $R$, called Eichler's transvection.
\begin{proof}
Let us define
\begin{equation}
t(e,a)\,:\,v\,\rightarrow\,v-(a,v)e+(e,v)a-\frac{1}{2}(a,a)(e,v)e.
\end{equation}
Its restriction on $e^\perp$ is $t'(e,a)$.
\end{proof} 
\end{lem}
\begin{lem}\label{lem:eich_u2}
Let us consider the lattice $R=U\oplus U_1$, where $U_1\cong U$. Then for all $v\in R$ there exists an isometry $g$ of $R$ generated by Eichler's transvections such that $g(v)\in U_1$.
\end{lem}
\begin{oss}
\Ref{lem}{eich_u2} can be easily extended to $R=U(n)\oplus U_1(n)$.
\end{oss}

The following is known as Eichler's criterion, see \cite[Proposition 3.3]{ghs1} for a proof. 

\begin{lem}
\label{lem:ghs_orbit}
Let $T$ be a lattice such that $T\cong U^2\oplus N$ for some lattice $N$ and let $v,w\in T$ be two primitive vectors such that the following hold:
\begin{itemize}
\item $v^2=w^2$.
\item $(v,T)\cong m\mathbb{Z}\cong (w,T)$.
\item $[\frac{v}{m}]=[\frac{w}{m}]$ in $A_T$.
\end{itemize}
Then there exists an isometry $g$ of $T$ such that $g(v)=w$.
\end{lem}
We will also need a slight generalization of it to prove \Ref{lem}{emblemma_3} and \Ref{lem}{emblemma_5}:
\begin{lem}\label{lem:ghs_orbit_gen}
Let $T\cong U(n)^2\oplus N$ for some lattice $N$ and some integer $n$ and let $v,w\in T$ be two primitive vectors such that the following hold:
\begin{itemize}
\item $v^2=w^2$,
\item $(v,T)\cong m\mathbb{Z}\cong (w,T)$,
\item There exists an isometry $h$ such that $[\frac{v}{m}]=h[\frac{w}{m}]$ in $A_T$.
\end{itemize}
Then there exists an isometry $g$ of $T$ such that $g(v)=w$.
\begin{proof}
First of all we use the generalized version of \Ref{lem}{eich_u2} to obtain two isometries $f$ and $f'$ such that $f(v)=v'$ and $f'(w)=w'$ are both orthogonal to the first copy of $U(n)$. By hypothesis the isometry $f^{-1}\circ h\circ f'^{-1}$ sends $[\frac{w'}{m}]$ to $[\frac{v'}{m}]$ and let $f^{-1}\circ h\circ f'^{-1}(w')=w"$. Let $d,e$ be a basis of the first copy of $U(n)$, then let $k$ be the following isometry:
\begin{equation}
v'\,\stackrel{t(e,\overline{v}')}{\rightarrow}\,(v-me)\,\stackrel{t(d,(v'-w")/m)}{\rightarrow}\,(w"-me)\stackrel{t(e,-\overline{w}")}{\rightarrow}\,w".
\end{equation}
Here $\overline{v}'$ is an element orthogonal to the first copy of $U$ such that $(v',\overline{v}')=m$ and the same goes for $\overline{w}"$ and $w"$.
Therefore 
$h^{-1}\circ f\circ k\circ f$ sends $v$ to $w$. 
\end{proof}
\end{lem}

\section{Niemeier lattices and Leech-type lattices}\label{sec:niemeier}
In this section we recall Niemeier list of negative definite even unimodular lattices in dimension 24 and we introduce a class of lattices which will be of fundamental interest in the rest of the section.
Detailed information about these lattices can be found in \cite[Chapter 16]{con} and in \cite[Section 1.14]{nik2}.\\
\begin{defn}\label{defn:leech_group}
Let $M$ be a lattice and let $G\subset\,O(M)$. Then $M$ is a Leech type lattice with respect to $G$ if the following are satisfied:
\begin{itemize}
\item $M$ is negative definite.
\item $M$ contains no vectors of square $-2$.
\item $G$ acts trivially on $A_M$.
\item $S_G(M)=M$.
\end{itemize} 
Moreover we call $(M,G)$ a Leech couple and $G$ a Leech type group.
\end{defn}

Notice that $(\Lambda, Co_0)$ is a Leech couple.
Now we recall Niemeier's list of definite even unimodular lattices of dimension 24. Usually they are defined as positive definite lattices but we will use them as negative definite ones. All of these lattices can be obtained by specifying a 0 or 24 dimensional Dynkin diagram such that every semisimple component has a fixed Coxeter number, in \Ref{tab}{nieme} we recall the possible choices. Having the Dynkin lattice $A(-1)$ of the lattice $N$ we obtain it by adding a certain set of glue vectors, which are a subset $G(N)$ \index{Lattice, Group of glue vectors of $N$, $G(N)$} of $A^{\vee}/A$. The precise definition of the glue vectors can be found in \cite[Section 4]{con} and we keep the same notation contained therein. Notice that the set of glue vectors forms an additive subgroup of $A^{\vee}/A$.\\
Another fundamental data is what we call maximal Leech-type group $Leech(N)$, \ie the maximal subgroup $G$ of $Aut(N)$ such that $(S_G(N),G)$ is a Leech-type couple. It is a well known fact that this group is obtained as $Aut(N)/W(N)$, where $W(N)$ is the Weyl group generated by reflections on $-2$ vectors. These groups where first computed by Erokhin \cite{ero}.\\
This data is summarized in \Ref{tab}{nieme}, let us explain briefly the notation used therein: for the Leech-type group we used standard notation from \cite{atlas}, where $n$ denotes a cyclic group of order $n$, $p^n$ \index{Group, Elementary $p$-group of order $p^n$, $p^n$} denotes an elementary $p$-group of order $p^n$, $G.H$ \index{Group, Extension of $G$ with $H$, $G.H$} denotes any group $F$ with a normal subgroup $G$ such that $F/G=H$ and $L_m(n)$ \index{Group, Special linear, $L_m(n)$, $PSL_m(\mathbb{Z}_{/(n)})$} denotes the group $PSL_m$ over the finite field with $n$ elements. $M_n$ \index{Group, Mathieu's, $M_n$} denotes the Mathieu group on $n$ elements and $Co_n$ denotes Conway groups.\\   
Regarding the glue codes we kept the notation of \cite{con}, hence a glue code $[abc]$ means a vector $(g,h,f)$ where $g$ is the glue vector of type $a$, $h$ is the one of type $b$ and $f$ of type $c$. Moreover $[(abc)]$ indicates all glue vectors obtained from cyclic permutations of $\{a,b,c\}$, hence $[abc],[bca],[cab]$.\\
\begin{table}[ht]\label{tab:nieme}
\caption{Niemeier lattices and their Leech automorphisms}
\begin{tabular}{|c|c|c|c|c|}
 \hline
Name & Dynkin & Leech-type& Coxeter& Generating glue code\\ 
 & diagram & Group &  Number & \\

\hline
$N_{1}$\index{Lattice, Niemeier, $N_i$} & $D_{24}$ & $1$ & $46$ & $[1]$\\ \hline
$N_{2}$ &$D_{16}E_8$ & $1$ & $30$ & $[10]$\\ \hline
$N_{3}$ &$E_8^3$ & $S_3$ & $30$ & $[000]$\\ \hline
$N_{4}$ &$A_{24}$ & $2$ & $25$ & $[5]$\\ \hline
$N_{5}$ &$D_{12}^2$ & $2$ & $22$ & $[12],[21]$\\ \hline
$N_{6}$ &$A_{17}E_7$ & $2$ & $18$ & $[31]$\\ \hline
$N_{7}$ &$D_{10}E_7^2$ & $2$ & $18$ & $[110],[301]$\\ \hline
$N_{8}$ &$A_{15}D_9$ & $2$ & $16$ & $[21]$ \\ \hline
$N_{9}$ &$D_8^3$ & $S_3$ & $14$ & $[(122)]$ \\ \hline
$N_{10}$ &$A_{12}^2$ & $4$ & $13$ & $[15] $\\ \hline
$N_{11}$ &$A_{11}D_7E_6$ & $2$ & $12$ & $[111]$\\ \hline
$N_{12}$ &$E_6^4$ & of order $48$ & $12$ & $[1(012)]$ \\ \hline
$N_{13}$ &$A_9^2D_6$ & $2^2$  & $10$ & $[240],[501],[053]$ \\ \hline
$N_{14}$ &$D_6^4$ & $S_4$ & $10$ & $[$even perm. of $\{0,1,2,3\}]$ \\ \hline
$N_{15}$ &$A_8^3$ & $S_3\times 2$  & $9$ & $[(114)]$\\ \hline
$N_{16}$ &$A_7^2D_5^2$ & $2^3$ & $8$ & $[1112],[1721]$\\ \hline
$N_{17}$ &$A_6^4$ & of order $24$ & $7$ & $[1(216)]$ \\ \hline
$N_{18}$ &$A_5^4D_4$ & as $N_{12}$ & $6$ & $[2(024)0],[33001],[30302],[30033] $ \\ \hline
$N_{19}$ &$D_4^6$ & $3\times S_6$ &$6$ & $[111111],[0(02332)]$ \\ \hline
$N_{20}$ &$A_4^6$ & $2.L_2(5).2$ & $5$ & $[1(01441)]$\\ \hline
$N_{21}$ &$A_3^8$ & $2^3.L_2(7).2$ & $4$ & $[3(2001011)]$ \\ \hline
$N_{22}$ &$A_2^{12}$ & $2.M_{12}$ &$3$ & $[2(11211122212)]$\\ \hline
$N_{23}$ &$A_1^{24}$ &  $M_{24}$ & $2$ & $[1(00000101001100110101111)]$ \\ \hline
$\Lambda$ &$\emptyset$ & $Co_0$ &$0$ & $\emptyset$ \\ \hline
\end{tabular}
\end{table}

By \Ref{lem}{nik_immerge1} all of the Niemeier lattices can be defined as primitive sublattices of $\Pi_{1,25}\cong U\oplus E_8(-1)^3$ \index{Lattice, $\Pi_{1,25}$} by specifying a primitive isotropic vector $v$ and setting $N=(v^\perp\cap\Pi_{1,25})/v$.
\begin{ex}\label{ex:leech_26}
Let $\Pi_{1,25}\subset\mathbb{R}^{26}$ (the first coordinate of $\mathbb{R}^{26}$ is the positive definite one) be as before and let
\begin{align}\nonumber
v = &(17,1,1,1,1,1,1,1,1,3,3,3,3,3,3,3,3,3,5,5,5,5,5,5,5,5)\\\nonumber
w = &(70,0,1,2,3,4,5,\dots,24)
\end{align}
be two isotropic vectors in the standard basis of $\mathbb{R}^{26}$. Then 
\begin{equation}\nonumber
\Lambda\cong (w^\perp\cap\Pi_{1,25})/w
\end{equation}
 and 
\begin{equation}\nonumber 
N_{15}\cong (v^\perp\cap\Pi_{1,25})/v.
\end{equation}
\end{ex}

\subsection{The "holy" construction}\label{ssec:holy}
In this subsection we give a few different constructions of the Leech lattice $\Lambda$ arising from the other Niemeier lattices. These constructions will be instrumental in the proof of \Ref{thm}{prime_autom_k32}.\\ The detailed construction is contained in \cite[Section 24]{con}, in the present paper we just sketch it:
Let $A_n$ be a Dynkin lattice defined by
\begin{equation}\nonumber
A_n=\{(a_1,\dots,a_{n+1})\,\in\,\mathbb{Z}^{n+1},\,\,\sum{a_i}=0\}.
\end{equation}
And let $f_j$ be the vector with $-1$ in the $j-$th coordinate and $1$ in the $(j+1)-$th, zero otherwise. Let moreover $f_0=(1,0,\dots,0,-1)$. In general the $f_i$ form a set of extended roots for the Dynkin lattice.\\ Let $g_0=h^{-1}(-\frac{1}{2}n,-\frac{1}{2}n+1,\dots,\frac{1}{2}n)$ where $h$ is the Coxeter number of $A_n$ and let the $g_i$'s be a cyclic permutation of coordinates of $g_0$.
Now let $A_n(-1)^m$ be a 24 dimensional lattice and let $h_k=(g_{i_1},\dots,g_{i_m})$ where $[i_1i_2\dots i_m]$ is a glue code obtained from \Ref{tab}{nieme}. Let $f^j_i=(0,\dots,0,f_i,0,\dots,0)$ where $f_i$ belongs to the $j-$th copy of $A_n$. Let $m_i^j$ and $n_w$ be integers.\\
Then the following holds: the set of vectors satisfying
\begin{equation}\label{holy_nieme}
\sum_{j=1}^m\,\sum_{i} m_i^jf_i^j+\sum_w n_wh_w,\,\,\,\,\sum_{w}n_w=0 
\end{equation}
is isometric to the Niemeier lattice with Dynkin diagram $A_n^m$. 
While the set of vectors
\begin{equation}\label{holy_leech}
\sum_{j=1}^m\,\sum_{i} m_i^jf_i^j+\sum_w n_wh_w,\,\,\,\,\sum_{w}n_w+\sum_{i,j}m_i^j=0
\end{equation}
is isometric to the Leech lattice $\Lambda$. We call the set defined by \eqref{holy_leech} the holy construction of $\Lambda$ with hole \eqref{holy_nieme}.\\
Moreover the glue code provides several automorphisms of the Leech lattice, where the action of $t\in G(N)$ is given by sending $h_w$ to $h_{w+t}$.\\
\begin{oss}\label{oss:congruenza_busona}
For all sets of extended roots of a Dynkin lattice there exists a linear combination $\sum_{i}a_if_i=0$ such that $\sum_{i}a_i=h$, the coxeter number. This implies that \eqref{holy_nieme} and \eqref{holy_leech} can be rewritten as congruences modulo $h$. Notice moreover that this implies also that the lattice $N\cap\Lambda$ has index $h$ both inside $N$ and inside $\Lambda$.
\end{oss}
This construction is really useful to explicit the action of some elements of $Co_0$ on $\Lambda$, namely in the following examples:
\begin{ex}\label{ex:p3E83}
Let us apply this construction to the lattice $E_8(-1)^3$ and let $\varphi$ be an order 3 permutation of the 3 copies of $E_8(-1)$. With the holy construction with hole $N_3$ it induces an automorphism $\varphi$ of $\Lambda$ of order 3 which fixes the only glue vector $g_0$. A direct computation shows that $T_\varphi(N_3)\cong E_8(-3)$ and $S_\varphi(N_3)\cong S_\varphi(\Lambda)=\{a-\varphi(a),\,a\in E_8(-1)^3\}$. Let us call this lattice $S_{3.exo}$ \index{Lattice, $S_{3.exo}$}, it is $S_g(\Lambda)$ for any $g\in Co_0$ in the conjugacy class $3D$ (in the notation of \cite{atlas}).
\end{ex}
\begin{ex}\label{ex:A212}
Let us apply this construction to the lattice $A_2(-1)^{12}$, we then have $G(N)=(\mathbb{Z}_{/(3)})^6$ acting on $\Lambda$. The normalizer of this group (inside $Co_1$) is one of the maximal subgroups of $Co_1$, therefore its structure is analyzed in \cite{atlas}. The elements of $G(N)$ fall under three conjugacy classes labeled $3A,3B$ and $3C$. Each conjugacy class has respectively $24,262$ and $440$ representatives. Therefore we can compute the rank of the invariant lattice inside $\Lambda$ for each of these conjugacy classes. This rank is 6 for elements of class $3C$, 12 for elements of class $3B$ and 0 for elements of class $3A$.
\end{ex}
\begin{ex}\label{ex:A46}
Let us apply this construction to the lattice $A_4(-1)^{6}$, we then have $G(N)=(\mathbb{Z}_{/(5)})^3$ acting on $\Lambda$. The normalizer of this group (inside $Co_1$) is one of the maximal subgroups of $Co_1$, therefore its structure is analyzed in \cite{atlas}. The elements of $G(N)$ fall under three conjugacy classes labeled $5A,5B$ and $5C$. Each conjugacy class has respectively $40,60$ and $24$ representatives inside $G(N)$. Therefore we can compute the rank of the invariant lattice inside $\Lambda$ for each of these conjugacy classes. This rank is 4 for elements of class $5C$, 8 for elements of class $5B$ and 0 for elements of class $5A$.
\end{ex}
\begin{ex}\label{ex:A122}
Let us apply this construction to the lattice $A_{12}(-1)^{2}$, we then have $G(N)=\mathbb{Z}_{/(13)}$. Let $\phi$ be an  automorphism of $\Lambda$ of order 13 generated by a non trivial element $g$ of $G(N)$ on this holy construction. $\phi$ cyclically permutes the extended roots of both copies of $A_{12}$ and therefore has no fixed points in $\Lambda$.
\end{ex}
\begin{ex}\label{ex:p11A212}
Let us look back at \Ref{ex}{A212} and let us analyze an automorphism of order 11: it can be defined by leaving the first copy of $A_2(-1)$ fixed and by cyclically permuting the remaining 11, and the action is extended accordingly to the glue vectors. This automorphism is defined on both $N_{22}$ and $\Lambda$. Let $\varphi$ be this isometry on $A_2^{12}(-1)\otimes\mathbb{Q}$.\\
A direct computation shows $T_{\varphi}N_{22}$ is spanned by
\begin{equation*}
f_1^1,f_1^2,\sum_2^{12}f_1^i,\sum_2^{12}f_2^i,\sum_1^{11} g_j,
\end{equation*}
where $g_j$ are generators for the glue code as in \Ref{tab}{nieme}. Keeping the same notation as before one sees that $S_{\varphi}N_{22}$ has rank 20 and is spanned by
\begin{equation}
(f^k_1-\varphi f^k_1),(f^k_2-\varphi f^k_2),(g_j-\varphi g_j).
\end{equation}

Where $k$ runs from 2 to 12. This vectors satisfy \eqref{holy_leech}, therefore this lattice is contained in $\Lambda$ and, since they are both primitive, $S_{\varphi}(N_{22})=S_{\varphi}(\Lambda)$.
\end{ex}
\begin{ex}\label{ex:p11A124}
A similar computation can be done for $A_1(-1)^{24}$. We use a standard notation where the copies of $A_1(-1)$ are indexed by the set
\begin{equation}\nonumber
\{\infty,0,1,\dots,22\}=\mathbb{P}^1(\mathbb{Z}_{/(23)}).
\end{equation}
Here the isometry $\varphi$ of order 11 is defined by the following permutation on the coordinates:
\begin{equation}\label{p11}
(0)(15\,7\,14\,5\,10\,20\,17\,11\,22\,21\,19)(\infty)(3\,6\,12\,1\,2\,4\,8\,16\,9\,18\,13).
\end{equation}
As before this isometry preserves both $N_{23}$ and $\Lambda$ and the lattice $S_{\varphi}(N_{23})$ is generated by the following vectors:
\begin{equation}
(f^k_1-\varphi f^k_1),(f^l_1-\varphi f^l_1),(g_j-\varphi g_j).
\end{equation}

Here $k$ runs along the indexes contained in the first 11-cycle of \eqref{p11}, $l$ runs  along the second one and $j$ along the generators of the glue code contained in \Ref{tab}{nieme}.\\
Once again all of these generators lie also in $\Lambda$ hence
$S_{\varphi}(N_{23})=S_{\varphi}(\Lambda):=S_{11.K3^{[2]}}$ \index{Lattice, $S_{11.K3^{[2]}}$}.
A direct computation shows that the lattice $S_{\varphi}(N_{23})$ is given by the following quadratic form:

\begin{equation}
\left( 
\tiny{\begin{array}{cccccccccccccccccccc}
 -4 & 1 & -2 & -2 & -1 & 1 & -1 & 1 & -1 & -1 & 2 & 1 & -1 & 2 & -1 & -2 & -2 & 2 & 1 & -1 \\
 1 & -4 & -1 & -1 & -1 & -1 & -1 & 1 & -1 & 2 & -1 & -2 & 2 & 0 & -1 & 0 & 0 & -1 & -2 & 1 \\
 -2 & -1 & -4 & -2 & -1 & -1 & 0 & 1 & 0 & -1 & 1 & 0 & -1 & 2 & -2 & -1 & -1 & 0 & 0 & 1 \\
 -2 & -1 & -2 & -4 & 0 & 0 & -2 & 0 & -1 & 0 & 2 & 1 & 0 & 1 & 0 & 0 & -1 & 1 & 0 & -1 \\
 -1 & -1 & -1 & 0 & -4 & 1 & -1 & 2 & -2 & -1 & 1 & 0 & -1 & 0 & -2 & -2 & 0 & 1 & 1 & -1 \\
 1 & -1 & -1 & 0 & 1 & -4 & 0 & -1 & 0 & 1 & -2 & -1 & 0 & -1 & -1 & 0 & -1 & 0 & -1 & 1 \\
 -1 & -1 & 0 & -2 & -1 & 0 & -4 & 1 & -2 & 1 & 1 & 1 & 0 & -1 & 0 & -1 & 0 & 2 & 0 & -2 \\
 1 & 1 & 1 & 0 & 2 & -1 & 1 & -4 & 0 & 0 & -1 & 1 & 1 & 0 & 2 & 1 & 0 & -1 & 1 & 0 \\
 -1 & -1 & 0 & -1 & -2 & 0 & -2 & 0 & -4 & 0 & 0 & 1 & 1 & 0 & -1 & -2 & 0 & 2 & 0 & -2 \\
 -1 & 2 & -1 & 0 & -1 & 1 & 1 & 0 & 0 & -4 & 1 & 1 & -2 & 1 & 0 & 0 & 1 & 1 & 1 & 0 \\
 2 & -1 & 1 & 2 & 1 & -2 & 1 & -1 & 0 & 1 & -4 & -2 & 2 & -1 & 0 & 0 & 0 & -1 & -2 & 1 \\
 1 & -2 & 0 & 1 & 0 & -1 & 1 & 1 & 1 & 1 & -2 & -4 & 1 & 0 & -1 & 0 & -1 & -1 & -2 & 2 \\
 -1 & 2 & -1 & 0 & -1 & 0 & 0 & 1 & 1 & -2 & 2 & 1 & -4 & 0 & -1 & 0 & 0 & 1 & 2 & 0 \\
 2 & 0 & 2 & 1 & 0 & -1 & -1 & 0 & 0 & 1 & -1 & 0 & 0 & -4 & 1 & 1 & 1 & 0 & 0 & -1 \\
 -1 & -1 & -2 & 0 & -2 & -1 & 0 & 2 & -1 & 0 & 0 & -1 & -1 & 1 & -4 & -2 & -1 & 1 & 0 & 0 \\
 -2 & 0 & -1 & 0 & -2 & 0 & -1 & 1 & -2 & 0 & 0 & 0 & 0 & 1 & -2 & -4 & -2 & 2 & 0 & -1 \\
 -2 & 0 & -1 & -1 & 0 & -1 & 0 & 0 & 0 & 1 & 0 & -1 & 0 & 1 & -1 & -2 & -4 & 1 & 0 & 0 \\
 2 & -1 & 0 & 1 & 1 & 0 & 2 & -1 & 2 & 1 & -1 & -1 & 1 & 0 & 1 & 2 & 1 & -4 & 0 & 2 \\
 1 & -2 & 0 & 0 & 1 & -1 & 0 & 1 & 0 & 1 & -2 & -2 & 2 & 0 & 0 & 0 & 0 & 0 & -4 & 1 \\
 -1 & 1 & 1 & -1 & -1 & 1 & -2 & 0 & -2 & 0 & 1 & 2 & 0 & -1 & 0 & -1 & 0 & 2 & 1 & -4
\end{array}}\right).
\end{equation}

\end{ex}

\section{Prime order Leech automorphisms of Niemeier lattices}\label{sec:prime_nieme}
In this section we give a brief analysis of prime order Leech automorphisms on Niemeier lattices, which will be used for \Ref{thm}{prime_autom_k32}. Our analysis focuses on automorphisms of order 3,5 and 7, while order 11 automorphisms have already been analyzed in \Ref{ex}{p11A124} and \Ref{ex}{p11A212}.\\ 
The fact that the co-invariant lattices in \Ref{ex}{p11A124} and \Ref{ex}{p11A212} are isomorphic is part of a more general behaviour:
\begin{lem}\label{lem:co_invariant_busone}
Let $N\neq \Lambda$ be a Niemeier lattice and let \\$s\in\,Aut(N)/W(N)$ be a Leech isometry. Let moreover $h$ be the Coxeter number of $N$ and let $n$ be the order of $s$. Suppose that $h $ and $n$ are relatively prime. Then $S_s(N)\cong S_{\overline{s}}(\Lambda)$, where $\overline{s}$ is the automorphism of $\Lambda$ obtained by extending the action of $s$ to the holy construction of $\Lambda$ with hole $N$. Moreover this lattice consists of elements of the form $v-\varphi(v)$ for $v\in N$.
\begin{proof}
Let $f_i,\,i\in I$ be a set of root vectors for the holy construction corresponding to $N$ and let $f_i,\,i\in I'\subset I$ be a basis for the Dynkin lattice $R(N)$ contained in $N$. $s$ acts on $I'$ by permuting its elements, and $s$ is univoquely determined by such permutation. Moreover if we let $g_j,\,j\in J$ be the corresponding glue vectors we have that $s$ induces a permutation on $J$, univoquely determined by the permutation on $I'$. This implies that there exists an isometry $\overline{s}$ of $\Lambda$ defined by these permutations, hence by $s$ itself. Let us denote as $s$ both the isometry on $N$ and that on $\Lambda$.\\
Let us first prove that $S_s(R(N))$ is generated by elements of the form $f_i-s(f_i)$. By \Ref{oss}{G_tors} all elements of this form are contained in $S_s(R(N))$, let us suppose on the contrary that there exists an element $v$ in $S_s(R(N))$ which is not of this form. This is equivalent to saying that we can write $v=\sum_{k\in I'/s}\sum_{i\in k} a_if_i$, where $I'/s$ is the set of orbits of $I'$ under the permutation action of $s$ and there exists an orbit $k'$ such that $\sum_{i\in k'}a_i\neq 0$.
 Therefore we have $0\neq\sum_{l<n}s^l(v)\in T_s(R(N))$ since the coefficient of any $f_i,\,i\in k'$ is non zero. But this is absurd since $S_s(R(N))$ is nondegenerate and orthogonal to $T_s(R(N))$.
Now by \Ref{oss}{G_tors} we have that $nv$ is a sum of elements $w-s(w)$ for all $v\in S_s(N)$ or $v\in S_s(\Lambda)$. Moreover we have that $hv$ lies in $N\cap\Lambda$ for all $v$ in $N$ or in $\Lambda$ as in \Ref{oss}{congruenza_busona} and elements of the form $w-s(w)$ are all in $N\cap\Lambda$. This implies that all elements of $S_s(N)$ and of $S_s(\Lambda)$ lie in the intersection $N\cap\Lambda$ and they are equal. 
\end{proof}
\end{lem}
In some cases this lemma can be improved, as in the following:
\begin{ex}\label{ex:p3D46}
Let us consider the lattice $D_4(-1)^6\subset N_{19}$. In this case, for $v\in N_{19}$, we have $2v\in D_4(-1)^6$. Therefore we can modify the proof of \Ref{lem}{co_invariant_busone} so that it works for all automorphisms of prime order $p\neq 2$.
\end{ex}
\begin{oss}\label{oss:busone_EnDn}
If we analyze what happens for Niemeier lattices $N_i$ containing a summand of type $D_n$ or $E_n$ we obtain a refining of \Ref{lem}{co_invariant_busone} in the same spirit of \Ref{ex}{p3D46}: although the coxeter numbers of these components are usually quite large, there is a lower integer $n$ such that $nv\in R(N_i)$, its root lattice, for all $v\in N_i$. Let us see the values of $n$:
\begin{table}[ht]\label{tab:n_for_busone}
\begin{tabular}{|c|c|c|}
\hline 
$i$ & Dynkin diagram of $R(N_i)$ & $n$\\
\hline
1 & $D_{24}$ & 2\\
\hline
2 & $D_{16}E_8$ & 2\\
\hline
3 & $E_8^3$ & 1\\
\hline
5 & $D_{12}^2$ & 2\\
\hline
6 & $A_{17}E_7$ & 6\\
\hline
7 & $D_{10}E_7^2$ & 2\\
\hline
8 & $A_{15}D_9$ & 8\\
\hline
9 & $D_8^3$ & 2\\
\hline
12 & $E_6^4$ & 3\\
\hline
14 & $D_6^4$ & 2\\
\hline
19 & $D_4^6$ & 2\\
\hline
\end{tabular}
\end{table}
\end{oss}
In the general case we cannot give explicit generators of co-invariant lattices, anyhow the following holds:
\begin{lem}\label{lem:co_invariant_busone_plus}
Let $N\neq \Lambda$ be a Niemeier lattice whose Dynkin diagram contains only elements of type $A_n$ and let $G\subset Leech(N)$ be a group of Leech isometries. Then $S_G(N)\cong S_G(\Lambda)$, where the action of $G$ on $\Lambda$ is induced by the holy construction with hole $N$. 
\begin{proof}
It is enough to prove that $S_\varphi(N)\cong S_\varphi(\Lambda)$ for all $\varphi\in G$ of prime order $p$.
Let $f_i$, $i\in I$ be a basis of the (negative) root lattice $R(N)$. As in \Ref{lem}{co_invariant_busone}, $\varphi$ acts as a permutation on $I$, therefore we have $S_\varphi (R(N))=<f_i-\varphi f_i>_{i\in I}$. Analogously the (generalized) lattice $R(N)^\vee$ has an induced action of $\varphi$ and, since this action preserves its standard dual basis $(e_i)_{i\in I}$, $S_\varphi (R(N)^\vee)$ is generated by elements of the form $e_i-\varphi(e_i)$. Let us remind what is this dual basis: for a Dynkin lattice $A_n\subset \mathbb{R}^{n+1}$ we have $e_i=(i,\dots,i,-(h-i),\dots,-(h-i))/(h)$ with $h=n+1$ and $h-i$ coordinates with value $i/(h)$, $i\leq n$. Let us remark moreover that the glue vectors $g_j$, $j\in J$ for the holy construction lie in $R(N)^\vee$.\\
Let now $S'=\{v-\varphi(v)\}_{v\in N}$ and let $w\in S_{\varphi}(N)-S'$. Since we have the obvious inclusions $R(N)\subset N\subset R(N)^\vee$ we can write $w=\sum_{i\in I}a_i(e_i-\varphi(e_i))$. Moreover we can suppose $|a_i|< h$ otherwise we can consider $w-[a_i/h]h(e_i-\varphi(e_i))$. Analogously we can suppose that the $|a_i|$ are minimal, \ie $w$ in $R(N)\otimes \mathbb{Q}$ can be written in the basis $f_i$ with coordinates with absolute value less than $1$. However this implies that in the holy construction with hole $N$ it can be written only in terms of the $g_j$, $j\in J$, therefore $w\in \Lambda$.\\ Analogously let us consider $<hf_i>_{i\in I}\subset \Lambda\subset <\frac{e_i}{h}>_{i\in I}$ by \eqref{holy_leech} and let $S'=\{v-\varphi(v)\}_{v\in\Lambda}$. As above $\varphi$ preserves a basis of $<hf_i>$ and of its dual $<\frac{e_i}{h}>$, therefore the co-invariant lattice is generated by elements of the form $v-\varphi(v)$. Let $w\in S_\varphi(\Lambda)-S'$. We can write $w=\sum_{i\in I}a_i(e_i-\varphi(e_i)/h)$ and we can suppose that $|a_i|<h^2$, otherwise we can consider  $w-[a_i/h^2]h(e_i-\varphi(e_i))$. As above this implies that $w$ can be written only in terms of the $g_j$, therefore $w\in N$. 
\end{proof}
\end{lem}
\begin{ex}\label{ex:p3A83}
Let us apply the holy construction to the lattice $A_8(-1)^3$ and let $\varphi$ be an order 3 permutation of the 3 copies of $A_8(-1)$. With the holy construction with hole $N_{15}$ it induces an automorphism $\varphi$ of $\Lambda$ of order 3 which fixes a rank 8 lattice, therefore by \Ref{ex}{A212} this element is in conjugacy class $3D$ and $S_{\varphi}(\Lambda)\cong S_{3.exo}$ as in \Ref{ex}{p3E83}. Now by \Ref{lem}{co_invariant_busone_plus} we have $S_{\varphi}(N_{15})\cong S_{3.exo}$.
\end{ex}
Moreover many of the computations will be simplified by the use of the following lemma, making good use of what is known for $K3$ surfaces (see \Ref{cap}{k3_case}, in particular \Ref{thm}{nik_k33}):
\begin{lem}\label{lem:k3lat}
Let $(M,G)$ be a Leech couple such that there exists a primitive embedding $M\,\rightarrow\, U^3\,\oplus\,E_8(-1)^2$ 
Then $G\,\cong\,G'$ for some finite algebraic group $G'$ on a $K3$ surface $Y'$. Moreover if $G$ is abelian then $M$ is univoquely determined.
\end{lem} 
We remind that a sufficient condition for the existence of a primitive embedding $M\,\rightarrow\, U^3\,\oplus\,E_8(-1)^2$ is that $rank(M)+l(A_M)\leq 21$.
Let us introduce some notation: suppose we have a negative definite Dynkin lattice $A\subset N\subset A\otimes \mathbb{Q}$ for some Niemeier lattice and suppose $A=C^rD^m$, where $C$ and $D$ are different semisimple component. Then the isometry $\tau_{(a\,\dots\,b)}$ is the permutation $(a\,\dots\,b)$ acting on the $r$ copies of $C$, while $\tau'_{(a\,\dots\,b)}$ does the same thing on the $m$ copies of $D$. If $C$=$A_n(-1)$ we denote $\sigma_i$ the isometry obtained by a central simmetry on the Dynkin diagram of the $i$-th copy of $A_n$. If $C$=$D_4(-1)$ we denote $\gamma_i$ the isometry which rotates by 120° the Dynkin diagram of the $i$-th copy of $D_4$ and we denote $\lambda_{ij}^k$ the isometry exchanging the $i$-th and $j$-th root on the Dynkin diagram of the $k$-th copy of $D_4$.


\subsection{Leech automorphisms of order 3}
Let $N$ be a Niemeier lattice and let $\varphi\subset Leech(N)$ be a Leech automorphism of order 3.
\begin{prop}\label{prop:p3_nieme}
Let $N$, $\varphi$ be as above, then one of the following holds:
\begin{itemize}
\item $S_{\varphi}(N)=K_{12}(-2)$,
\item $S_{\varphi}(N)=W(-1)$,
\item $S_{\varphi}(N)=S_{3.exo}$ as in \Ref{ex}{p3E83}, 
\item $rank(S_{\varphi}(N))=24$.
\end{itemize}
\begin{proof}
By \Ref{lem}{co_invariant_busone} and \Ref{lem}{co_invariant_busone_plus} we need only to analyze what happens in the case of $N=\Lambda$ and $N=N_{i}$ for $i=3,12,18,19$. Moreover by \Ref{oss}{busone_EnDn} we can avoid considering the case of $N_3$ and $N_{19}$. Therefore the proof is a case by case analysis on this 3 Niemeier lattices. For ease of reference in this proof we will denote all Niemeier lattices with the Dynkin diagram they contain, apart for $\Lambda$.
\begin{itemize}
\item[$E_6^4$] There is only one conjugacy class of automorphism of order 3, namely that of $\tau_{(1\,2\,3)}$. Hence $rank(S_{\varphi}(N_{12}))=12$. Moreover $T_{\varphi}(N_{12})$ contains a copy of $E_6(-1)$, hence $rank(S_{\varphi}(N_{12}))+l(A_{S_\varphi})\leq19$ by \Ref{oss}{overl_group}. This implies $S_{\varphi}(N_{12})\cong K_{12}(-2)$ by \Ref{lem}{k3lat}. 
\item[$A_5^4D_4$] There is only one conjugacy class of automorphism of order 3, namely that of $\tau_{(1\,2\,3)}\gamma_1$. Hence $rank(S_{\varphi}(N_{18}))=12$. Moreover $T_{\varphi}(N_{18})$ contains a copy of $A_5(-1)$, hence  $rank(S_{\varphi}(N_{18}))+l(A_{S_\varphi})\leq20$ by \Ref{oss}{overl_group}. This implies $S_{\varphi}(N_{18})\cong K_{12}(-2)$ by \Ref{lem}{k3lat}.
\item[$\Lambda$] There are 4 conjugacy classes of automorphisms of order 3 in $Co_0$, denoted $3A$, $3B$, $3C$ and $3D$ in \cite{atlas}.\\
In \Ref{ex}{A212} we already computed the rank of $T_{\varphi}$ for $\varphi$ in conjugacy classes $3A,3B$ and $3C$. Therefore we know that the element of class $3A$ has no fixed points, moreover it can be used to define the complex Leech lattice (see \cite[chapter 10, section 3.6]{con} and \cite[page 131 and 181]{atlas} or use \Ref{oss}{int_to_cyclo} to endow $\Lambda$ with a $\mathbb{D}_3$-lattice structure of rank 12). The element of class $3D$ has already been computed in \Ref{ex}{p3E83} and it is isometric to $S_{3.exo}$. \Ref{lem}{co_invariant_busone} implies $S_{3B}(\Lambda)\cong S_{\psi_1}(N_{23})$, where $\psi_1$ is an isometry of $N_{23}$ of order 3 given by 6 cycles of length 3 inside $M_{24}$. Therefore $Rank(S_{3B}(\Lambda))=12$, moreover $3B$ fixes on $N_{23}$ a copy of $A_1(-1)^6\oplus A_1(-3)^6$, hence we have $rank(S_{3B}(\Lambda))+l(A_{S_{3B}})\leq18$ by \Ref{oss}{overl_group}. This implies $S_{3B}(\Lambda)\cong K_{12}(-2)$ by \Ref{lem}{k3lat}.  $S_{3C}(\Lambda)$ is given by the remaining case of \Ref{ex}{A212} and 
from that example we see that $Rank(S_{3C}(\Lambda))=18$. 
We also obtain that $S_{3C}\cong W(-1)$ as we argued in \Ref{ex}{slat_rk6}.
\end{itemize}
\end{proof}
\end{prop} 

\subsection{Leech automorphisms of order 5}
Let $N$ be a Niemeier lattice and let $\varphi\subset Leech(N)$ be a Leech automorphism of order 5.
\begin{prop}\label{prop:p5_nieme}
Let $N$, $\varphi$ be as above, then one of the following holds:
\begin{itemize}
\item $S_{\varphi}(N)=S_{5.K3}$ as in \Ref{ex}{S_5K3},
\item $S_{\varphi}(N)\cong S_{5.exo}$ as in \Ref{ex}{slat_5rk4}, 
\item $rank(S_{\varphi}(N))=24$.
\end{itemize}
\begin{proof}
By \Ref{lem}{co_invariant_busone} and \Ref{lem}{co_invariant_busone_plus} we need only to analyze what happens in the case of $N=\Lambda$. 
There are 3 conjugacy classes of order 5 elements in $Co_1$, and they can be obtained using the "holy" construction on $A_4^6$ as in \Ref{ex}{A46}. Keeping the same notation of that example we have that an element of class $5A$ fixes no elements of $\Lambda$. An element of class $5C$ fixes a lattice of rank 4, therefore we have $S_{5C}(\Lambda\cong S_{5.exo}$ by \Ref{ex}{slat_5rk4}. Finally if $\varphi$ is in class $5B$ there is a lattice $F$ of rank 4 and $l(A_F)=1$ inside $T_\varphi{\Lambda}$, therefore $rank(S_{\varphi}(\Lambda))+l(A_{S_\varphi})\leq21$. This implies our claim by \Ref{lem}{k3lat}.
\end{proof}
\end{prop}

\subsection{Leech automorphisms of order 7}
Let $N$ be a Niemeier lattice and let $\varphi\subset Leech(N)$ be a Leech automorphism of order 7.
\begin{prop}\label{prop:p7_nieme}
Let $N$, $\varphi$ be as above, then one of the following holds:
\begin{itemize}
\item $S_{\varphi}(N)=S_{7.K3}$ as in \Ref{ex}{S_7K3},
\item $rank(S_{\varphi}(N))=24$.
\end{itemize}
\begin{proof}
By \Ref{lem}{co_invariant_busone} we need only to analyze what happens in the case of $N=\Lambda$.
There are 2 conjugacy classes $7A,7B$ of elements of order 7 and they can be both obtained by applying the "holy" construction to the lattice $A_6^4$ and considering automorphisms given by the glue code $G(N_{17})$: One class, such as that of the glue code $[1\,2\,1\,6]$, has $rank(S_{\varphi}(\Lambda))=24$. If we take the other class, like that of $[2\,1\,3\,0]$, we obtain $rank(S_{\varphi}(\Lambda))=18$. Moreover in this case $T_{\varphi}(\Lambda)$ contains the lattice $(-6)^6$, hence by \Ref{oss}{overl_group}, we obtain $rank(S_{\varphi}(\Lambda))+l(A_{S_\varphi})\leq21$ therefore our claim holds by \Ref{lem}{k3lat}.

\end{proof}
\end{prop}

\setcounter{prop}{0}
\chapter{Known results on $K3$ surfaces}\label{cap:k3_case}
Since \hk manifolds in dimension 2 are nothing else but $K3$ surfaces it is worthwhile looking at what happens for automorphisms of $K3$ surfaces. This is an active field of research, however there are fundamental results encompassing most of the theory.\\ This chapter is meant as a short survey on the topic and emphasizes the similarities between $K3$ surfaces and manifolds of \ktipo.\\
Notice that the stronger results in this case are mainly due to the stronger statement of the global Torelli, which we recall in \Ref{thm}{global_torelli_k3}, and to the nature of the \kahl cone of a $K3$ surface, which we recall in \Ref{oss}{kahl_cone_k3}. 
For the general theory of K3 surfaces one can see \cite{bhpv}.\\
In this chapter we let $L=U^3\oplus E_8(-1)^2\cong H^2(K3,\mathbb{Z})$ be the $K3$ unimodular lattice.

\begin{thm}\label{thm:global_torelli_k3}
Let $S$ and $S'$ be two $K3$ surfaces and let $\psi\,:H^2(S,\mathbb{Z})\,\rightarrow\,H^2(S',\mathbb{Z})$ be an isometry respecting the Hodge structure and sending a \kahl class $\omega$ in $S$ to a \kahl class $\omega'$ in $S'$.\\ Then there exists a unique isomorphism $f\,:S\,\rightarrow\,S'$ such that $f^*=\psi$.
\end{thm}
\begin{oss}\label{oss:kahl_cone_k3}
Let $S$ be a $K3$ surface, then the \kahl cone $\mathcal{K}_S$ is cut out by -2 divisors, \ie
\begin{equation}
\mathcal{K}_S=\{\alpha\,\in\,\mathcal{C}_S\,|\,(\alpha,e)>0\,\forall\,e\,\,\in\,Pic(S),\,e\,\text{effective},\,e^2=-2\}.
\end{equation}
\end{oss}
\section{On Automorphisms and Cohomology}
In this section we analyze briefly two lattices linked to a group of automorphism on a $K3$ surface, namely the invariant and the co-invariant lattices.
\begin{defn}
Let $X$ be a K3 surface and let $G\subset Aut(X)$. Let $S_G(X)=S_G(H^2(X,\mathbb{Z}))$ \index{Lattice, Co-invariant lattice for a \hk manifold, $S_G(X)$} be the Co-invariant lattice and let $T_G(X)=T_G(H^2(X,\mathbb{Z}))$ \index{Lattice, Invariant lattice for a \hk manifold, $T_G(X)$} be the invariant lattice.
\end{defn}
These lattices share several properties with their higher dimensional analogues (see \Ref{lem}{gaction_gen} and \Ref{lem}{algaction_general} for a comparison).
First of all let us remind that for every finite group $G$ of automorphisms on a $K3$ surface $X$ there is an exact sequence
\begin{equation}
1\,\rightarrow\,G_{0}\,\rightarrow\,G\,\rightarrow\,\mathbb{Z}_{/(n)}\,\rightarrow\,1,
\end{equation}
for some $n$. Here $G_0\subset Aut_{s}(X)$ is a group of symplectic automorphisms, the first map being the natural inclusion and the last map sends an automorphism to its eigenvalue on $H^{2,0}(X)$.

\begin{lem}[Nikulin, \cite{nik1}]\label{lem:k3_algaction}
Let X be a $K3$ surface and let $G\subset Aut(X)$ be a finite group. Let $G_0$ and $n$ be as above. Then the following hold:
\begin{enumerate}
\item If $n>1$, $X$ is algebraic.
\item $g\in G$ acts trivially on $T(X)$ if and only if $g\in G_0$.
\item The representation of $G/G_0$ on $T(X)\otimes\mathbb{Q}$ is isomorphic to a direct sum of irreducible representations each of which has maximal rank $\phi(n)$.
\end{enumerate}
\begin{proof}
\begin{itemize}
\item Suppose that $n>1$. $X/G$ is a normal complex space and let $Y$ be its minimal resolution of singularities. We have that $H^2(Y,\mathbb{C})=H^2(X/G,\mathbb{C})\oplus E$, where $E$ is generated by the exceptional divisors. However $H^2(X/G,\mathbb{C})$ is generated by divisors, hence $h^{2,0}(Y)=0$, \ie $Y$ is algebraic.
\item Let $g\in G_0$, let $\sigma_X$ be a holomorphic 2-form on $X$ and let $\sigma$ be the map from $T(X)$ to $\mathbb{C}$ sending $\alpha$ to $(\sigma_X,\alpha)$. $g$ preserves the intersection form and the Hodge structure, therefore for $x\in T(X)$ we have
\begin{equation}
\sigma(x)=(g\sigma_X,gx)=g\sigma(gx)=\sigma(gx).
\end{equation}
Therefore $x-gx$ lies in $Ker(\sigma)=T(X)\cap S(X)$ and $g$ is the identity on $T(X)/Ker(\sigma)$. This kernel is either 0 or 1 dimensional, in the first case we are done, otherwise let \\$<c>=Ker(\sigma)$. By Riemann-Roch either $c$ or $-c$ is represented by an effective divisor ($c$ has square zero), therefore $gc=c$. This implies that all eigenvalues of $g$ on $T(X)$ are 1 but, since it has finite order, this means that it is actually the identity. Conversely let $g\in G$ act as the identity on $T(X)$, therefore it acts as the identity also on $T(X)\otimes\mathbb{C}$ which contains $\sigma_X$.
\item To prove this we must show that every nontrivial element of $G/G_0$ has no eigenvalue equal to 1 on $T(X)$, so let $g\in G-G_0$, \ie $g\sigma_X=\lambda\sigma_X$, $\lambda\neq1$. Since we now know that $X$ is algebraic the map $\sigma\,:T(X)\,\rightarrow\,\mathbb{C}$ is an embedding and, for all nonzero $x\in T(X)$, $\sigma(x)\neq0$. This implies
\begin{equation}
(\lambda^{-1}\sigma_X,x)=(g^{-1}\sigma_X,x)=(\sigma_X,gx),
\end{equation} 
\ie $gx\neq x$.
\end{itemize}
\end{proof}
\end{lem}
As a consequence of this lemma we have some limitation to the possible order of nonsymplectic automorphisms:
\begin{cor}
Let $X$ and $n$ be as before, then $\phi(n)\leq 21$ and $n\leq66$.
\begin{proof}
By \Ref{lem}{k3_algaction} the representation of $\mathbb{Z}_{(n)}$ over $T(X)\otimes\mathbb{Q}$ is irreducible of maximal rank, \ie of rank $\phi(n)$. Since $X$ is algebraic $rk(T(X))\leq 21$, hence our claim.
\end{proof}
\end{cor}
It is interesting to remark that this bound is attained by an example of I. Dolgachev, moreover most of the intermediate cases also exist, see the recent work of Keum \cite{keum}.\\
Now let us specialize to the symplectic case, our proof of the following lemma differs a little from the original one of Nikulin, but it is almost identical to the higher dimensional case of \Ref{lem}{algaction}:
\begin{lem}\label{lem:k3_algaction2}
Let $X$ be a $K3$ surface and let $G=G_0$ a finite symplectic group of automorphisms of $X$. Then the following hold:
\begin{itemize}
\item $S_G(X)$ is nondegenerate and negative definite.
\item $S_G(X)$ does not contain elements with square -2.
\item $S_G(X)\subset S(X)$ and $T(X)\subset T_G(X)$.
\item $G$ acts trivially on the discriminant group $A_{S_G(X)}$.
\end{itemize}
\begin{proof}
The third assertion is an immediate consequence of \Ref{lem}{k3_algaction} because $G$ acts as the identity on $\sigma_X$ and therefore on all of $T(X)$.\\
To prove that $S_G(X)$ and $T_G(X)$ are nondegenerate let $H^2(X,\mathbb{C})=\oplus_{\rho}U_\rho$ be the decomposition in orthogonal representations of $G$, where $U_{\rho}$ contains all irreducible representations of $G$ of character $\rho$ inside $H^2(X,\mathbb{C})$. Obviously $T_G(X)=U_{Id_{|\mathbb{Z}}}$ and $S_G(X)=H^2(X,\mathbb{Z})\cap \oplus_{\rho\neq Id}U_\rho$, which implies they are orthogonal and of trivial intersection. Hence they are both nondegenerate.\\
Since $G$ is finite there exists a $G$-invariant K\"{a}hler class $\omega_G$ given by $\sum_{g\in G}g\omega$, where $\omega$ is any K\"{a}hler class on $X$.  
Therefore we have: 
\begin{equation}\nonumber
\sigma_X\mathbb{C}\oplus\overline{\sigma}_X\mathbb{C}\oplus\omega_G\mathbb{C}\,\subset\,T_G(X)\otimes\mathbb{C}.
\end{equation}
Hence the lattice $S_G(X)$ is negative definite.\\
To prove the last assertion we use the natural $G$-equivariant isomorphism between $A_{S_G(X)}$ and $A_{T_G(X)}$ given by \Ref{oss}{nik_overlattice2}. On the latter $G$ acts as the identity, therefore it does the same on the former.\\
Let us prove that there are no $-2$ vectors inside $S_G(X)$. Assume on the contrary that we have an element $c\in S_G(X)$ such that $(c,c)=-2$. Then by Riemann-Roch either $c$ or $-c$ is represented by an effective divisor $D$ on $X$. Let $D'=\sum_{g\in G}gD$ which is also an effective divisor on X, but $[D']\in\,S_G(X)\cap T_G(X)=\{0\}$. This implies $D'$ is linearly equivalent to 0, which is impossible.
\end{proof}
\end{lem}
\section{Main results on symplectic automorphisms}
In this section we state the most important results on finite symplectic groups of automorphisms of $K3$ surfaces, let us start with the results of Nikulin:
\begin{thm}\cite[Theorem 4.3]{nik1}\label{thm:nik_k31}
Let $G\subset O(L)$ be a finite group and suppose the following are satisfied:
\begin{itemize}
\item $S_G(L)$ is negative definite.
\item $S_G(L)$ does not contain any element with square $-2$.
\item $rank(S_G(L))\leq 18$.
\end{itemize}
Then there exists a $K3$ surface $S$ and $G'= Aut_s(S)$ such that $G'\cong G$ and $S_G(L)\cong S_{G'}(S)$. 
\end{thm}
Let $\mathcal{G}_{K3}^{alg}$ \index{Set of groups of symplectic automorphisms on some K3 surface, $\mathcal{G}_{K3}^{alg}$} be the set whose elements are isomorphisms classes of finite symplectic subgroups of $Aut(S)$ for some $K3$ surface $S$. Let moreover $\mathcal{G}_{K3}^{alg,ab}$ \index{Set of abelian groups of symplectic automorphisms on some K3 surface, $\mathcal{G}_{K3}^{alg}$} be the subset of $\mathcal{G}_{K3}^{alg}$ obtained by considering only abelian groups.
\begin{thm}\cite[Theorem 4.5]{nik1}\label{thm:abel_k3}
The following assertions hold:
\begin{itemize}
\item $\mathcal{G}_{K3}^{alg}$ is closed under the operations of taking a subgroup of one of its elements or taking a quotient.
\item Let $G\in\mathcal{G}_{K3}^{alg}$. Then every abelian subgroup of $G$ belongs to $\mathcal{G}_{K3}^{alg,ab}$. If $[G,G]$ is its commutator then $G/[G,G]\in\mathcal{G}_{K3}^{alg,ab}$.
\end{itemize}
Moreover
\begin{align}\nonumber 
\mathcal{G}_{K3}^{alg,ab} =&   \{(\mathbb{Z}_{/(2)})^k,\,k\leq\,4;\,\mathbb{Z}_{/(4)},\,\mathbb{Z}_{/(2)}\times\mathbb{Z}_{/(4)};\,(\mathbb{Z}_{/(4)})^2;\,\mathbb{Z}_{/(8)};\,\mathbb{Z}_{/(3)};\\\nonumber  & \mathbb{Z}_{/(5)};\,\mathbb{Z}_{/(7)};\,\mathbb{Z}_{/(6)};\,\mathbb{Z}_{/(2)}\times\mathbb{Z}_{/(6)}\}.
\end{align}
\end{thm}
The following is an important statement on co-invariant lattices which happens to be false for non abelian groups (see \cite{hashi}):
\begin{thm}\cite[Theorem 4.7]{nik1}\label{thm:nik_k33}
Let $G\,\in\mathcal{G}_{K3}^{alg,ab}$ and let $G\subset Aut_s(S)$ and $G\subset Aut_s(S')$ for two $K3$ surfaces $S$ and $S'$. Then the action of $G$ on $H^2(S,\mathbb{Z})$ is isomorphic to the action of $G$ on $H^2(S',\mathbb{Z})$.
\end{thm} 
which is almost equivalent to the following:
\begin{thm}\cite[Theorem 4.8]{nik1}
Let $i\,:G\,\rightarrow\,O(L)$ and $j\,:G\,\rightarrow\,O(L)$ be two embeddings of a finite abelian group $G$ into $O(L)$. Suppose moreover that $G$ satisfies the conditions of \Ref{thm}{nik_k31} for both embeddings. Then $G\,\in\mathcal{G}_{K3}^{alg,ab}$ and there exists $\phi\in O(L)$ such that $i(g)=\phi\circ j(g)\circ \phi^{-1}$ for all $g\in G$.
\end{thm}
Concerning non abelian groups of symplectic automorphisms and the Mathieu group $M_{23}$ there is the beautiful result of Mukai \cite{muk}:
\begin{thm}\label{thm:muk_k3}
Let $S$ be a $K3$ surface and let $G\subset Aut(S)$ be a finite group of symplectic automorphisms. Then $G\subset M_{23}$ and the natural $G$-action as a subset of $M_{23}$ on the set with 24 elements has at least 5 orbits.
\end{thm}
Mukai also classified all elements of $\mathcal{G}_{K3}^{alg}$ without computing the co-invariant lattices, however Kondo's proof \cite{kon} of this result allows an explicit computation of $S_G(S)$ for all $G\in\mathcal{G}_{K3}^{alg}$ and $G\subset Aut(S)$, namely in the following way:
\begin{thm}\label{thm:k3_kondo}
Let $G\in\mathcal{G}_{K3}^{alg}$ and let $S$ be a $K3$ surface on which $G$ acts symplectically. Then $S_G(X)\cong S_{G'}(N)$ where $N$ is one of the 24 negative definite Niemeier lattices and $G'\cong G$ is a subgroup of $Leech(N)$.
\end{thm}
We remark that Kondo's proof allows to exclude only the case $N=\Lambda$ and allows also to impose $G\subset M_{23}$.
\section{Fixed locus of an automorphism and the abelian case}\label{sec:fix_k3}
In this section we compute the fixed locus of a finite abelian symplectic automorphism group on a $K3$ surface $S$ using a simple topological argument.
Let $x\in S$ be a fixed point of $G$, \ie the stabilizer $Stab_G(x)=G_x$ is non trivial. Let us choose local coordinates around $x$ such that the action of $G_x$ is linear. By hypothesis we have that $G_x$ preserves the symplectic form $\sigma_S$ and moreover $\sigma_S(x)\neq 0$, therefore we can set $G_x\subset Sl(2,\mathbb{C})$. It is a well known fact that finite subgroups of $Sl(2,\mathbb{C})$ are cyclic, therefore in suitable local coordinates a generator of $G_x$ can be written as
\begin{equation}
\left( \begin{array}{cc}\xi & 0\\ 0 & \overline{\xi} \end{array}\right),
\end{equation}
where $\xi$ is a primitive root of unity of order $m_x=|Stab_G(x)|$. Since $G$ is abelian the orbit of $x$ consists of $|G/Stab_G(x)|$ points and each of them has stabilizer $Stab_G(x)$. Locally in a neighbourhood of $x$ the quotient $X/G$ has a singularity of type $A_{m_x-1}$ whose resolution yields $m_x-1$ rational curves whose intersection matrix is given by the Dynkin lattice $A_{m_x-1}(-1)$. We wish to remark that the minimal resolution of singularities of $X/G$ is still simply connected and has trivial canonical class, hence it is again a $K3$ surface. Let now $G_i,\,i=1\,\dots\,N$ be all nontrivial cyclic subgroups of $G$ and let $m_i=|G_i|$, $m=|G|$. Let $k_i$ be the number of points with stabilizer $G_i$ and let $k=\sum k_i$. We have that $m/m_i$ divides $k_i$. Let $x_{i,j},\,j=1\,\dots\,k_i$ be the points with stabilizer $G_i$. Let $X'=X-\{x_{i,j}\}_{i,j}$, we know $\chi(X')=24-k$. Let $Y$ be the minimal resolution of $X/G$ and let $Y'$ be $Y$ without the exceptional divisor. Since by removing the resolution of a $A_l$ singularity the Euler characteristic decreases by $l+1$ we have 
\begin{equation}
\chi(Y')=24-\sum_{i=1}^N\frac{k_im_i^2}{m}.
\end{equation}
Notice moreover that $X'/G\cong Y'$ and the restriction of the quotient map $X'\rightarrow X'/G$ is a topological $|G|:1$ cover. Therefore
\begin{equation}\label{char_k3}
\frac{24-k}{m}=24-\sum_{i=1}^N \frac{k_im_i^2}{m}.
\end{equation}
To obtain the result in \Ref{thm}{abel_k3} it is enough to work out a few cases, let us analyze some:
\begin{itemize}
\item Let $G=\mathbb{Z}_{/(p)}$, where $p$ is a prime number. In this case $k$ is just the number of fixed points, each of them has stabilizer $G$. Applying \eqref{char_k3} we obtain $k=\frac{24}{p+1}$, therefore $p\leq11$. Moreover the case $p=11$ can be eliminated since we would have 
\begin{equation}\nonumber
A_{10}(-1)\oplus A_{10}(-1)\subset NS(Y),
\end{equation}
which is clearly impossible since $NS(Y)$ has signature $(3,19)$.
\item Let $G=\mathbb{Z}_{/(p^2)}$, where $p$ is a prime number. Let $t_p$ and $t_{p^2}$ the number of points with stabilizer $\mathbb{Z}_{/(p)}$ and $\mathbb{Z}_{/(p^2)}$ respectively.
Using \eqref{char_k3} and substituting by the previous case $t_p=\frac{24}{p+1}-t_{p^2}$ we have
\begin{equation}\nonumber
24(p^2-1)=(\frac{24}{p+1}-t_{p^2})(p-1)+t_{p^2}(p^4-1).
\end{equation} 
\ie $24=\frac{24}{p+1}+t_{p^2}p^2$. This implies $p\leq 3$, but if $p=3$ we have $t_3=4$ and $t_9=2$ which is impossible since $3$ must divide $t_3$, therefore $p=2$.
\item Let $G=\mathbb{Z}_{/(pq)}$, where $p$ and $q$ are both prime numbers. Let $t_p,t_q$ and $t_{pq}$ be as above. We have $t_p=\frac{24}{p+1}-t_{pq}$ and $t_q=\frac{24}{q+1}-t_{pq}$, therefore \eqref{char_k3} yields
\begin{align*}
24(pq-1)&=(\frac{24}{p+1}-t_{pq})(p^2-1)\\
&+(\frac{24}{q+1}-t_{pq})(q^2-1)+t_{pq}(p^2q^2-1).
\end{align*}
Therefore either $p=3, q=5$ or $p=2,q=3$. However the first case is impossible since we would have 
\begin{equation}\nonumber
A_2(-1)\oplus A_4(-1)\oplus A_{14}(-1)\subset NS(Y).
\end{equation}
\end{itemize}
Proceeding with all possible cases one sees also that every abelian group $G$ can act symplectically in a unique way for what concerns the topology of $X/G$. In the nonsymplectic case the situation is quite the opposite, since already non symplectic involutions form several topologically distinct families, see \cite{nik3} and \cite{nik4}.\\
\section{The nonabelian case}
The aim of this section is to prove \Ref{thm}{k3_kondo}, our proof is slightly different from that given by Kondo \cite{kon} but it is almost identical to the proof of \Ref{thm}{sporadic}.
Let $X$ be a $K3$ surface and let $G\subset Aut(X)$ be a finite symplectic group of automorphisms. Let moreover $H^2(X,\mathbb{Z})\rightarrow L'=U^4\oplus E_8(-1)^2$ be a primitive embedding of the $K3$ lattice such that $(H^2(X,\mathbb{Z}))^\perp$ is just the first hyperbolic summand. Let $R_G(X)=S_G(X)^{\perp_{L'}}$, it has signature $(4,20-rk(S_G(X)))$. By \Ref{lem}{nik_esiste} there exists a negative definite lattice $T'$ of rank $24-rk(S_G(X))$ with the same discriminant group of $R_G(X)$ and the same discriminant form. Since $R_G(X)$ and $S_G(X)$ are unimodular complements they have the same discriminant group and opposite discriminant forms by \Ref{oss}{nik_overlattice2}. Hence $T'\oplus S_G(X)\subset N$, where $N$ is unimodular, even, negative definite and of rank 24. Moreover $G$ acts trivially on $A_{S_G(X)}$, hence the equivariant morphism of \Ref{lem}{nik_overlattice} allows us to extend $G$ to a group of isometries of $N$ such that it acts as the identity on $T'$. As we saw in \Ref{sec}{niemeier} $N$ is a negative definite Niemeier lattice and, since $S_G(X)$ contains no element of square $-2$, $G\subset Leech(N)$.
Up to now we have simply proved that $G\subset Co_0$, we need to eliminate the case $N=\Lambda$. Here comes Kondo's clever trick: it is sufficient to prove that $T'$ can be chosen in such a way that it contains a $-2$ vector. Obviously there are $-2$ vectors inside $R_G(X)$ since it contains a copy of $U$, therefore we let $R'_G(X)$ be the orthogonal complement inside $R_G(X)$ of one of these vectors and we let $T"$ be a negative definite lattice with its discriminant form, group and rank by \Ref{lem}{nik_esiste}. As before $S_G(X)\oplus(-2)\oplus T"\subset N$, but this time $N$ contains a $-2$ vector, hence we can choose $N\neq \Lambda$ and moreover $G\subset M_{24}$. The last step for \Ref{thm}{muk_k3} is now easy: $G\subset M_{23}$ since it fixes at least a $-2$ vector (remember that the action of $M_{24}$ on $A_1^{24}$ is by permutations) and it has at least 5 orbits since the rank of $T'$ is equal to the number of orbits and $rank(S_G(X))\leq 19$.
\section{A few examples}
In this section we give a few interesting examples of $K3$ surfaces with symplectic automorphisms. The interested reader can consult \cite{muk} for a full list of $K3$ surfaces endowed with a maximal symplectic group.
\begin{ex}
Let $X_t$ be the zero locus of the polynomial $\sum x_i^4+tx_0x_1x_2x_3$ in $\mathbb{P}^3$. The group of permutation $S_4$ on the coordinates of $\mathbb{P}^3$ induces automorphisms of $X_t$, however not all of them preserve the symplectic form, but the alternating subgroup $A_4$ does.
\end{ex}  
\begin{ex}
Let $X$ be the complete intersection in $\mathbb{P}^5$ given by $0=\sum_{i=1}^6x_i=\sum_{i=1}^6x_i^2=\sum_{i=1}^6x_i^3$. Again the group of permutations $S_6$ of the coordinates induces automorphisms of $X$, but only its alternating subgroup $A_6$ is symplectic.
\end{ex}
\begin{ex}
Let $X$ be the zero locus of Klein's quartic polynomial $x_0^3x_1+x_1^3x_2+x_2^3x_0+x_3^4$. This is a cyclic $4:1$ cover of $\mathbb{P}^2$ ramified along the curve $C=V(x_0^3x_1+x_1^3x_2+x_2^3x_0)$, it is a classical fact that $Aut(C)=PSL_2(\mathbb{Z}_{/(7)})$ and a direct computation shows that these induce symplectic automorphisms on $X$. 
\end{ex}
\begin{ex}[Sarti, Van Geemen, \cite{vs}]
Let $X$ be a $K3$ surface with an elliptic fibration. Suppose moreover that $X$ has a zero section $\delta$ and a section $\tau$ of order 2. If $X$ is general with respect to these conditions we can suppose
\begin{equation}
X=V(x(x^2+a(t)x+b(t))-y^2),
\end{equation}
where $t\in\mathbb{P}^1$, $a$ has degree 4 and $b$ has degree 8. Moreover $\delta(t)$ is the point at infinity and $\tau(t)=(0,0)$. The section $\tau$ induces an automorphism $\varphi$ of order $2$ on $X$ which is symplectic, moreover if we let $Y$ be the minimal resolution of singularities of $X/\varphi$ we have:
\begin{equation}
Y=V(x(x^2-2a(t)x+(a(t)^2-4b(t)))-y^2).
\end{equation} 
This kind of involutions were first considered by Van Geemen and Sarti and in the literature they are often referred to as Van Geemen-Sarti involutions.
\end{ex}
\setcounter{prop}{0}
\chapter{Examples of Symplectic automorphisms}\label{cap:exoexa}
This chapter is devoted to providing examples of symplectic automorphisms on \hk manifolds. 
There is a natural way to extend an automorphism of a $K3$ surface $S$ to an automorphism of its Douady space $S^{[n]}$ and the same holds for automorphism of abelian surfaces inducing automorphisms on generalized Kummer manifolds. We will call an automorphism \emph{standard} if it can be deformed to an automorphism induced in this way and we will call an automorphism \emph{exotic} otherwise. A precise definition will be given in \Ref{defn}{standard}. We wish to remark that if the fixed locus of a finite order automorphism $\psi$ is topologically different from the fixed locus of a standard automorphism of the same order then $\psi$ is exotic.
Examples \ref{ex:fano3} and \ref{ex:fanoA6} concerning Fano schemes of lines on a cubic fourfold are already present in the literature, we only study in greater detail their group of symplectic automorphisms.

\begin{oss}\label{oss:simple_sympl}
Let $G\subset PGL_6(\mathbb{C})$ and let $[f]$ be a $G$-invariant class of a nonsingular cubic homogeneous polynomial on 6 variables. Let $X=F(V(f))$ be the Fano scheme of lines of the cubic fourfold associated to $f$ and let $Y=VSP(f,10)$ be the variety of sums of powers. Then $G\subset Aut(X)$ and $G\subset Aut(Y)$. Notice however that $Y$ might not be a \hk manifold if $f$ is not general. Moreover if $G$ is simple then $G\subset Aut_s(X)$ and $G\subset Aut_s(Y)$ (whenever this is well defined) using \eqref{exactgroup}. We must remark moreover that the natural polarization of $X$ is $G$-invariant.
\end{oss}

\section{Involutions}
\begin{ex}\label{ex:invol_stand}
Let $S$ be a $K3$ surface and let $\varphi\in Aut_s(X)$ be a symplectic involution. Let $X=S^{[n]}$, $\varphi$ induces a symplectic involution $\varphi^{[n]}$ on it. If we analyze the fixed locus $X^{\varphi^{[n]}}$ we see that, in case $n=2$, it consists of $28$ isolated points and 1 $K3$ surface $Y$. The 28 points are given by pairs $(a,b)$, where $a,b\in S^\varphi$. The fixed $K3$ surface is the closure of the analytic subsets $(x,\varphi(x))$ where $x\in S$ and $\varphi(x)\neq x$. Therefore $Y$ is the resolution of singularities of $X/\varphi$. If $n\geq 3$ the fixed locus consists in a series of points, $K3$ surfaces isomorphic to $Y$ and their Douady schemes. 
\end{ex}

\begin{ex}\label{ex:fano_invol}
This example appeared in a paper of Camere \cite{cam}, let
\begin{equation}
\small{f=x_0^2L_1(x_2,\dots,x_5)+x_0x_1L_2(x_2,\dots,x_5)+x_1^2L_3(x_2,\dots,x_5)+ G(x_2,\dots,x_5)}
\end{equation}
 be a cubic polynomial in six variables, where $L_i$ are linear forms and $G$ is a cubic polynomial. Let $Y=V(f)$ and let $\varphi$ be the involution induced on it by the projectivity sending $[x_0,\dots,x_5]$ in $[-x_0,-x_1,\dots,x_5]$. By \Ref{oss}{residue_fano} this involution induces a symplectic involution $\psi$ on the Fano scheme of lines $F(Y)$. Moreover this family is 12 dimensional and the fixed locus of $\psi$ consists of 28 isolated points and 1 $K3$ surface. 
\end{ex}

\begin{ex}\label{ex:kum_invol}
Let $T$ be an abelian surface and let $X=K_n(T)$ be a generalized Kummer manifold. Then the automorphism $-Id$ of $T$ induces an automorphism of $T^{[n+1]}$ preserving $X$. This involution acts trivially on $H^2(X)$, therefore it is also symplectic.
\end{ex}

\section{Automorphisms of order 3}\label{sec:exo3}
\begin{ex}\label{ex:3_stand}
Let $S$ be a $K3$ surface and let $\varphi\in Aut_s(S)$ be an automorphism of order 3. Let $X=S^{[n]}$ and let $\psi=\varphi^{[n]}$. If $n=2$ then the fixed locus of $\psi$ on $X$ consists of 27 isolated points given by 15 points of the form $(a,b)$ with $a,b\in S^\varphi$, $a\neq b$ and 12 points which are given as the nonreduced points corresponding to the two eigenspaces in $T_aS$ for all fixed points $a\in S^{\varphi}$. If $n=3$ then the fixed locus consists in some isolated points and a $K3$ surface $Y$, given as the closure of the surface $\{(x,\varphi(x),\varphi^2(x)),\,x\in S,\,\varphi(x)\neq x\}$. This is precisely the $K3$ surface obtained from the resolution of singularities of $X/\varphi$. If $n\geq 4$ then the fixed locus consists in a series of isolated points, $K3$ surfaces isomorphic to $Y$ and their Douady schemes. 
\end{ex}
\begin{ex}\label{ex:fano_3stand}
Let $\varphi$ be the projectivity of $\mathbb{P}^5$ sending $[x_0,\dots,x_5]$ to $[\omega x_0,\omega x_1,\overline{\omega} x_2,\overline{\omega} x_3,x_4,x_5]$, where $\omega=e^{\frac{2\pi i}{3}}.$ There exists a nonsingular cubic polynomial $f$ invariant for the induced action of $\varphi$. Then $Y=V(f)$ has an automorphism of order 3 induced by $\varphi$. Moreover by \Ref{oss}{residue_fano} it induces a symplectic automorphism of the Fano variety of lines $F(Y)$ with 27 isolated fixed points. We remark that these examples form a family with $8$ moduli. We will prove in the next chapters that this example is standard.
\end{ex}

\begin{ex}\label{ex:fano_3stand2}
Let $\varphi$ be the projectivity of $\mathbb{P}^5$ sending $[x_0,\dots,x_5]$ to $[\omega x_0,\overline{\omega} x_1,x_2,x_3,x_4,x_5]$, where $\omega=e^{\frac{2\pi i}{3}}.$ There exists a nonsingular cubic polynomial $f$ invariant for the induced action of $\varphi$, it has the form $ax_0^3+bx_1^3+L(x_2,x_3,x_4,x_5)x_0x_1+C(x_2,x_3,x_4,x_5)$, where $a,b\in\mathbb{C}$, $L$ is linear and $C$ is a cubic polynomial. Then $Y=V(f)$ has an automorphism of order 3 induced by $\varphi$. Moreover by \Ref{oss}{residue_fano} it induces a symplectic automorphism of the Fano variety of lines $F(Y)$ with 27 isolated fixed points which are precisely the 27 lines on the cubic surface $V(C)\subset\mathbb{P}^3$. We remark that these examples form a family with $8$ moduli. We will prove in the next chapters that this example is standard.
\end{ex}

\begin{ex}\label{ex:epw_3stand}
Let $\varphi$ be as in \Ref{ex}{fano_3stand} and let $A\in \mathbb{LG}(\Lambda^3\mathbb{C}^6)$ be a $\varphi$-invariant Lagrangian subspace. Then $\varphi$ induces an automorphism of order 3 on the EPW-sextic $Y_A$, moreover this automorphism acts trivially on $K_{Y_A}$ if $Y_A\neq\mathbb{P}^5$. Notice that a $\varphi$-invariant lagrangian subspace is generated by eigenvectors for the action of $\varphi$ on $\Lambda^3\mathbb{C}^6$. Since $A\in\mathbb{LG}^0$ is an open condition we must only find a $\varphi$-invariant lagrangian in $\mathbb{LG}^0$ to obtain a family of Double-EPW Sextics with an order 3 symplectic automorphism. To avoid tedious computation we just refer to \Ref{ex}{epw_11} which satisfies these conditions. We remark that these manifolds form a family with $8$ moduli and we will prove in the next chapters that this example is standard.  
\end{ex}

\begin{ex}\label{ex:epw_3exo}
Let $\varphi$ be the following automorphism of $V=\mathbb{C}^6$:
\begin{equation}
(x_0,x_1,x_2,x_3,x_4,x_5)\,\rightarrow\,(x_0,x_1,x_2,x_3,\omega x_4,\overline{\omega}x_5),
\end{equation}
where $\omega=e^{\frac{2\pi i}{3}}$ And let $V_i$ be the eigenspace with eigenvalue $i$ for $\varphi$. Let $\sigma$ be the symplectic form on $\Lambda^3V$ induced by the standard volume form $vol(e_0\wedge\,\dots\,\wedge e_5)=1$ and let $A$ be a $\varphi$-invariant lagrangian.
Let us remark that $\varphi$ fixes a sextic surface inside $Y_A$.
The action of $\varphi$ preserves the canonical class of $Y_A$, thus it induces an automorphism on $X_A$ which is still trivial on the canonical class. Therefore if there exists a \hk resolution of $X_A$ where the action of $\varphi$ can be extended we would have found an order 3 symplectic automorphism with a fixed surface and an invariant polarization of square $2$. We will prove in \Ref{sec}{exa2} that this is impossible.

\end{ex}

\begin{ex}\label{ex:fano3}
Let $C,D\subset \mathbb{P}^2$ be two elliptic curves given as the zero locus of the cubic polinomials $f$ and $g$ respectively and let $X\subset\mathbb{P}^5$ be the zero locus of the polinomial $f(x_0,x_1,x_2)+g(x_3,x_4,x_5)$. Let $F=F(X)$ be the Fano scheme of lines of $X$. It was first shown in \cite{nami} that $F$ has a symplectic automorphism of order 3 which is not standard, here we show that indeed $F$ has a bigger group of symplectic automorphisms.
 Without loss of generality we can suppose that $f$ and $g$ are in Hesse's normal form, hence the equation of $X$ is
\begin{equation}
x_0^3+x_1^3+x_2^3+x_3^3+x_4^3+x_5^3+\lambda_1x_0x_1x_2+\lambda_2x_3x_4x_5.
\end{equation}
We can moreover consider $C=X\cap\{x_3=x_4=x_5=0\}$ and $D=X\cap\{x_0=x_1=x_2=0\}$.
Let $\psi$ be Namikawa's automorphism, which is defined by 
\begin{equation}
\{x_0,x_1,x_2,x_3,x_4,x_5\}\rightarrow \{\omega x_0,\omega x_1,\omega x_2,x_3,x_4,x_5\},
\end{equation}
where $\omega=e^{\frac{2\pi i}{3}}$. A direct computation shows that the fixed locus on $F$ of the automorphism it induces is isomorphic to an abelian surface (namely $C\times D$), hence it is exotic.\\
But there are several more automorphisms of $F$, we wish to see which automorphisms of $C$ and $D$ extend to automorphisms of $X$ given by projectivities. Let $p_0$ be an inflection point of $C$, we have 
\begin{equation}\label{eq_o1_fano3}
\mathcal{O}_C(1)=\mathcal{O}_C(3p_0).
\end{equation}
Let $f$ be a translation on $C$ given by a point $q$ or order $n$, $f$ is induced by a projectivity of $\mathbb{P}^2$ if and only if $f^*\mathcal{O}_C=\mathcal{O}_C$. By \eqref{eq_o1_fano3} we must have $n=3$, therefore the group of points of order $3$ of $C$ and $D$ induce a group of automorphisms of $X$ isomorphic to $\mathbb{Z}_{/(3)}^4$. Here we list four generators:
\begin{align}
\{x_0,x_1,x_2,x_3,x_4,x_5\}&\rightarrow  \{x_0,\omega x_1,\omega^2 x_2,x_3,x_4,x_5\},\\
\{x_0,x_1,x_2,x_3,x_4,x_5\}&\rightarrow  \{x_2,x_0,x_1,x_3,x_4,x_5\},\\
\{x_0,x_1,x_2,x_3,x_4,x_5\}&\rightarrow  \{x_0,x_1,x_2,x_3,\omega x_4,\omega^2 x_5\},\\
\{x_0,x_1,x_2,x_3,x_4,x_5\}&\rightarrow  \{x_0,x_1,x_2,x_5,x_3,x_4\}.
\end{align} 
Notice that there are several automorphisms of $X$ inducing the same automorphisms on $C$ and $D$ but they are all conjugate through the action of $\psi$. By \Ref{oss}{residue_fano} these automorphisms are all symplectic, furthermore also the involution $\sigma_1\sigma_2$ is, where
\begin{align}
\{x_0,x_1,x_2,x_3,x_4,x_5\} &\stackrel{\sigma_1}{\rightarrow}  \{x_0,x_2,x_1,x_3,x_4,x_5\},\\
\{x_0,x_1,x_2,x_3,x_4,x_5\} &\stackrel{\sigma_2}{\rightarrow}  \{x_0,x_1,x_2,x_3,x_5,x_4\}.
\end{align}
Therefore we have $\mathbb{Z}_{/(3)}^5.\mathbb{Z}_{/(2)}\subset Aut_{s}(F)$. Notice that this examples form a family with 2 moduli.
\end{ex}

\begin{ex}\label{ex:fanoA6}
Let $X\subset\mathbb{P}^5$ be Fermat's cubic, \ie the zero locus of $x_0^3+\dots+x_5^3$. Let $F$ be its Fano scheme of lines. Obviously the permutation group $S_6$ acts on $X$ and, by \Ref{oss}{residue_fano}, it is easy to see that its alternating subgroup $A_6$ induces symplectic automorphisms on $F$. Furthermore Kawatani \cite{kaw} found other symplectic automorphisms $\psi_{i,j,k}$ given by
\begin{align}
\psi_{i,j,k}(x_l)&=x_l &  \text{if }\,\,\,l\notin \{i,j,k\},\\
\psi_{i,j,k}(x_l)&=\omega x_l &  \text{else}.
\end{align}
Here $\omega=e^{\frac{2\pi i}{3}}$. These automorphisms generate a group isomorphic to $\mathbb{Z}_{/(3)}^4$, hence we have $\mathbb{Z}_{/(3)}^4.A_6\subset Aut_{s}(F)$. 
\end{ex}

\section{Automorphisms of order 5}
\begin{ex}\label{ex:5_stand}
Let $S$ be a $K3$ surface and let $\varphi\in Aut_s(S)$ be an automorphism of order 5. Let $X=S^{[n]}$ and let $\psi=\varphi^{[n]}$. If $n\leq4$ then the fixed locus of $\psi$ on $X$ consists only of isolated points (14 if $n=2$). If $n=5$ then the fixed locus consists in some isolated points and a $K3$ surface $Y$, given as the closure of the surface $\{(x,\varphi(x),\dots,\varphi^4(x)),\,x\in S,\,\varphi(x)\neq x\}$. This is precisely the $K3$ surface obtained from the resolution of singularities of $X/\varphi$. If $n\geq 6$ then the fixed locus consists in a series of isolated points, $K3$ surfaces isomorphic to $Y$ and their Douady schemes. 
\end{ex}

\begin{ex}\label{ex:fano5}
Let $\varphi$ be the projectivity of $\mathbb{P}^5$ sending $(e_0,\dots,e_5)$ to $(\omega e_0,\omega^2 e_1,\overline{\omega} e_2,\overline{\omega^2} e_3,e_4,e_5)$, where $\omega=e^{\frac{2\pi i}{5}}.$ Let $f$ be a nonsingular cubic polynomial invariant for the induced action of $\varphi$. Then $Y=V(f)$ has an automorphism of order 5 induced by $\varphi$. Moreover by \Ref{oss}{residue_fano} it induces a symplectic automorphism of the Fano variety of lines $F(Y)$ with 14 isolated fixed points. We remark that these examples form a family with $4$ moduli and we will prove in the next chapters that this example is standard.
\end{ex}

\begin{ex}\label{ex:epw_5}
Let $\varphi$ be as in \Ref{ex}{fano5} and let $A\in \mathbb{LG}(\Lambda^3\mathbb{C}^6)$ be a $\varphi$-invariant Lagrangian subspace. Then $\varphi$ induces an automorphism of order 5 on the EPW-sextic $Y_A$, moreover this automorphism acts trivially on $K_{Y_A}$ if $Y_A\neq\mathbb{P}^5$. We wish to remark that a $\varphi$-invariant lagrangian subspace must be generated by eigenvectors for the action of $\varphi$ on $\Lambda^3\mathbb{C}^6$. These eigenspaces are all $4$ dimensional and decomposable tensors inside them span a $2$ dimensional subspace, therefore it is always possible to choose a lagrangian $A$ without decomposable eigenvectors. As an example let $V_i\subset\mathbb{C}^6$ be the $i$ eigenspace for $\varphi$. Then the 1 eigenspace $(\Lambda^3\mathbb{C}^6)_1$ on $\Lambda^3\mathbb{C}^6$ is $(V_{\omega}\otimes V_{\overline{\omega}}\otimes V_1) \oplus (V_{\omega^2}\otimes V_{\overline{\omega^2}}\otimes V_1)$. Projectivizing we have two lines $\mathbb{P}(V_{\omega^2}\otimes V_{\overline{\omega^2}}\otimes V_1)$ and $\mathbb{P}(V_{\omega}\otimes V_{\overline{\omega}}\otimes V_1)$ of decomposable tensors inside $\mathbb{P} (\Lambda^3\mathbb{C}^6)_1\cong\mathbb{P}^3$. We can therefore choose a line $\mathbb{P}(<a_1,a_2>)$ with empty intersection with them. Since $A\in\mathbb{LG}^0$ is an open condition we must only find a $\varphi$-invariant lagrangian in $\mathbb{LG}^0$ to obtain a family of Double-EPW Sextics with an order 5 symplectic automorphism. To avoid tedious computation we just refer to \Ref{ex}{epw_11} which satisfies these conditions. We will prove in the next chapters that this example is standard.   
\end{ex}

\section{Automorphisms of order 7}
\begin{ex}\label{ex:7_stand}
Let $S$ be a $K3$ surface and let $\varphi\in Aut_s(S)$ be an automorphism of order 7. Let $X=S^{[n]}$ and let $\psi=\varphi^{[n]}$. If $n\leq 6$ then the fixed locus of $\psi$ on $X$ consists only of isolated points (9 if $n=2$). If $n=7$ then the fixed locus consists in some isolated points and a $K3$ surface $Y$, given as the closure of the surface $\{(x,\varphi(x),\dots,\varphi^6(x)),\,x\in S,\,\varphi(x)\neq x\}$. This is precisely the $K3$ surface obtained from the resolution of singularities of $X/\varphi$. If $n\geq 8$ then the fixed locus consists in a series of isolated points, $K3$ surfaces isomorphic to $Y$ and their Douady schemes. 
\end{ex}

\begin{ex}
Let $\varphi$ be the projectivity of $\mathbb{P}^5$ sending $(e_0,\dots,e_5)$ to $(\omega e_0,\omega^2 e_1,\omega^3 e_2,\overline{\omega} e_3,\overline{\omega^2} e_4, \overline{\omega^3} e_5)$, where $\omega=e^{\frac{2\pi i}{5}}.$ There exist nonsingular cubic polynomials invariant for the induced action of $\varphi$, let $f$ be one of them. Then $Y=V(f)$ has an automorphism of order 7 induced by $\varphi$. Moreover by \Ref{oss}{residue_fano} it induces a symplectic automorphism of the Fano variety of lines $F(Y)$ with 9 isolated fixed points. We remark that these examples form a family with $2$ moduli.
\end{ex}

\begin{ex}\label{ex:kum_trivial}
Let $T$ be an abelian surface and let $t\in T$ be a point of order 7. Let $X=K_n(T)$ be the generalized Kummer of $T$ and let $n=7m-1$. Then we can consider the automorphism given by adding $t$ to any point of $T$: it induces an automorphism of $T^{(n+1)}$ which can be used to induce an order $7$ automorphism on $T^{[n+1]}$ preserving $K_n(T)$. Hence we have an order 7 automorphism on $X$. It is a well known fact that this automorphism acts trivially on $H^2(X)$, therefore it is also symplectic. This kind of examples can be given for any $n$, obtaining a symplectic automorphism of order $n+1$.
\end{ex}

\section{Automorphisms of order 11}
 Symplectic automorphisms of order 11 are not present in the case of $K3$ surfaces, let us  state some examples.\\
\begin{ex}\label{ex:fano11}
Let $X\subset\mathbb{P}^5$ be the zero locus of $x_0^3+x_1^2x_5+x_2^2x_4+x_3^2x_2+x_4^2x_1+x_5^2x_3$ and let $F_{Kl}=F(X)$ be the Fano scheme of lines of $X$. The group of symplectic automorphisms of $F_{Kl}$ induced by projectivities on $X$ is particularly interesting. Let $\varphi$ be the automorphism given by $Diag(1,\omega,\omega^3,\omega^4,\omega^5,\omega^9)$, where $\omega=e^{\frac{2\pi i}{11}}$. This automorphism is symplectic on $F_{Kl}$ by \Ref{oss}{residue_fano} and has order 11, therefore it is automatically exotic. Let $KA\subset\mathbb{P}^4$ be the zero locus of $x_0^2x_4+x_1^2x_3+x_2^2x_3+x_3^2x_0+x_4^2x_2$, it is shown in \cite{adl} and \cite{kle} that the group $PSL_2(\mathbb{Z}_{/(11)})=L_2(11)$ acts through projectivities on $KA$. $X$ is a 3 to 1 cover of $\mathbb{P}^4$ ramified along $KA$ through the map $(x_0,\dots,x_5)\,\,\rightarrow\,\,(x_1,\dots,x_5)$, therefore the group $L_2(11)$ acts also on $X$ and on $F_{Kl}$.
Let $(1\,4\,2\,3\,5)$ be a permutation and let $\beta$ be the automorphism it induces on $\mathbb{P}^4$ by permuting the coordinates $[x_1,\dots,x_5]$. $\beta$ leaves $KA$ invariant, hence it induces an automorphism $\beta$ of order 5 on $F_{Kl}$.\\
Using \Ref{oss}{residue_fano} one obtains that $\beta$ is symplectic on $F_{Kl}$. Furthermore a direct computation on the Jacobian ring of $X_{Kl}$ shows that $rk\,(S_\beta(F_{Kl}))=16$.\\
Let us consider the following exact sequence: 
\begin{equation}
1\,\,\rightarrow\,\,H\,\,\rightarrow\,\,\mathbb{Z}_{/(3)}\times L_2(11)\,\,\rightarrow\,\,\mathbb{C}^*,
\end{equation}
where the last map is given by the action of $\mathbb{Z}_{/(3)}\times L_2(11)$ on $H^{2,0}(F_{Kl})$ and $H$ is the quotient of $\mathbb{Z}_{/(3)}\times L_2(11)$ by the image in $\mathbb{C}^*$. Therefore $H$ is a normal subgroup of $L_2(11)$, which is simple. Since $\beta\,\in\,H$ we have $H=L_2(11)$, therefore $L_2(11)$ acts symplectically on $F_{Kl}$.  
Furthermore $F_{Kl}$ has a $\varphi$ invariant polarization of square $6$ and divisibility $2$, hence it must lie in $T^2_{11}$.
\end{ex}


\begin{ex}\label{ex:epw_11}
Let $V=\mathbb{C}^6=<e_0,e_1,e_2,e_3,e_4,e_5>$ and let $vol(e_0\wedge e_1\wedge e_2\wedge e_3\wedge e_4\wedge e_5)=1$ be a volume form inducing a symplectic form on $\Lambda^3\mathbb{C}^6$. Let us consider as in \Ref{ex}{fano11} a representation of $G:=PSL_2(\mathbb{Z}_{/(11)})$ on $V$: it splits as the direct sum of a trivial representation on $<e_0>$ and an irreducible representation of dimension 5 on $V'=<e_1,e_2,e_3,e_4,e_5>$. Let us keep calling $\varphi$ the element of order 11 given by $Diag(1,\omega,\omega^3,\omega^4,\omega^5,\omega^9)$, where $\omega=e^{\frac{2\pi i}{11}}$.
$G$ has elements of order $2,3,5,6$ and $11$. Apart for the elements of order 11 their action on $V$ depends only on their order and can be given in a basis of eigenvectors by $Diag(1,1,1,1,-1,-1)$, $Diag(1,1,\eta^2,\eta^2,\eta^4,\eta^4)$, $Diag(1,1,\nu,\nu^2,\nu^3,\nu^4)$ and $Diag(1,1,\eta,\eta^2,\eta^4,\eta^5)$ respectively, where $\nu^5=\eta^6=1$ are primitive roots of unity (see the character tables in \cite{atlas} for more details). The induced $G$-representation on $\Lambda^6V$ is trivial, hence $G$ acts on the set of lagrangians of $\Lambda^3 V$. We wish to find a $G$-invariant lagrangian $A$ and to prove that there exists a double EPW-sextic $X_A$ which is \hk.
The induced $G$-representation on $\Lambda^3V$ splits as the direct sum of 2 isomorphic (and lagrangian) irreducible representations of dimension 10, given respectively by $F_{e_0}$ and $\Lambda^3V'$. We remark that therefore there is no $G$-invariant element inside $\Lambda^3V$.
Let now $f$ be a $G$-equivariant isomorphism between $F_{e_0}$ and $\Lambda^3V'$, let us denote with $F$ the involution of $\Lambda^3V$ given by
\begin{align}
F(x) &=f(x),\,\,\,\text{if }x\in F_{e_0},\\
F(x) &=f^{-1}(x),\,\,\,\text{if }x\in \Lambda^3V'.
\end{align}
We remark that $vol(x,y)=-vol(F(x),F(y))$. Let $A:=\{(x,f(x)),\,x\in F_{e_0}\}$. Notice that we have
\begin{align*}
(x,f(x))\wedge (y,f(y))&=x\wedge y+f(x)\wedge f(y)+x\wedge f(y)+f(x)\wedge y\\
&=x\wedge f(y)+f(x)\wedge y=x\wedge f(y)+f(x)\wedge y\\
&=x\wedge f(y)-F(f(x))\wedge F(y)=0.
\end{align*}
Therefore $A$ is lagrangian. 
Let us give explicitly the lagrangian $A$:
\begin{align}\label{gen_A}
A= & <e_0\wedge e_2\wedge e_5\,-\, e_2\wedge e_3\wedge e_4,\,e_0\wedge e_3\wedge e_5\,+\, e_1\wedge e_2\wedge e_5,\\ \nonumber
 &\,e_0\wedge e_4\wedge e_5\,+\, e_1\wedge e_3\wedge e_5,\,e_0\wedge e_1\wedge e_2\,-\, e_1\wedge e_4\wedge e_5,\\ \nonumber
 &\, e_0\wedge e_1\wedge e_3\,+\, e_2\wedge e_3\wedge e_5,\,e_2\wedge e_4\wedge e_5\,-\, e_0\wedge e_1\wedge e_4,\\\nonumber
 &\, e_0\wedge e_2\wedge e_3\,-\,e_3\wedge e_4\wedge e_5\,,\,e_1\wedge e_2\wedge e_3\,-\, e_0\wedge e_2\wedge e_4,\\\nonumber
 &\,e_1\wedge e_2\wedge e_4\,-\, e_0\wedge e_3\wedge e_4,\,e_1\wedge e_3\wedge e_4\,+\, e_0\wedge e_1\wedge e_5>
\end{align}
Then $Y_A[3]$ is empty. In fact by \Ref{ssec}{epw} this is equivalent to
\begin{equation}\nonumber
min_{v\in V}\,Rank(\lambda_M(v))\geq 8.
\end{equation}
This condition can be easily checked with some computer algebra and holds true.\\
Moreover with the help of computer algebra we also establish that $Gr(3,6)\cap A=\emptyset$. We would like a more theoretic proof of this fact but we could find none.
Therefore $A\in LG(\Lambda^3V)^0$ and the double cover $X_A$ of $Y_A$ is a manifold of \ktipo. A direct computation shows that $X_A^\varphi$ consists of the 5 points $[e_1],[e_2],[e_3],[e_4],[e_5]$. Moreover $G\subset Aut_s(X_A)$ by \Ref{oss}{simple_sympl}. 

\end{ex}

\section{Alternating groups}
\begin{ex}\label{ex:fano15}
Let $f=x_0^2x_1+x_1^2x_2+x_2^2x_3+x_3^2x_0$ and let $C\subset\mathbb{P}^3$ be its zero locus. Let $g=x_4^3+x_5^3$ and let $Y$ be the zero locus of $f+g$ inside $\mathbb{P}^5$. Let us consider $C\subset Y$ in the obvious way. 
A direct computation shows that $Y$ is nonsingular. Let $\omega=e^{\frac{2\pi i}{15}}$ be a $15$-th primitive root of unity and let $\psi=diag(\omega,\omega^{13},\omega^4,\omega^7,\omega^5,1)$ be a projectivity. Recall that $\psi^5_{|C}=Id_C$.  Let $X=F(Y)$ be the Fano scheme of lines of $Y$ and let $\varphi$ be the automorphism induced by $\psi$ on $X$. Applying \Ref{oss}{residue_fano} one quickly sees that $\varphi$ is symplectic and has order 15. Moreover if we consider the permutation $(0\,1\,2\,3)(4\,5)$ acting on the standard coordinates of $\mathbb{P}^5$ we have that it induces an automorphism $\nu$ of $X$ which has order 4 and, again by \Ref{oss}{residue_fano}, it is symplectic on $X$. 
A natural question would be to determine all possible automorphisms of $X$, let us restrict ourselves to its automorphisms induced by projectivities of $Y$. It is obvious that all automorphisms of the cubic surface $C$ can be extended to automorphisms of $Y$, so let us use the classical work of Segre \cite{segr} and the more recent computations of Hosoh \cite{hos}.
$Aut(C)$ contains an element of order 4 given by $\nu_C$ and one of order 5 given by $\psi^3_C$, therefore its order is a multiple of $20$. Looking at the list of possible automorphism group we see that $Aut(C)=S_5$ is the only possibility and that $C$ is isomorphic to Clebsch's cubic surface $C_l\subset\mathbb{P}^4$ \cite{cleb} given by 
\begin{equation}
\sum_{i=0}^{4}x_i^3=\sum_{i=0}^{4}x_i=0,
\end{equation}  
where $S_5$ acts by permutations on the standard coordinates. However only the elements of $A_5$ induce symplectic automorphisms, but if we compose the others with the permutation sending $[x_0,x_1,x_2,x_3,x_4,x_5]$ to $[x_0,x_1,x_2,x_3,x_5,x_4]$ we have $S_5\subset Aut_s(X)$. We remark that $\psi^5$ commutes with the subgroup $A_5\subset S_5$, we therefore have $(\mathbb{Z}_{/(3)}\times A_5).\mathbb{Z}_{/(2)}\subset Aut_s(X)$.
\end{ex}

\begin{ex}\label{ex:fanoA7}
Let $f$ be the following cubic polynomial:
\begin{equation}
x_0^3+x_1^3+x_2^3+x_3^3+x_4^3+x_5^3-(x_0+x_1+x_2+x_3+x_4+x_5)^3.
\end{equation}
Let $Y=V(f)\subset\mathbb{P}^5$ and let $X=F(Y)$ be its Fano scheme of lines. A direct computation shows that $Y$ is nonsingular and therefore so is $X$, moreover the symmetric group $S_7$ naturally acts as permutations on the set $\{e_0,\dots,e_5,-(e_0+\dots+e_5)\}$ and it preserves $Y$. Another direct computation shows that its alternating subgroup $A_7$ induces symplectic automorphisms of $X$. Considering the natural covering morphism as in \Ref{oss}{cover_fano} we obtain $S_7.\mathbb{Z}_{/(3)}\subset Aut(X)$ in the notation of \cite{atlas}.
\end{ex}

\setcounter{prop}{0}
\chapter{Deformations of automorphisms}\label{cap:deformations}
The aim of this chapter is to analyze the behaviour of automorphisms of a \hk manifold on deformations of the same manifold. The main result contained here is the density of points corresponding to Hilbert squares of points on a K3 and an automorphism induced from the K3 surface inside certain moduli spaces of manifolds of \ktipo with a symplectic automorphism of order 2,3 or 5.\\
First of all let us start with the basic:
\begin{defn}
Let $X$ be a manifold and let $G\subset Aut(X)$. We call a $G$-deformation of $X$ (or a deformation of the couple $(X,G)$) the following data:
\begin{itemize}
\item A flat family $\mathcal{X}\rightarrow B$ and a map $\{0\}\rightarrow B$ such that $\mathcal{X}_0\cong X$.
\item A faithful action of the group $G$ on $\mathcal{X}$ inducing fibrewise faithful actions of $G$. 
\end{itemize}
\end{defn}
From this we give an equivalent to \Ref{defn}{def_equiv}, \ie two couples $(X,G)$ and $(Y,H)$ are deformation equivalent if $(Y,H)$ lies in a $G$-deformation of $X$.
If $G$ is a cyclic group whose action is generated by the automorphism $\varphi$ we will call all $G$-deformations as deformations of the couple $(X,\varphi)$.
The first interesting remark is that, to some extent, all symplectic automorphism groups of a \hk manifold can be deformed:\\
\begin{oss}\label{oss:twistor_deform}
Let $X$ be a \hk manifold such that $G\subset Aut_s(X)$ and $|G|<\infty$. Let $\omega$ be a $G$ invariant \kahl class. Then $TW_{\omega}(X)$ is a $G$ deformation of $X$.
\end{oss}
There is a natural question whenever $X$ is of $K3^{[n]}$-type or of $K_n(T)$-Type which is the following:
\begin{ques}\label{ques:stand_k3}
Let $X$ be a \hk manifold of $K3^{[n]}$-type and let $G\subset Aut_s(X)$. Is it possible to deform the couple $(X,G)$ to a manifold $(S^{[n]},G)$ such that $G\subset Aut_s(S)$ and its action on $S^{[n]}$ is induced by that on $S$?
\end{ques}
The same can be phrased for generalized Kummer manifolds:
\begin{ques}\label{ques:stand_kum}
Let $X$ be a \hk manifold of $K_{n}(T)$-type and let $G\subset Aut_s(X)$. Is it possible to deform the couple $(X,G)$ to a manifold $(K_{n}(T),G)$ such that $G\subset Aut_s(T)$ and its action on $K_{n}(T)$ is induced by that on $T$?
\end{ques}
This answer is false in general, see \Ref{ex}{fano11} for a counter-example. However there are many cases where these questions hold true, let us give the following:
\begin{defn}\label{defn:standard}
Let $(X,G)$ be a couple consisting in a \hk manifold and a finite group $G$ such that \Ref{ques}{stand_k3} or \Ref{ques}{stand_kum} holds true. Then we call $(X,G)$ a standard couple and $G$ a standard automorphism group.
\end{defn}
We call $(X,G)$ exotic otherwise. These definition is equivalent to \Ref{defn}{standard_intro}. In \Ref{cap}{exoexa} we have given several examples of manifolds endowed with an exotic automorphism group.
\section{The universal deformation of $(X,G)$}
In this section we will give a representative to the functor of small deformations of the couple $(X,G)$, where $X$ is a \hk manifold and $G\subset Aut_s(X)$ a finite group. Our construction uses the universal family of deformations $\mathcal{X}\rightarrow Def(X)$.

Let us choose a small ball $U$ representing $Def(X)$ whose tangent space at the origin is given by $H^1(\mathcal{T}_X)$.\\ Let us extend locally the action of $G$ on $U$ using its natural action on $H^1(\mathcal{T}_X)$. Let us shrink $U$ if needed, therefore we can suppose $G(U)=U$.\\ The action of $G$ on $X$ and on $U$ extends to an automorphism of the versal deformation family $\mathcal{X}\,\rightarrow\, U$ as follows:
\begin{equation}\nonumber
\begin{array}{ccc}
G\times\mathcal{X} & \stackrel{M}{\longrightarrow} & \mathcal{X}\\
\downarrow &  & \downarrow\\
G\times U & \stackrel{M_U}{\longrightarrow} & U
\end{array}.
\end{equation}
Moreover $M_U$ induces an action of $G$ on $\mathcal{X}$ which yields fibrewise isomorphisms between $\mathcal{X}_t$ and $\mathcal{X}_{g(t)}$ for all $g\in G$. The differential of $g$ at $0$ is given by the action of $g$ on $H^1(\mathcal{T}_X)$.
On the other hand $U^G$ is smooth since $G$ is linearizable and hence 
\begin{equation}\nonumber
dim (U^G)=dim(H^1(\mathcal{T}_X)^G)=dim(H^2(X)^G)-2,
\end{equation}
which is always positive by \Ref{oss}{twistor_deform}. We wish to obtain a deformation of the couple $(X,G)$, hence we need to restrict to $U^G$ to get a fibrewise action of $G$. Therefore we obtain the following diagram:
\begin{equation}\label{deform}
\begin{array}{cccc}
G\times\mathcal{Y}= &G\times\mathcal{X}_{|U^G} & \stackrel{M}{\longrightarrow} & \mathcal{X}\\
&\downarrow &  & \downarrow\\
&G\times U^G & \stackrel{M_U}{\longrightarrow} & U,
\end{array}
\end{equation}
where $\mathcal{Y}\,\rightarrow\, U^G$ represents the functor of deformations of the couple $(X,G)$, \ie all small deformations of this couple must embed in $\mathcal{Y}\,\rightarrow\, U^G$. We remind that the action $G\times U^G$ is trivial.
The automorphisms $g_{t}$ are given by $M_{|X_t}(g,-)$.\\
It is obvious that this deformation space is "maximal" in some sense, let us make this more precise using the period map.

\begin{defn}
Given a finite group $G$ acting faithfully on a lattice $M$, we call $\Omega_G$ \index{Group-invariant periods, $\Omega_G$}the set of points $(X,f)$ in the period domain $\Omega_N$ such that $f(\sigma_X)\,\in\,T_G(M)$. 
\end{defn}
\begin{defn}
Given $(X,f)$ with a finite group $G$ acting faithfully on it via symplectic bimeromorphisms we call the following a maximal family of deformations of $(X,G_{Bir})$\index{Couple of \hk manifold and birational morphisms, $(X,G_{Bir})$}
\begin{equation}\nonumber
\begin{array}{ccc}
X & \stackrel{i}{\longrightarrow} & \mathcal{X}_U\\
\downarrow &  & \downarrow\\
\{0\} & \stackrel{i}{\longrightarrow} & U,
\end{array}
\end{equation}
where the family $\mathcal{X}$ over $U$ is endowed with a fibrewise faithful bimeromorphic action of $G$ and the period map $\mathcal{P}$, given a compatible marking, sends surjectively a neighbourhood of $0\in\,U$ inside a neighbourhood of $\mathcal{P}(X,f)\cap \Omega_{G}$.\\
We give the same definition for maximal families $(X,G_{Aut})$\index{Couple of \hk manifold and automorphisms, $(X,G_{Aut})$} or $(X,G_{Hod})$\index{Couple of \hk manifold and Hodge isometries, $(X,G_{Hod})$} having $G$ acting as symplectic automorphisms or Hodge isometries on $H^2(X,\mathbb{Z})$ respectively.\\
Notice that the family $\mathcal{Y}\,\rightarrow\,U^G$ we stated before is a maximal family for the couple $(X,\varphi)$.\\

\end{defn}

\begin{oss}\label{oss:generic_point}
We remark that the set 
\begin{equation}
\Omega'_G=\bigcup_{v\in T_G(L)}\{x\in \Omega_G\,:\,(x,v)=0\}
\end{equation}\index{Group-invariant periods, very general subset, $\Omega'_G$}
is the union of countable codimension 1 subsets and consists of Hodge structures on marked varieties $(X,f)$ over $\Omega_G$ such that the inclusion $f(T(X))\hookrightarrow T_G(L)$ is proper. Moreover outside this set $T(X)$ is irreducible. 
\end{oss}

\begin{prop}\label{prop:fixed_abel}
Let $X$ be a manifold of \ktipo and let $\varphi\subset Aut(X)$ be a symplectic automorphism of finite order. Suppose $\varphi$ fixes at least one complex torus $T$. Then $T_{\varphi}(X)$ has rank at most 6.
\begin{proof}
Suppose on the contrary that $T_{\varphi}(X)$ has rank $\geq7$. Let us consider small deformations of the couple $(X,\varphi)$ over a representative $U$ of $Def(X)$ given by
\begin{equation}
\begin{array}{ccc}
\mathcal{X}_{|U^G} & \stackrel{\Phi}{\longrightarrow} & \mathcal{X}\\
\downarrow &  & \downarrow\\
U^G & \stackrel{M_U(Id,-)}{\longrightarrow} & U,
\end{array}
\end{equation}
as shown in \eqref{deform}, where $G=\langle\varphi\rangle$. We let $\sigma_t$ be the symplectic form on $\mathcal{X}_t$.\\ We remark that, by linear algebra, the fixed locus $X^{\varphi}$ is smooth and consists only of symplectic varieties since the symplectic form $\sigma$ restricts to a nonzero symplectic form on all connected components of $X^\varphi$. Moreover it is stable for small deformations of the couple $(X,\varphi)$, \ie the fixed loci $\mathcal{X}^{\Phi}$ is a small deformation of the fixed locus $X^{\varphi}$.
Therefore we have a well defined map of integral Hodge structures $H^2(\mathcal{X}_t,\mathbb{C})^{\Phi_t}\,\rightarrow\, H^2(T_t,\mathbb{C})$ sending a class on $H^2(\mathcal{X}_t)$ to its restriction to $T_t$, where $T_t$ is a small deformation of $T$ fixed by $\Phi_t$ (\ie is a component of the fibre over $t$ of $\mathcal{X}^{\Phi}$). Since $\Phi_t(\sigma_t)=\sigma_t$ and $\sigma_{t|T_t}\neq 0$ this map is not the zero map and, being a map of Hodge structures, its kernel is again a Hodge structure.\\ 
Given a marking $F$ over $\mathcal{X}$ we have that $(\mathcal{X},F)$ is a maximal family of deformations of the couple $(X,\varphi)$. Let $V=\{\mathcal{P}(\mathcal{X}_t,F_t),\,t\in U^G\}\,\subset\,\Omega_\varphi$, by \Ref{oss}{generic_point} there exists $u\in V\backslash \Omega'_\varphi$ and this period corresponds to a marked manifold $(\mathcal{X}_t,F_t)$ such that $T(\mathcal{X}_{t})=T_{\Phi_t}(\mathcal{X}_{t})$, \ie this Hodge structure is irreducible. Therefore we have that the map $H^2(\mathcal{X}_t,\mathbb{C})^{\Phi_t}\rightarrow H^2(T_t,\mathbb{C})$ is an injection. But this is absurd if $T_\varphi(X)$ has rank greater than 6 since $H^2(T_t)$ has dimension 6.
\end{proof}

\end{prop}

\section{Standard automorphisms}
In this section we will prove that symplectic automorphisms $\varphi$ of order $2,3$ and $5$ are standard if some conditions on the lattices $S_\varphi(X)$ are met, however \Ref{thm}{prime_autom_k32} allows us to improve considerably the statement of \Ref{thm}{standard_morph}. 
The technique of the proof is the same in all these cases, we will prove that given any couple $(X,\varphi)$ with the above properties there exists a sequence of couples $(S_n^{[2]},\psi_n^{[2]})$ converging to $(X,\varphi)$ in an appropriate moduli space.\\ 
To prove this result we will need a series of technical lemmas. Let us first fix some notation:
Let $L$ be as in \eqref{latticeK3n}, $M_2$ and $L'$ be as in \Ref{ex}{latticeLam}. Let $M_3$ be as in \Ref{ex}{discr_autom3},  
 $K_{12}(-2)$ be as in \Ref{ex}{cox-todd}, $M_5$ be as in \Ref{ex}{discr_autom5} 
  and $S_{5.K3}$ be as in \Ref{ex}{S_5K3}.

\begin{lem}\label{lem:marking_2lattice}
Let $M,R\,\subset L$ such that $M\cong R\cong E_8(-2)$. Then there exists $f\,\in\,O(L)$ such that $f(M)=R$. 
\begin{proof}
By \Ref{ex}{discrinvol} we know the discriminant form and group of $E_8(-2)$. Therefore we can apply \Ref{lem}{nik_immerge}, obtaining that embeddings of $E_8(-2)$ into $L$ are given by quintuples $(H,H',\gamma,K,\gamma_K)$. Moreover two such embeddings $(H,H',\gamma,K,\gamma_K)$ and $(N,N',\gamma',K',\gamma'_{K'})$ are conjugate if and only if we have $H$ conjugate to $N$ through an automorphism of $(\mathbb{Z}_{/(2)})^8$ sending $\gamma$ into $\gamma'$. In our case the computations are particularly simple: due to the values of $q_{E_8(-2)}$ (all elements have square $0$ or $1$) and $q_L$ (all non zero elements have square $\frac{1}{2}$) the only possible choices of $H$ and $H'$ are given by the one element group and so we obtain our claim.\\
Moreover this implies that we can always choose a marking of $(X,\varphi)$, where $\varphi$ is a symplectic involution such that the induced action of $\varphi$ on $L$ is given by leaving $(-2)\oplus U^3$ invariant and exchanging the two remaining $E_8(-1)$, so that $S_{\varphi}$ is given by the differences $a-\varphi(a)$ for $a\in\,E_8(-1)$.    
\end{proof}
\end{lem}
\begin{lem}\label{lem:marking_3lattice}
Let $M,R\,\subset L$ such that $M\cong R\cong K_{12}(-2)$. Then there exists $f\,\in\,O(L)$ such that $f(M)=R$.
\begin{proof}
By \Ref{ex}{cox-todd} we know the discriminant form and group of $K_{12}(-2)$. Therefore we can apply \Ref{lem}{nik_immerge}, obtaining that embeddings of $K_{12}(-2)$ into $L$ are given by quintuples $(H,H',\gamma,K,\gamma_K)$. Moreover two such embeddings $(H,H',\gamma,K,\gamma_K)$ and $(N,N',\gamma',K',\gamma'_{K'})$ are conjugate if and only if we have $H$ conjugate to $N$ through an automorphism of $(\mathbb{Z}_{/(3)})^6$ sending $\gamma$ into $\gamma'$. In our case the computations are particularly simple: due to the structure of $A_{K_{12}(-2)}$ (all non zero elements have order $3$) and $A_L$ (all non zero elements have order 2) the only possible choices of $H$ and $H'$ are given by the one element group and so we obtain our claim.\\

\end{proof}
\end{lem}

\begin{lem}\label{lem:marking_5lattice}
Let $M,R\,\subset L$ such that $M\cong R\cong S_{5.K3}$. Then there exists $f\,\in\,O(L)$ such that $f(M)=R$.
\begin{proof}
The proof goes as in \Ref{lem}{marking_3lattice}, this time uniqueness up to isometry is a consequence of $A_{S_{5.K3}}=\mathbb{Z}_{/(5)}^4$.
\end{proof}
\end{lem}

Now we define some moduli spaces, namely $\mathcal{M}_{2}=\mathcal{M}_{E_8(-2),L}$\index{Moduli of $E_8(-2)$-polarized \hk fourfolds, $\mathcal{M}_{2}$} as in \Ref{defn}{hk_gen_polar} and analogously $\mathcal{M}_{3}=\mathcal{M}_{K_{12}(-2),L}$\index{Moduli of $K_{12}(-2)$-polarized \hk fourfolds, $\mathcal{M}_{3}$} and
$\mathcal{M}_{5}=\mathcal{M}_{S_{5.K3},L}$\index{Moduli of $S_{5.K3}$-polarized \hk fourfolds, $\mathcal{M}_{5}$}. Notice that lemmas \ref{lem:marking_2lattice}, \ref{lem:marking_3lattice} and \ref{lem:marking_5lattice} imply that these are univoquely determined.

\begin{defn}
Let $\Omega_i=\mathcal{P}(\mathcal{M}_i),\,\, i=2,3,5$\index{Periods of \hk manifolds in $\mathcal{M}_i$, $\Omega_{i}$} and furthermore let $\Omega_{v,i}$\index{Periods of \hk manifolds in $\mathcal{M}_i$ orthogonal to $v$, $\Omega_{v,i}$} denote the set of $\omega\,\in\Omega_i$ such that $(v,\omega)=0$ for $v\in L$. 
\end{defn}

Let $M_2$, $M_3$ and $M_5$ be as before, there is a sublattice $M_{i,0}$ of $L$ isomorphic to $M_i$, $i=2,3,5$ given by $f(T_{\varphi^{[2]}}(S^{[2]}))$ where $(S^{[2]},f)$ is a marked Hyperk\"{a}hler manifold and $\varphi$ is a symplectic automorphism of order $i$ on $S$. 
 Moreover, by \Ref{lem}{marking_2lattice}, \Ref{lem}{marking_3lattice} and \Ref{lem}{marking_5lattice}, all such lattices are conjugate through an isometry of $L$, hence without loss of generality we fix $M_{5,0},M_{3,0},M_{2,0}\subset L$,  $M_i\cong M_{i,0}$, $i=2,3,5$ 
  and we can impose
\begin{equation}\nonumber
\mathcal{P}(X,f)\in\mathbb{P}(M_{i,0}\otimes\mathbb{C}),\,\,\,i=2,3,5
\end{equation}
for all couples $(X,\varphi)$ where $\varphi$ is a symplectic automorphism of order $i$ and $f$ is an appropriate marking. 

\begin{lem}\label{lem:emblemma_2}
Let $0\neq w\in M_2$ be a primitive isotropic vector, then there exist a sublattice $w\in T\subset M_2$ and a $(-2)$ vector $p$ such that:
\begin{itemize}
\item $p$ is 2-divisible in $M_2$,
\item $q_{M_2|T}$ is nondegenerate,
\item $R:=T^{\perp_{M_2}}\cong U\oplus <p> \oplus R'$ for some lattice $R'$.
\end{itemize} 

\begin{proof}
Since $M_2=U^2\oplus (U\oplus E_8(-2)\oplus (-2))$ we can apply \Ref{lem}{ghs_orbit}. Therefore we can analyze up to isometry all isotropic vectors inside $M_2$ knowing only their divisibility $m$ (\ie $(w,M_2)=m\mathbb{Z}$) and their image $[\frac{w}{m}]$ in $A_{M_2}$. Let us give a basis of $M_2$ as follows:
\begin{equation}
\{e_1,f_1,e_2,f_2,e_3,f_3,a_1,a_2,a_3,a_4,a_5,a_6,a_7,a_8,t\},
\end{equation}
where $\{e_i,f_i\}$ is a standard basis of $U$, $\{a_1,\dots,a_8\}$ is a standard basis of $E_8(-2)$ and $t$ is a generator of the lattice $(-2)$.\\
The first key remark is that since $A_{M_2}$ is of 2-torsion $m$ can either be $1$ or $2$.
Therefore if $m=1$ we have that $\frac{w}{m}$ lies in $M_2$, which implies $[\frac{w}{m}]=0$ in $A_{M_2}$. Thus by \Ref{lem}{ghs_orbit} there exists an isometry $g$ of $M_2$ sending $w$ to $e_1$. To obtain our claim we let $T=g^{-1}(<e_1,f_1>)$, $p=g^{-1}(t)$ and $R=g^{-1}(<e_2,f_2,e_3,f_3,a_1,a_2,a_3,a_4,a_5,a_6,a_7,a_8,t>)$.\\
If $m=2$ we have that $\frac{w}{2}$ is a square zero element of $M_2^\vee$, \ie $[\frac{w}{2}]$ has square zero in $A_{M_2}$. Looking at \Ref{ex}{discrinvol} it is easy to see that square zero elements must lie in $A_{E_8(-2)}\subset A_{M_2}$ and they are given by $[\frac{v}{2}]$ where $v$ is a primitive vector of square $c\equiv 0\,\,\,mod\,\,8$ inside $E_8(-2)$. 
Therefore by \Ref{lem}{ghs_orbit} there exists an isometry $g$ of $M_2$ sending $w$ to $r=2e_1+\frac{c}{4}f_1+v$. Thus we set $T=g^{-1}(<r,f_1>)$, $p=g^{-1}(t)$, $K=v^{\perp_{E_8(-2)}}$ and $R=g^{-1}(<e_2,f_2,e_3,f_3,K,t>)$.   
\end{proof}
\end{lem}

\begin{lem}\label{lem:emblemma_3}
Let $0\neq w\in M_3$ be a primitive isotropic vector, then there exist a sublattice $w\in T\subset M_3$ and a $(-2)$ vector $p$ such that:
\begin{itemize}
\item $p$ is 2-divisible in $M_3$,
\item $q_{M_3|T}$ is nondegenerate,
\item $R:=T^{\perp_{M_3}}\cong <p> \oplus R'$ for some lattice $R'$.
\end{itemize} 
\begin{proof}
We can apply \Ref{lem}{nik_spezza} to obtain that $M_3\cong U^2\oplus P$ for some lattice $P$. Hence we can apply \Ref{lem}{ghs_orbit}. Therefore we can analyze up to isometry all isotropic vectors inside $M_3$ knowing only their divisibility $m$ (\ie $(w,M)=m\mathbb{Z}$) and their image $[\frac{w}{m}]$ in $A_{M_3}$. Let us give a basis of $M_3$ as follows:
\begin{equation}
\{e,f,e_1,f_1,e_2,f_2,a_1,a_2,b_1,b_2,t\},
\end{equation}
where $\{e,f\}$ is a standard basis of $U$, $\{e_i,f_i\}$ is a standard basis of $U(3)$, $\{a_1,a_2\}$ and $\{b_1,b_2\}$  are a standard basis of $A_2(-1)$ and $t$ is a generator of the lattice $(-2)$.\\
The first key remark is that since $A_{M_3}$ is of 6-torsion $m$ can either be $1$, $2$, $3$ or $6$. Moreover a direct computation shows $m\neq2$ and $m\neq 6$ due to the values the discriminant form (see \Ref{ex}{discr_autom3}). 
Therefore if $m=1$ we have that $\frac{w}{m}$ lies in $M_3$, which implies $[\frac{w}{m}]=0$ in $A_{M_3}$. Thus by \Ref{lem}{ghs_orbit} there exists an isometry $g$ of $M_3$ sending $w$ to $e$. To obtain our claim we set $T=<g^{-1}(e),g^{-1}(f)>$, $p=g^{-1}(t)$ and $R'=g^{-1}(<e_1,\dots,b_2>)$.\\
Now suppose $m=3$: by \Ref{lem}{ghs_orbit} we have that there exists an isometry $g$ sending $w$ inside $D=<e_1,f_1,e_2,f_2,a_1,a_2,b_1,b_2>$, therefore we can set $p=g^{-1}(t)$, $T=<w,g^{-1}(e),g^{-1}(f)>$ and $R'=g^{-1}(g(w)^{\perp_{D}})$.

\end{proof}
\end{lem}

\begin{lem}\label{lem:emblemma_5}
Let $0\neq w\in M_5$ be a primitive isotropic vector, then there exist a $(-2)$ vector $p$ such that:
\begin{itemize}
\item $p$ is 2-divisible in $M_5$,
\item $(w,p)=0$.
\end{itemize} 
\begin{proof}
Let $p$ be an element of divisibility 2 and square $-2$. By \Ref{lem}{nik_immerge} these elements form a single orbit under the action of $O(M_5)$. By \Ref{lem}{ghs_orbit_gen} we need only to prove that $[w/div(w)]$ is orthogonal to $[p/2]$ in $A_{M_5}$. By \Ref{ex}{discr_autom5} we know the discriminant form and group of $M_5$. Since $w^2=0$ we have also that $[w/div(w)]^2\equiv 0\,mod\,2$. Let $e_1,f_1$ and $e_2,f_2$ be two standard generators respectively of the first and second copy of $U(5)\subset M_5$. Then $A_{M_5}$ is generated by $x_1=[e_1/5],y_1=[f_1/5],x_2=[e_2/5],y_2=[f_2/5]$ and $z=[p/2]$ with the following bilinear form:
\begin{equation}\nonumber
\left( \begin{array}{ccccc} 
0 & \frac{1}{5} & 0 & 0 &0\\
\frac{1}{5}& 0&0&0&0\\
0&0&0&\frac{1}{5} &0\\
0&0&\frac{1}{5}&0&0\\
0&0&0&0&-\frac{1}{2}
\end{array}\right).
\end{equation}
Let $[w/m]=a_1x_1+b_1y_1+a_2x_2+b_2y_2+cz$, $a_i,b_i\in\mathbb{Z}_{/(5)}$ for $i=1,2$ and $c\in\mathbb{Z}_{/(2)}$. A direct computation shows that $w^2=0$ implies $c=0$.
\end{proof}
\end{lem}

\begin{lem}\label{lem:denselemma}
Let $0\neq w_0\in M_{i,0},\,\,i=2,3,5$ be a primitive vector of square 0.\\ There exists an element $q\in M_{i,0}$ of square $-2$ and divisibility $2$ in $L$ such that $w_0\perp q$.
\begin{proof}
We keep the same notation as before and we fix an isometry $\eta_i\,: M_{i,0}\,\rightarrow\,M_i$.\\ First of all let us prove that there exists such a $q$ with divisibility $2$ inside $M_{i,0}$. The proof goes identically for all $i$, let us do it for $i=2$.
Let $w=\eta_2(w_0)$, since it satisfies the hypothesis of \Ref{lem}{emblemma_2} we have an element $p$ orthogonal to $w$, where $p$ is a 2-divisible $(-2)$ vector. hence we can impose $q=\eta_2^{-1}(p)$. Now we need to prove that $div(q)=2$ also inside $L$, \ie that $q\oplus q^{\perp_L}=L$. We know that $p^{\perp_{M_2}}\cong U^3\oplus E_8(-2)$ hence $q^{\perp_L}$ is an overlattice of $U^3\oplus E_8(-2)^2$ which, by \Ref{lem}{evenembed}, implies $\eta^{-1}(p)$ is 2-divisible in $L$. 
\end{proof}  
\end{lem}

\begin{defn}
Let $\mathcal{P}_{i,exc}=\{f\,\in\,M_{i,0}\,:\,f^2=-2\,\,,\,(f,L)=2\mathbb{Z}\}$\index{Exceptional classes, Set of exceptional primitive classes inside $\Omega_{i}$, $\mathcal{P}_{i,exc}$} be the set of exceptional primitive classes inside $M_{i,0}$.
\end{defn}
Notice that $\mathcal{P}^{-1}(v)$ contains the Hilbert square of a $K3$ surface for all periods $v$ orthogonal to some element of $\mathcal{P}_{i,exc}$.
\begin{lem}\label{lem:denselemma_235}
$\cup_{v\in\mathcal{P}_{i,exc}}\Omega_{v,i}$ is dense in $\Omega_i$.
\begin{proof}
The proof goes the same for $i=2,3,5$.
It is enough to prove that $\cup_{v\in\mathcal{P}_{i,exc}}\Omega_{v,i}$ is dense in $\Omega_i\cap\mathbb{P}(M_{i,0}\otimes\mathbb{C})$ by lemmas \ref{lem:marking_2lattice}, \ref{lem:marking_3lattice} and \ref{lem:marking_5lattice}.\\ Let $Q_{M_{i,0}}$ be the subset of isotropic vectors inside $\mathbb{P}(M_{i,0}\otimes\mathbb{C})$. Let $Q_{M_{i,0}}(\mathbb{R})$ and $Q_{M_{i,0}}(\mathbb{Q})$ be the subsets of isotropic vectors spanned by real (respectively rational) isotropic vectors. Let $\omega$ be in $\Omega_i\cap\mathbb{P}(M_{i,0}\otimes\mathbb{C})$, we have $\omega^{\perp_{M_{i,0}}}\,\cap\,Q_{M_{i,0}}(\mathbb{R})=(\alpha\omega+\overline{\omega\sigma})^{\perp_{M_{i,0}}}\,\cap\,Q_{M_{i,0}}(\mathbb{R})$.\\ But since $(\alpha\omega+\overline{\alpha\omega})^{\perp_{M_{i,0}}}$ has signature $(1,j)$, $j>3$ we have that $\exists\,\,u\,\in\,Q_{M_{i,0}}(\mathbb{R})\cap (\alpha\omega+\overline{\alpha\omega})^{\perp_{M_{i,0}}}$. Since $Q_{M_{i,0}}(\mathbb{Q})$ is non-empty it is dense inside $Q_{M_{i,0}}(\mathbb{R})$, therefore $\exists\,\{v_n\}$ such that $[v_n]\rightarrow [u]$ in $\mathbb{P}(M_{i,0}\otimes\mathbb{C})$, where the $v_n$ are primitive isotropic vectors inside $M_{i,0}$. Thus we can apply \Ref{lem}{denselemma} to find a sequence $\{w_n\}$ of elements of $\mathcal{P}_{i,exc}$ such that $[v_n]\rightarrow [u]$ and $w_n\perp v_n$.
\end{proof}
\end{lem}

\begin{thm}\label{thm:standard_morph}
Let $(X,\varphi)$ be a couple consisting in a manifold of \ktipo and a symplectic automorphism of order $i=2,3$ or 5. Suppose moreover that $S_\varphi(X)\cong E_8(-2)$ if $i=2$, $S_\varphi(X)\cong K_{12}(-2)$ if $i=3$ and $S_\varphi(X)\cong S_{5.K3}$ if $i=5$. Then the couple $(X,\varphi)$ is standard.

\begin{proof}
Keeping notation as above we have $T_{\varphi}(X)\cong M_i$. Let $f$ be a marking of $X$ such that $\mathcal{P}(X,f)\subset\mathbb{P}(M_{i,0}\otimes\mathbb{C})$ and $f(S_\varphi(X))\perp M_{i,0}$.  
Moreover let $\mathcal{X}\rightarrow U$ be a maximal family of deformations of the couple $(X,\varphi)$ as in \eqref{deform} and let $F$ be a marking of $\mathcal{X}$ compatible with $f$ such that $V=\{\mathcal{P}(\mathcal{X}_t,F_t),\,t\in U\}$ is a small neighbourhood of $\mathcal{P}(X,f)$. By \Ref{lem}{denselemma_235} there exist a point $v\in V$ and a 2-divisible primitive vector $e$ of square $(-2)$ such that $v\perp e$. Since the global Torelli theorem holds we can use \Ref{thm}{graph} on the manifold $\mathcal{X}_u$ such that $\mathcal{P}(\mathcal{X}_u,F_u)=v$. This gives that $\mathcal{X}_u$ is bimeromorphic to the Hilbert square of a certain K3 surface $S$.\\ Thus we get a bimeromorphic morphism $\varphi$ on $S^{[2]}$ such that $S_{\varphi}(S^{[2]})\subset Pic(S)\subset Pic(S^{[2]})$, where 
\begin{equation}\nonumber
Pic(S)=\{t\in Pic(S^{[2]}), e\perp t\}.
\end{equation}
By \Ref{thm}{nik_k31} and \Ref{thm}{nik_k33} we have a symplectic morphism $\psi$ of order $i$ on $S$ given by the action of $\varphi$ on $e^\perp\cong H^2(S,\mathbb{Z})$ which induces a symplectic automorphism $\psi^{[2]}$ on $S^{[2]}$. Furthermore the birational map $(\psi^{[2]})^{i-1}\circ\varphi$ induces the identity on $H^2(S^{[2]},\mathbb{Z})$, therefore it is biregular (sends any K\"{a}hler class into itself), and it is also the identity. 
 This means $\varphi=\psi^{[2]}$, which implies our claim.
\end{proof}
\end{thm}
\setcounter{prop}{0}
\chapter{Fixed loci of automorphisms}\label{cap:fixed}
In this chapter we will use an approach similar to \Ref{sec}{fix_k3} to compute the fixed locus of automorphisms on \hk manifolds. In this more general setting the computations are harder and we are able only to provide partial results, namely only in the setting of manifolds of \ktipo.\\

\section{Fixed point formulas}
Our main tool is a formula first devised by Atiyah and Singer \cite{aty} in the analytic context and then specialized by Donovan \cite{don} in the algebraic case. It is usually referred to as Holomorphic Lefschetz-Riemann-Roch formula.\\ The formula works in a broader context, but we will use it only for coherent sheaves endowed with an automorphism of finite order induced by an automorphism $\psi$ of the manifold $Y$.
\begin{defn}
Let $Y$ be a Complex manifold and let $\psi\in Aut(Y)$ be an automorphism of finite order. Then we define the following:
\begin{itemize}

\item $ct(\mathcal{F},\psi)\,\in\,H^*(Y,\mathbb{Q})\otimes\mathbb{C}$\index{Chern trace, $ct(\,,\,)$} is the chern trace with respect to $\psi$ of $\mathcal{F}$, where $\mathcal{F}$ is a coherent sheaf and, with an abuse of notation,  $\psi\,:\mathcal{F}\,\rightarrow\,\mathcal{F}$ is the automorphism induced by $\psi$ on $Y$. It is given as follows: suppose $\mathcal{F}$ decomposes as the direct sum of $\mathcal{L}_s$ which are eigensheaves of eigenvalue $s$ for $\psi$, then \\ $ct(\mathcal{F},\psi)=\sum s\times ch(L_s)$ where $ch$ is the usual chern character.\\
\item $(\oplus \Lambda^t N^{\vee}_Z,\lambda_Z)$\index{Exterior normal sheaf and induced morphism, $(\oplus \Lambda^t N^{\vee}_Z,\lambda_Z)$} is a couple consisting in a sheaf, defined for every variety $Z$ fixed by $\psi$ and an automorphism given as follows: with a further abuse of notation let $\psi$ be also the natural automorphism on $N^{\vee}_Z$ induced by $\psi$ on $Y$ and define $\lambda_Z=(-1)^t\Lambda^t\psi$ on $\Lambda^t N^{\vee}_Z$. 
\end{itemize}
\end{defn}
Now the formula can be written as:
\begin{equation}\label{donovan_formula}
(-1)^t tr(H^t(\mathcal{F}))=\sum\int \frac{Todd(Z)ct(\mathcal{F}_{|Z},\psi_{|Z})}{ct(\oplus_{t} \Lambda^t N^{\vee}_Z,\lambda_Z)},
\end{equation}
where the sum is taken on all $Z$ varieties inside the fixed locus of $\psi$.\\
Another useful tool is a formula developed recently by Boissi\`{e}re, Nieper Wi\ss kirchen and Sarti. We must stress that this formula applies only to manifolds of \ktipo, although potentially it might be possible to extend it to some more \hk manifolds, see \cite{boiniesar2}. Let us introduce a few more notations: let $X$ be a \hk manifold and let $G$ be a group of automorphisms, let $a_G(X)=l(A_{S_G(X)})$\index{Lattice, Length of the discriminant group, $a_G(X)$} and let $m_G(X)=rank(S_G(X))/(|G|-1)$\index{Lattice, Normalized rank of $S_G(X)$, $m_G(X)$}. Notice that \Ref{oss}{G_tors} implies that these are both integers when $G$ has prime order.
\begin{thm}\cite[Theorem 1.2]{boiniesar2}
Let $X$ be a manifold of \ktipo and let $G$ be a group of automorphisms of prime order $p$, $3\leq p\leq 19$, $p\neq 5$. Then the following holds:
\begin{align}\label{bns_fixed}
Dim(H^*(X^G,\mathbb{Z}_{/(p)}))& = 324-2a_G(X)(25-a_G(X))-\\\nonumber
& +(p-2)m_G(X)(25-2a_G(X))+\\\nonumber
& +\frac{1}{2}m_G(X)((p-2)^2m_G(X)-p).
\end{align}
\end{thm}

\section{Automorphisms on fourfolds of \ktipo}
In this section we specialize the computations to a manifold $X$ of \ktipo and a symplectic automorphism $\varphi$ of prime order $p$.\\ We wish to remark that the results contained herein will be instrumental in the proof of \Ref{thm}{prime_autom_k32}. Moreover we will only compute what is needed in the proof of the above cited theorem since it will already give a full classification of all possible fixed loci $X^\varphi$.
Next is a remark taken from \cite{cam}:
\begin{oss}\label{oss:fixed_components}
A Symplectic automorphism of finite order on a manifold $X$ of \ktipo has a smooth fixed locus, moreover its connected components are one of the following:
\begin{itemize}
\item An isolated point.
\item An abelian surface.
\item A K3 surface.
\end{itemize}
\begin{proof}
The statement on smoothness is proven in \cite{don} in the more general case of actions by finite groups, while the latter statement is due to the classification of \kahl Symplectic surfaces and to the fact that $TX_{Z}=U\,\oplus\,V\,\oplus\,V^{\vee}$ where Z is a connected component of the fixed locus, $U=TZ$ is the $1$-eigenspace of the action and the symplectic form is locally defined on $\Lambda^2U\oplus\,V\otimes V^{\vee}$.
\end{proof}
\end{oss}
\subsection{p=2}
Let $\tau$ be the trace on $H^2(X,\mathbb{C})$ of a symplectic involution $\varphi$, the following is a result due to Camere \cite{cam}:
\begin{prop}\label{prop:camere_fixedp2}
Let $X,\varphi$ and $\tau$ be as before. Then the fixed locus of $\varphi$ and the values of $\tau$ are one of the following:
\begin{itemize}
\item $X^{\varphi}=$ 28 isolated points and 1 K3 surface, $\tau=5$, 
\item $X^{\varphi}=$ 12 isolated points and at least 1 abelian surface, $\tau=-3$,
\item $X^{\varphi}=$ 36 isolated points and at least 1 abelian surface, $\tau=3$.
\end{itemize}
\end{prop} 
We give an improvement of this result by eliminating the last two items:
\begin{thm}\label{thm:fixed_p2}
Let $X$ be a Hyperk\"{a}hler manifold of \ktipo with a symplectic involution $\varphi$. Then  the fixed locus $X^{\varphi}$ consists of 28 isolated points and one K3 surface. 
Moreover the lattice $S_{\varphi}(X)$ has rank 8.
\begin{proof}
By \Ref{prop}{camere_fixedp2} we have that $rank(T_{\varphi}(X))\geq 11$. By \Ref{prop}{fixed_abel} we therefore have that symplectic involutions cannot fix complex tori, hence we have our claim.
\end{proof}
\end{thm}

\subsection{p=3}
As before let $X$ be a manifold of \ktipo and let now $\varphi$ be a symplectic automorphism of order 3. 
We proceed to classify $X^\varphi$ using \eqref{donovan_formula}. 
We will work this formula in detail for $\mathcal{F}=\mathcal{O}_X,\,\Omega^1_X,\,\Omega^2_X$ which are the sheaves whose cohomology generates all of $H^*(X)$.\\ In this subsection let  $\omega=e^{\frac{2\pi i}{3}}$ and let $a$ be the dimension of the $\omega$-eigenspace on 
$H^{1,1}(Y)$ (notice that $a=m_G(X)$).
\begin{thm}\label{thm:fixedp3}
Let $X$ and $\varphi$ be as before, then one of the following holds:
\begin{itemize}
\item $X^\varphi$ consists of 27 isolated points and $a=6$.
\item $X^\varphi$ consists of at least 1 abelian surface, $a=9$.
\item $X^\varphi$ consists of 6 isolated points and 2 K3 surfaces, $a=5$.
\end{itemize}
\end{thm}
To my knowledge these computations where also independently done by Camere \cite{cam2}.
\begin{oss}
We can use proposition \Ref{prop}{fixed_abel} to conclude that no fixed abelian surface exists in the last and first cases but we cannot use it to exclude the middle case, infact an example of such an action exists (see \cite{kaw}, \cite{nami} and \Ref{ex}{fano3}).\\
\end{oss}

We begin the proof by evaluating $Todd(Z)/ct(\oplus \Lambda^t N^{\vee}_Z,\lambda_Z)$ for all possible connected components $Z$ of $X^\varphi$.
\begin{prop}\label{prop:ctordine3}
Let $X,\varphi$ be as above and let $Z$ be a connected component of $X^\varphi$. Then the following hold:
\begin{itemize}
\item $Todd(Z)/ct(\oplus \Lambda^t N^{\vee}_Z,\lambda_Z)=\frac{1}{3}+\frac{i\sqrt{3}}{9}c_1(N_Z^{\vee\omega})-\frac{5}{36}c_2(Z)+\frac{c_2(X)[Z]}{6}$ if Z is a surface. 
\item $Todd(Z)/ct(\oplus \Lambda^t N^{\vee}_Z,\lambda_Z)=\frac{1}{9}$ if Z is an isolated point.
\end{itemize}
Here $N_Z^{\vee\omega}$ is the $\omega$ eigensheaf inside $N^{\vee}_Z$.
\begin{proof}
First of all we have $Todd(Z)=1$ for an isolated point and $Todd(Z)=1+\frac{1}{12}c_2(Z)$ for a K3 or abelian surface.\\ Locally around an isolated fixed point the automorphism can be written as $\left( \begin{array}{cccc}
\omega & 0 & 0 & 0\\ 0 & \omega & 0 & 0\\ 0& 0& \overline{\omega} & 0\\ 0 & 0 & 0& \overline{\omega}\ \end{array} \right)$ since the matrix must be inside $Sp(2,\mathbb{C})$, of order 3 and must not have 1 as an eigenvalue (otherwise the automorphism would locally fix at least a curve), moreover this is also the local form on $N^{\vee}_Z$ hence it decomposes as a trivial $\omega$-eigensheaf of rank 2 and a trivial $\overline{\omega}$-eigensheaf of rank 2, which implies that $\oplus \Lambda^t N^{\vee}_Z$ decomposes as trivial eigensheaves of eigenvalues $1,-\omega,-\overline{\omega},\omega,\overline{\omega}$ of rank respectivelly $6,4,4,1,1$. Therefore $ct(\oplus \Lambda^t N^{\vee}_Z,\lambda_Z)=9$.\\
Locally around a fixed surface the automorphism can be written as $\left( \begin{array}{cccc}
\omega & 0 & 0 & 0\\ 0 & \overline{\omega} & 0 & 0\\ 0& 0& 1 & 0\\ 0 & 0 & 0& 1 \end{array} \right)$ which decomposes as $N^{\vee}_Z\oplus TZ$ where $TZ$ is the 1-eigensheaf. This means we need only the chern classes of $N^{\vee\omega}_Z$ and $N^{\vee\overline{\omega}}_Z$ to evaluate $ct$.\\ We obtain all the chern classes of $N^{\vee}_Z$ by the exact sequence 
\begin{equation}\nonumber
0\rightarrow TZ\rightarrow TY_{|Z}\rightarrow N_Z\rightarrow 0,
\end{equation}
which are $c_1(N^v_Z)=0$ and $c_2(N^v_Z)=c_2(X)[Z]-c_2(Z)$.\\ Now we consider the exact sequence $0\rightarrow N^{\vee\omega}_Z\rightarrow N^{\vee}_Z\rightarrow N^{\vee\overline{\omega}}_Z\rightarrow 0$.\\ Thus we have
\begin{equation} 
c_1(N_Z^{\vee})=0=c_1(N_Z^{\vee\omega})+c_1(N_Z^{\vee\overline{\omega}})
\end{equation} and 
\begin{equation} c_2(N^{\vee}_Z)=c_2(X)[Z]-c_2(Z)=c_1(N_Z^{\vee\overline{\omega}})c_1(N_Z^{\vee\omega}).
\end{equation}
This gives 
\begin{equation}
ct(\oplus \Lambda^t N^{\vee}_Z,\lambda_Z)=3-i\sqrt{3}c_1(N_Z^{\vee\omega})+\frac{c_2(Z)}{2}-\frac{c_2(X)[Z]}{2}.
\end{equation}
And inverting it we obtain our claim.
\end{proof}
\end{prop}
We need now only to evaluate $ct(\Omega^1_{X|Z})$ and $ct(\Omega^2_{X|Z})$:
\begin{prop}\label{prop:ctomega3}
Let $X,\varphi$ be as above and let $Z$ be a connected component of $X^\varphi$. Then the following hold:
\begin{itemize}
\item $ct(\Omega^1_{X|Z})=-2$ if Z is an isolated fixed point.
\item $ct(\Omega^1_{X|Z})=1+i\sqrt{3}c_1(N_Z^{\vee\omega})+\frac{c_2(X)[Z]}{2}-\frac{3}{2}c_2(Z)$ if Z is a fixed surface.
\item $ct(\Omega^2_{X|Z})=3$ if Z is an isolated fixed point.
\item $ct(\Omega^2_{X|Z})=2i\sqrt{3}c_1(N_Z^{\vee\omega})+c_2(X)[Z]$ if Z is a fixed surface.
\end{itemize}
\begin{proof}
These computations are easier: we have $\Omega^1_{X|Z}=N^{\vee}_Z$ for Z an isolated fixed point, hence $ct(\Omega^1_{X|Z})=-2$ and $\Omega^1_{X|Z}=N^{\vee}_Z\oplus TZ$ for Z a fixed surface. Therefore $ct(\Omega^1_{X|Z})=ct(N^{\vee}_Z)+ct(TZ)$ which yields our claim.\\The sheaf $\Omega^2_{X|Z}$ is locally the exterior product of the preceeding, thus we obtain the desired result.
\end{proof}
\end{prop}
Now we have all we need to obtain the theorem, let us denote $N$ the number of isolated fixed points, $K$ the number of fixed K3s and $A=\sum\int_Z c_2(X)[Z]$. Let us further remark that $\int_Z c_2(Z)=0$ for an abelian surface and $\int_Z c_2(Z)=24$ for a K3.
By applying Donovan's formula to the sheaves $\mathcal{O}_X,\,\Omega^1_X,\,\Omega^2_X$ we obtain the following system:
\begin{subequations}
\begin{align}
3=& \frac{N}{9}-\frac{10K}{3}+\frac{A}{6},\\
6a-42 = & -\frac{2N}{9}-\frac{70K}{3}+\frac{2A}{3},\\
\frac{N}{3}-16K+A = &\frac{9a^2}{2}-\frac{129a}{2}+234.
\end{align} 
\end{subequations}
We use the first equation to eliminate $A$ from the other two and we obtain:
\begin{subequations}
\begin{align}
6a-54 = & -\frac{2N}{3}-10K,\\
-\frac{N}{3}+4K = & \frac{9a^2}{2}-\frac{129a}{2}+216.
\end{align}
\end{subequations}
Since $N,K\geq0$ we have $a\leq9$, and by eliminating $N$ from the last equation we obtain
\begin{equation}
9a^2-135a+(486-18K)=0,
\end{equation}
whose integer solutions with $a\leq9$ give us the three cases described in the theorem.
For the proof of \Ref{thm}{prime_autom_k32} we need also to specialize to one particular case:

\begin{prop}\label{prop:fixed_p3_sarti}
Let $\varphi$ and $X$ be as before, let moreover \\$S_\varphi(X)=W(-1)$. Then $X^\varphi$ consists of one abelian surface.
\begin{proof}
Let $G$ be the group of automorphisms generated by $\varphi$. We have $a_G(X)=5$ and $m_G(X)=9$. Using \eqref{bns_fixed} we obtain that Dim$(H^*(X^\varphi))=16$ and, by \Ref{thm}{fixedp3}, this implies that $X^\varphi$ consists of one abelian surface.
\end{proof}
\end{prop}

\subsection{p=5}

Let $X$ be a manifold of \ktipo and let $\varphi$ be a symplectic automorphism of order 5.
As before we use \eqref{donovan_formula} on $\mathcal{F}=\mathcal{O}_X,\,\Omega^1_X,\,\Omega^2_X$.\\
In this subsection let $\omega=e^{\frac{2\pi i}{5}}$ and let $a$ be the dimension of the $\omega$-eigenspace of $\varphi$ on $H^{1,1}(X)$.
\begin{thm}\label{thm:fixedp5}
Let $X$, $\varphi$ be as before. Then $X^\varphi$ consists of 14 isolated points and $a=4$.
\end{thm}
\begin{oss}
In this case there is an important change: the local action of the automorphism on the connected components of the fixed loci is not univocally determined by its topological type but there are several choices:
\begin{enumerate}
\item Fixed points with local action of $\varphi$ given by $\left( \begin{array}{cccc}
\omega & 0 & 0 & 0\\ 0 & \omega & 0 & 0\\ 0& 0& \overline{\omega} & 0\\ 0 & 0 & 0& \overline{\omega}\ \end{array} \right)$ and we call those \emph{points of the first kind}.
\item Fixed points with local action given by $\left( \begin{array}{cccc}
\omega^2 & 0 & 0 & 0\\ 0 & \omega^2 & 0 & 0\\ 0& 0& \overline{\omega}^2 & 0\\ 0 & 0 & 0& \overline{\omega}^2\ \end{array} \right)$ and we call those \emph{points of the second kind}.
\item Fixed points with local action given by $\left( \begin{array}{cccc}
\omega & 0 & 0 & 0\\ 0 & \omega^2 & 0 & 0\\ 0& 0& \overline{\omega} & 0\\ 0 & 0 & 0& \overline{\omega}^2\ \end{array} \right)$ and we call those \emph{points of the third kind}.
\item Fixed surfaces with local action given by $\left( \begin{array}{cccc}
\omega & 0 & 0 & 0\\ 0 & \overline{\omega} & 0 & 0\\ 0& 0& 1 & 0\\ 0 & 0 & 0& 1 \end{array} \right)$ and we call those \emph{surfaces of the first kind}. 
\item Fixed surfaces with local action given by $\left( \begin{array}{cccc}
\omega^2 & 0 & 0 & 0\\ 0 & \overline{\omega}^2 & 0 & 0\\ 0& 0& 1 & 0\\ 0 & 0 & 0& 1 \end{array} \right)$ and we call those \emph{surfaces of the second kind}.
\end{enumerate}
\end{oss}
Moreover there is another important
\begin{oss}\label{oss:symmetry}
Given a Symplectic automorphism $\varphi$ of order 5 with $N_1$ fixed points of the first kind, $N_2$ fixed points of the second kind, $N_3$ fixed points of the third kind, $S_1$ fixed surfaces of the first kind and $S_2$ fixed surfaces of the second kind we have that $\varphi^2$ is a symplectic automorphism of order 5 with $N_2$ fixed points of the first kind, $N_1$ fixed points of the second kind, $N_3$ fixed points of the third kind, $S_2$ fixed surfaces of the first kind and $S_1$ fixed surfaces of the second kind.\\ Furthermore among the surfaces we have the same number of K3's and abelian surfaces being of the first kind in one case and of the second in the other.
\end{oss}
We now need to evaluate the same characteristic classes as before and we start with the following:
\begin{prop}
Let $X,\varphi$ be as above and let $Z$ be a connected component of $X^\varphi$. Then the following hold:
\begin{itemize}
\item \small{$ct(\oplus \Lambda^t N^{\vee}_Z,\lambda_Z)=\frac{15-5\sqrt{5}}{2}$} for $Z$ a point of the first kind.
\item \small{$ct(\oplus \Lambda^t N^{\vee}_Z,\lambda_Z)=\frac{15+5\sqrt{5}}{2}$} for $Z$ a point of the second kind.
\item \small{$ct(\oplus \Lambda^t N^{\vee}_Z,\lambda_Z)=5$} for $Z$ a point of the third kind.
\item \small{$ct(\oplus \Lambda^t N^{\vee}_Z,\lambda_Z)=\frac{5-\sqrt{5}}{2}-c_1(N_Z^{\vee\omega})\frac{\sqrt{10+2\sqrt{5}}}{2}+(c_2(X)[Z]-c_2(Z))(\frac{-1-\sqrt{5}}{4})$} for $Z$ a surface of the first kind.
\item \small{$ct(\oplus \Lambda^t N^{\vee}_Z,\lambda_Z)=\frac{5+\sqrt{5}}{2}-c_1(N_Z^{\vee\omega^2})\frac{\sqrt{10-2\sqrt{5}}}{2}+(c_2(X)[Z]-c_2(Z))(\frac{-1+\sqrt{5}}{4})$} for $Z$ a surface of the second kind.
\end{itemize}
\begin{proof}
The evaluation goes as in \Ref{prop}{ctordine3}, we need only to change the eigenvalues and we obtain our claim. 
\end{proof}
\end{prop}
\begin{prop}
Let $X,\varphi$ be as above and let $Z$ be a connected component of $X^\varphi$. Then the following hold:
\begin{itemize}
\item \small{$Todd(Z)/ct(\oplus \Lambda^t N^{\vee}_Z,\lambda_Z)=\frac{5+\sqrt{5}}{10}+c_1(N_Z^{\vee\omega})\frac{\sqrt{10+2\sqrt{5}}}{15-5\sqrt{5}}-\frac{9+3\sqrt{5}}{20}c_2(X)[Z]+\frac{59+19\sqrt{5}}{120}c_2(Z)$} for $Z$ a surface of the first kind.
\item \small{$Todd(Z)/ct(\oplus \Lambda^t N^{\vee}_Z,\lambda_Z)=\frac{5-\sqrt{5}}{10}+c_1(N_Z^{\vee\omega^2})\frac{\sqrt{10-2\sqrt{5}}}{15+5\sqrt{5}}+\frac{-5+2\sqrt{5}}{10}c_2(X)[Z]+\frac{13-5\sqrt{5}}{24}c_2(Z)$} for $Z$ a surface of the second kind.
\end{itemize}
\end{prop}
The final computation is the evaluation of $ct(\Omega^1_{Z|Z})$ and $ct(\Omega^2_{Z|Z})$:
\begin{prop}
Let $X,\varphi$ be as above and let $Z$ be a connected component of $X^\varphi$. Then the following hold:
\begin{itemize}
\item $ct(\Omega^1_{Z|Z})=-1+\sqrt{5}$ if Z is a point of the first kind.
\item $ct(\Omega^1_{Z|Z})=-1-\sqrt{5}$ if Z is a point of the second kind.
\item $ct(\Omega^1_{Z|Z})=-1$ if Z is a point of the third kind.
\item \small{$ct(\Omega^1_{Z|Z})=\frac{3+\sqrt{5}}{2}+c_1(N_Z^{\vee\omega})\frac{\sqrt{10+2\sqrt{5}}}{2}-\frac{1+\sqrt{5}}{4}c_2(X)[Z]-\frac{3+\sqrt{5}}{4}c_2(Z)$} if $Z$ is a surface of the first kind.
\item \small{$ct(\Omega^1_{X|Z})=\frac{3-\sqrt{5}}{2}+c_1(N_Z^{\vee\omega^2})\frac{\sqrt{10-2\sqrt{5}}}{2}+\frac{-1+\sqrt{5}}{4}c_2(X)[Z]+\frac{-3+\sqrt{5}}{4}c_2(Z)$} if $Z$ is a surface of the second kind.
\item $ct(\Omega^2_{X|Z})=\frac{7-\sqrt{5}}{2}$ if Z is a point of the first kind.
\item $ct(\Omega^2_{X|Z})=\frac{7+\sqrt{5}}{2}$ if Z is a point of the second kind.
\item $ct(\Omega^2_{X|Z})=1$ if Z is a point of the third kind.
\item $ct(\Omega^2_{X|Z})=1+\sqrt{5}+c_1(N_Z^{\vee\omega^2})\sqrt{10+2\sqrt{5}}-\frac{1+\sqrt{5}}{2}c_2(X)[Z]+c_2(Z)$ if $Z$ is a surface of the first kind.
\item $ct(\Omega^2_{X|Z})=1-\sqrt{5}+c_1(N_Z^{\vee\omega^2})\sqrt{10-2\sqrt{5}}+\frac{-1+\sqrt{5}}{2}c_2(X)[Z]+c_2(Z)$ if $Z$ is a surface of the second kind.
\end{itemize}
\begin{proof}
This computation mimics that of \Ref{prop}{ctomega3} only with different eigenvalues.
\end{proof}
\end{prop}
Now let us call $N_1$ the number of points of the first kind, $N_2$ the number of points of the second kind, $N_3$ the number of points of the third kind, $K_1$ the number of K3's of the first kind, $K_2$ the number of K3's of the second kind, $A_1=\int{c_2(X)[Z]}$ over surfaces of the first kind and $A_2=\int{c_2(X)[Z]}$ over surfaces of the second kind.\\
Summing all these propositions and dividing the rational part from the irrational one we have: 
\begin{subequations}
\small{
\begin{align}
60 = & 6N_1+6N_2+4N_3-9A_1-10A_2+216K_1+260K_2,\\ 
0 = & 2N_1-2N_2-3A_1+4A_2+76K_1-100K_2, \label{antiAp5}\\
100a-420 = & -8N_1-8N_2-2N_3-23A_1-25A_2+452K_1+440K_2,\\ 
0 = & 4N_1-4N_2-10A_1+11A_2+194K_1-168K_2 \label{A1p5},\\
4680-2150a+250a^2 = & 16N_1+16N_2+4N_3-74A_1-70A_2+1816K_1+1960K_2,\\ 
0 = & 4N_1-4N_2-34A_1+36A_2+744K_1-792K_2.\label{A2p5}
\end{align}}
\end{subequations}
We use \eqref{antiAp5} in \eqref{A1p5} and \eqref{A2p5} to eliminate $A_1$ and $A_2$ and we obtain
\begin{align}
A_1  = & \frac{225}{7}K_1-\frac{332}{7}K_2,\nonumber\\
A_2  = & \frac{62}{7}K_1-\frac{144}{7}K_2.\nonumber
\end{align}
Using this and multiplying by 7 we get the following system
\begin{subequations}
\begin{align}
420 = & 42N_1+42N_2+28N_3-1133K_1+6248K_2, \label{sym1p5}\\
0 = & 14N_1-14N_2-105K_1+280K_2, \label{N1p5}\\
700a-2940 = & -56N_1-56N_2-14N_3-3561K_1+14316K_2, \\
32760-15050a+1750a^2 = & 112N_1+112N_2+28N_3-8278K_1+48368K_2, \\
7A_1  = & 225K_1-332K_2,\\
7A_2  = & 62K_1-144K_2.
\end{align}
\end{subequations}
We can now use \Ref{oss}{symmetry} on \eqref{sym1p5} to obtain
\begin{align}
420 = & 42N_1+42N_2+28N_3-1133K_1+6248K_2, \nonumber\\
420 = & 42N_2+42N_1+28N_3-1133K_2+6248K_1. \nonumber
\end{align}
\ie $42(N_2-N_1)+3113(K_2-K_1)=0$ which gives us another equation. We now use \eqref{N1p5} to eliminate $N_1$:
\begin{subequations}
\begin{align}
210 = & 42N_2+14N_3-409K_1+2704K_2, \label{N3p5}\\
3953K_2 = & 3428K_1,\\
700a-2940 = & -112N_2-14N_3-3981K_1+15436K_2, \label{sym2p5}\\
32760-15050a+1750a^2 = & 224N_2+28N_3-8338K_1+46128K_2, \\
7A_1 = & 225K_1-332K_2,\\
7A_2 = & 62K_1-144K_2,\\
14N_1 = & 14N_2+105K_1-280K_2. 
\end{align}
\end{subequations}
We finally use \eqref{N3p5} to eliminate $N_3$ and \Ref{oss}{symmetry} on \eqref{sym2p5} to obtain 
\begin{align*}
3953K_2 = & 3428K_1,\\
853K_2 = & 1728K_1,\\
700a-2730 = & -70N_2-4390K_1+18140K_2, \\
33180-15050a+1750a^2 = & 140N_2-7520K_1+40720K_2, \\
7A_1 = & 225K_1-332K_2,\\
7A_2 = & 62K_1-144K_2,\\
14N_1 = & 14N_2+105K_1-280K_2, \\
14N_3 = & 210-42N_2+409K_1-2704K_2. 
\end{align*}
We can easily see this implies $K_1=K_2=A_1=A_2=0$ \ie the result of \Ref{thm}{fixedp5}.
\subsection{p=7}
We will not use Donovan's formula in this case, the only tool we will need for \Ref{thm}{prime_autom_k32} is an easy application of \eqref{bns_fixed}.
\begin{prop}\label{prop:fixed_p7_sarti}
Let $X$ be a manifold of \ktipo and let $\varphi$ be a symplectic automorphism of order 7 such that $S_{\varphi}(X)=S_{7.K3}$ as defined in \Ref{ex}{S_7K3}. Then $X^\varphi$ consists of 9 isolated points.
\begin{proof}
Let $G$ be the group of automorphisms generated by $\varphi$. We have $a_G(X)=m_G(X)=3$, therefore by \eqref{bns_fixed} we obtain Dim$(H^*(X^G))$=9, which implies by \Ref{oss}{fixed_components} our claim.
\end{proof}
\end{prop}

\subsection{p=11}
Again we avoid using Donovan's formula and we only make a simple computation with \eqref{bns_fixed}.
\begin{prop}\label{prop:fixed_p11_sarti}
Let $\varphi$ be a symplectic automorphism of order 11 of a \hk manifold $X$ of \ktipo such that $a_\varphi=m_\varphi=2$. Then $X^{\varphi}$ consists of 5 isolated points.
\begin{proof}
Using \eqref{bns_fixed} we see that $dim(H^*(X^{\varphi}))=5$. Since it consists of symplectic varieties we obtain our claim.
\end{proof} 
\end{prop}
\setcounter{prop}{0}
\chapter{Sporadic groups and Symplectic Automorphisms}\label{cap:groups_sporadic}
This chapter is devoted to obtain an analogue of \Ref{cap}{k3_case} for what concerns the link between symplectic automorphisms on a \hk manifold and isometries of its second cohomology.

\section{Automorphisms and cohomology on \hk manifolds}\label{sec:gen_lattices}
In this section we prove some useful general properties concerning automorphisms of \hk manifolds and then we specialize to manifolds which are not of \kntipo. We are thus able to provide limitations on the order of finite symplectic automorphisms on those manifolds. We also provide a way to compute the coinvariant lattice for those automorphisms using isometries of certain well known unimodular lattices. These results are not effective, \ie there exist isometries of the above cited lattices which do not come from automorphisms of \hk manifolds.
We wish to remark that some among these results are already contained in \cite{beau2}, such as most of \Ref{lem}{gaction_gen} and \eqref{exactgroup}.\\
Throughout this section $\overline{G}$ \index{Group, Automorphism group, $\overline{G}$} will denote a finite group of automorphisms on a \hk manifold $X$.
\begin{defn}\label{defn:inv_locus}
Let $\overline{G}$ be a group acting faithfully on a \hk manifold $X$, we define $T_{\overline{G}}(X)$ inside $H^2(X,\mathbb{Z})$ to be the subgroup fixed by the induced action of $\overline{G}$ on $H^2(X,\mathbb{Z})$. Moreover we define the co-invariant locus $S_{\overline{G}}(X)\subset H^2(X,\mathbb{Z})$ as $T_{\overline{G}}(X)^{\perp}$. 
The fixed locus of $\overline{G}$ on $X$ will be denoted $X^{\overline{G}}$ as before.
\end{defn}

We wish to remark that the map
\begin{equation}\label{numap}
Aut(X)\stackrel{\nu}{\rightarrow} O(H^2(X,\mathbb{Z}))\index{Projection from automorphisms to Hodge isometries, $Aut(X)\stackrel{\nu}{\rightarrow} O(H^2(X,\mathbb{Z}))$}
\end{equation}
might have nontrivial kernel if $X$ is not of $K3^{[n]}$-type, as in \Ref{ex}{kum_invol} and \Ref{ex}{kum_trivial}. We will soon discuss in greater detail the injectivity of $\nu$. We will call $G$ the image of $\overline{G}$ by $\nu$. Obviously $S_G(H^2(X,\mathbb{Z}))=S_{\overline{G}}(X)$, therefore we will not distinguish between the two notations.
 Moreover we have the following exact sequence for any finite group $G$ of Hodge isometries on $H^2(X,\mathbb{Z})$:
\begin{equation}\label{exactgroup}
1\,\rightarrow\, G_0\,\rightarrow\,G\,\stackrel{\pi}{\rightarrow}\,\Gamma_m\,\rightarrow\,1,
\end{equation}
where $\Gamma_m\subset U(1)$ is a cyclic group of order $m$. In fact the action of $G$ on $H^{2,0}$ is the action of a finite group on $\mathbb{C}$.
Let $\gamma_X$ be the following useful map:
\begin{equation}\label{sigmamap}
\gamma_X\,:\,T(X)\,\rightarrow\,\mathbb{C}.\index{Evaluation map on Transcendental part, $\gamma_X$}
\end{equation}
Here $\gamma_X(x)=(\sigma,x)_X$, which has kernel $T(X)\cap S(X)=0$.\\
\begin{oss}
We wish to remark that recently Oguiso \cite{ogu4} proved that in the case of manifolds of Kummer $n$-type the map
\begin{equation}
Aut(X)\rightarrow Aut(H^*(X,\mathbb{Z}))
\end{equation}
has trivial kernel.
\end{oss}

\begin{lem}\label{lem:h2_trivial_k3n}
Let $X$ be a manifold of \kntipo. Then the map $Aut(X)\stackrel{\nu(X)}{\rightarrow} O(H^2(X,\mathbb{Z}))$ is injective.
\begin{proof}
Hassett and Tschinkel \cite[Theorem 2.1]{ht3} proved that $Ker(\nu(X))$ is invariant under smooth deformations. Beauville \cite[Lemma 3]{beau2} proved that, if $S$ is a $K3$ surface with no nontrivial automorphisms then $Aut(S^{[n]})=Id$, therefore $Id=Ker(\nu(S^{[n]}))=Ker(\nu(X))$.
\end{proof}

\end{lem}

\begin{lem}\label{lem:gaction_gen}
Let $X$ be a \hk manifold and let $\overline{G}\subset Aut(X)$ be a group such that $\nu(\overline{G})=:G$ is finite. Then the following hold:
\begin{enumerate}
\item $g\in \overline{G}$ acts trivially on $T(X)\,\iff\,g\in \overline{G}_0$.
\item The representation of $\Gamma_m$ on $T(X)\otimes\mathbb{Q}$ splits as the direct sum of irreducible representations of the cyclic group $\Gamma_m$ having maximal rank (\ie of rank $\phi(m)$).
\end{enumerate}
\begin{proof}
First of all let us remark that without loss of generality we can consider only elements of $G$ instead of $\overline{G}$.
\begin{enumerate}
\item Let $g\in G_0$. Let us show that $g^*$ acts trivially on $T(X)\otimes\mathbb{Q}$. We start by considering the kernel of the map $g^*-Id_{T(X)}$ which is a lattice (and a Hodge substructure) $R$ inside $T(X)$. Hence, by minimality of $T(X)$, $R\otimes\mathbb{Q}$ is either 0 or $R\otimes\mathbb{Q}=T(X)\otimes\mathbb{Q}$. Considering the map \eqref{sigmamap}, since $g^*$ is a Hodge isometry we have
\begin{equation}\nonumber
\gamma_X(x)=(g^*\sigma,g^*x)=(\sigma,g^*x).
\end{equation}
Since $g^*\sigma=\sigma$ we have that $g^*x-x\,\in\,ker(\gamma_X)=T(X)\cap S(X)=0$. Thus $R$ is all of $T(X)$.\\ 
To obtain the converse we prove that $g^*\sigma=\lambda\sigma$ with $\lambda\neq1$ implies that 1 is not an eigenvalue of $g^*$ on $T(X)$. In fact
\begin{equation}
\gamma_X(x)=(g^*\sigma,g^*x)=\lambda\gamma_X(g^*x), \nonumber
\end{equation}
\ie $g^*x\neq x$.
\item The preceeding arguments show that every nontrivial element of $G/G_0$ has no eigenvalue 1 on $T(X)$ and hence also on $T(X)\otimes\mathbb{Q}$, this implies our claim.
\end{enumerate}

\end{proof} 
\end{lem}

As a consequence we have the following:
\begin{cor}[Oguiso, Schr\"{o}er, \cite{ogu2}]
Let $X$ be a \hk manifold and let $\varphi\in Aut(X)$ be an automorphism of finite order $m$ such that $\varphi(\sigma_X)=\omega\sigma_X$, where $\omega$ is a primitive $m$-th root of unity. Then the following hold:
\begin{itemize}
\item $m\leq 66$ and $\phi(m)\leq 22$ if $X$ is of $K3^{[n]}$-type.
\item $m\leq 18$ and $\phi(m)\leq 6$ if $X$ is of Kummer $n$-type.
\item $m\leq 18$ and $\phi(m)\leq 7$ if $X$ is deformation equivalent to $Og_6$.
\item $m\leq 66$ and $\phi(m)\leq 23$ if $X$ is deformation equivalent to $Og_{10}$.
\end{itemize}
\end{cor}
\begin{oss}
Let us stress that if $X$ is of $K3^{[n]}$-type and $\varphi$ is as above then only $m=23$ and $m=46$ cannot be obtained via standard automorphisms.
\end{oss}
To obtain stronger results we will need one more definition:
\begin{defn}
Let $X$ be a \hk manifold and let $\overline{G}\subset Aut_s(X)$. We say that $\overline{G}$ is \emph{quadratically nontrivial} if $\overline{G}=\nu(\overline{G})=G$. Furthermore we say that $\overline{G}$ is \emph{discriminant preserving} if its induced action on the discriminant group of $H^{2}(X,\mathbb{Z})$ is trivial.
\end{defn}
Before proceeding further let us briefly analyze what this two conditions imply in the known cases:\\
If $X$ is a manifold of $K3^{[n]}$-type then there are no quadratically trivial automorphisms by \Ref{lem}{h2_trivial_k3n}, moreover the discriminant group of $L_n\cong H^2(X,\mathbb{Z})$ has only a few isometries given by multiplying $1\in\mathbb{Z}_{/(2n-2)}$ by a square root of $1$ in $\mathbb{Z}_{/(2n-2)}$. Therefore all groups of odd order are discriminant preserving.\\
If $X$ is a manifold of Kummer $n$-type then quadratically trivial automorphisms form a group isomorphic to the semidirect product of $(\mathbb{Z}_{/(n+1)})^4$ and $\mathbb{Z}_{/(2)}$, see \cite{boiniesar}. Again isometries of the discriminant groups of $L_{K.n}\cong H^2(X,\mathbb{Z})$ have order 1 or 2 and are multiplication by a square root of 1 in  $\mathbb{Z}_{/(2n+2)}$. Therefore all groups of odd order are discriminant preserving.\\
If $X$ is deformation equivalent to $Og_6$ or $Og_{10}$ then it is not known which automorphisms can act trivially on the second cohomology. However in the six dimensional case there are only two isometries of the discriminant group of $L_{O.6}\cong H^2(Og_6,\mathbb{Z})$, namely the identity and the one induced by exchanging the two copies of $(-2)$ inside $L_{O.6}$. This implies that once more all odd order automorphism groups are discriminant preserving. Finally for the 10 dimensional case there are no nontrivial isometries of the discriminant group of $L_{O.10}\cong H^2(Og_{10},\mathbb{Z})$, therefore all groups are discriminant preserving.

Let now $\overline{G}$ be a group of automorphisms such that $G=G_0$ and $G$ is finite.
\begin{lem}\label{lem:algaction_general}
Let $X$ be a \hk manifold and let $\overline{G}\subset Aut_s(X)$ be a group such that $G$ is finite. Then the following assertions are true:
\begin{enumerate}
\item $S_G(X)=S_{\overline{G}}(X)$ is nondegenerate and negative definite.
\item $T(X)\subset T_G(X)$ and $S_G(X)\subset S(X)$.
\item Suppose $\overline{G}$ is discriminant preserving. Then $G$ acts trivially on $A_{S_G(X)}$.
\end{enumerate}
\begin{proof}
The second assertion is an immediate consequence of \Ref{lem}{gaction_gen} because $G$ acts as the identity on $\sigma$ and therefore on all of $T(X)$.\\
To prove that $S_G(X)$ and $T_G(X)$ are nondegenerate let $H^2(X,\mathbb{C})=\oplus_{\rho}U_\rho$ be the decomposition in orthogonal representations of $G$, where $U_{\rho}$ contains all irreducible representations of $G$ of character $\rho$ inside $H^2(X,\mathbb{C})$. Obviously $T_G(X)=U_{Id_{|\mathbb{Z}}}$ and $S_G(X)=H^2(X,\mathbb{Z})\cap \oplus_{\rho\neq Id}U_\rho$, which implies they are orthogonal and of trivial intersection. Hence they are both nondegenerate.\\
 Since $G$ is finite there exists a $G$-invariant K\"{a}hler class $\omega_G$ given by $\sum_{g\in G}g\omega$, where $\omega$ is any K\"{a}hler class on $X$.  
Therefore we have: 
\begin{equation}\nonumber
\sigma\mathbb{C}\oplus\overline{\sigma}\mathbb{C}\oplus\omega_G\mathbb{C}\,\subset\,T_G(X)\otimes\mathbb{C}.
\end{equation}
Hence the lattice $S_G(X)$ is negative definite.\\
To prove the last assertion let us proceed as in \Ref{lem}{biggen_g}, \ie let us choose a primitive embedding of $H^2(X,\mathbb{Z})$ into an unimodular lattice $M$ of signature $(4,r)$, where $r\geq b_2(X)-3$. And let us extend the action of $G$ trivially outside the image of $H^2(X,\mathbb{Z})$. Therefore $S_G(X)\cong S_G(M)$ and $A_{S_G(M)}\cong A_{T_G(M)}$, where the isomorphism is $G$-equivariant. $G$ acts trivially on $T_G(M)$, thus its induced action on $A_{T_G(M)}$ is trivial. Using the $G$-equivariant isomorphism we have that $G$ acts trivially also on $A_{S_G(M)}=A_{S_G(X)}$.\\
Let us specify that $M=U^4\oplus E_8(-1)^2$ if $X$ is of $K3^{[n]}$-type, $M=U^4$ if $X$ is of Kummer $n$-type, $M=U^4\oplus E_8(-1)$ if $X$ is deformatione equivalent to $Og_6$ and $M=U^4\oplus E_8(-1)^3$ if $X$ is deformation equivalent to $Og_{10}$.

\end{proof}
\end{lem}

In the rest of this section we will not consider anymore manifolds of $K3^{[n]}$-type, they will be analyzed in greater detail in \Ref{sec}{k32_case} and \Ref{sec}{k3n_case}.
Now we wish to provide some restriction on possible finite groups $G=\nu(\overline{G})$, $\overline{G}\subset Aut_s(X)$:
\begin{prop}
Let $X$ be a \hk manifold and let $\overline{G}\subset Aut_s(X)$ be a group such that $G=\nu(\overline{G})$ is finite and discriminant preserving. Then there exists an embedding $S_G(X)\rightarrow P$ and $G$ extends to a group of isometries of $P$ acting trivially on $S_G(X)^{\perp_P}$. Here $P$ is as follows:
\begin{itemize}
\item $E_8(-1)$ if $X$ is of Kummer $n$-type.
\item $E_8(-1)^2$ or $D_{16}^+(-1)$ if $X$ is deformation equivalent to $Og_6$.
\item $U\oplus\Lambda$ if $X$ is deformation equivalent to $Og_{10}$.
\end{itemize}
\begin{proof}
Since the group $G$ is discriminant preserving we have that $S_G(X)$ is negative definite and $G$ acts trivially on $A_{S_G(X)}$. We will embed $H^2(X,\mathbb{Z})$ in an unimodular lattice $W$, let us look separately at the 3 cases:
\begin{itemize}
\item[$K_n(T)$] Let us give an embedding $L_{K,n}\rightarrow U^4$ and let $T=S_G(X)^{\perp_{U^4}}$. Since $T$ is the unimodular complement of $S_G(X)$ it has the same discriminant group and $G$ acts trivially on it, therefore we can estend $G$ to a group of isometries of $U^4$ acting trivially on $T$. Let $r\leq4$ be the rank of $S_G(X)$ and let $m=l(A_{S_G(X)})\leq4$. By \Ref{lem}{nik_esiste} there exists a negative definite lattice  $T'$ of rank $8-r$ and discriminant group $A_{S_G(X)}$ with the opposite discriminant form. Therefore by \Ref{lem}{nik_immerge1} there exists a primitive embedding $S_G(X)\rightarrow P$, where $P$ is an even negative definite unimodular lattice of rank 8 (\ie $E_8(-1)$) and $G$ extends to a group of isometries of $E_8(-1)$ acting trivially on the orthogonal complement $T'$. 
\item[$Og_6$] The proof is similar, this time we embed $L_{O.6}$ in $U^5$ and we have that $G$ extends to $O(U^5)$ satisfying $S_G(U^5)=S_G(X)$ and $T_G(U^5)=T$. Let $r\leq5$ be the rank of $S_G(X)$ and let $m=l(A_{S_G(X)})\leq5$. Again by \Ref{lem}{nik_esiste} we have a negative definite lattice $T'$ of discriminant group $A_{S_G(X)}$ and opposite discriminant form. However this time if $m=5$ it has rank $16-r$. Therefore by \Ref{lem}{nik_immerge1} there exists a primitive embedding $S_G(X)\rightarrow P$, where $P$ is an even negative definite unimodular lattice of rank 16 (\ie $E_8(-1)^2$ or $D_{16}^+(-1)$) and $G$ extends to a group of isometries of $P$ acting trivially on the orthogonal complement $T'$. 
\item[$Og_{10}$] The proof goes the same, this time we embed $L_{O,10}$ into $R=U^5\oplus E_8(-1)^2$ and we remark that $T_G(R)$ contains the lattice $A_2$, therefore by \Ref{oss}{overl_group} we have $rank(S_G(X))+l(A_{S_G(X)})\leq 25$, which implies by \Ref{lem}{nik_immerge1} that $S_G(X)$ embeds into an unimodular lattice of rank $26$ and signature $(1,25)$. Since all of these lattices are isometric we obtain our claim.
\end{itemize}
\end{proof}
\end{prop}
Let us remark that we choose to embed $S_G(Og_{10})$ into an indefinite lattice for two reasons: the first is that otherwise we would have had to choose a definite lattice of rank $32$, which number in the millions. The second is that the isometry group of $U\oplus\Lambda$ has been studied in greater detail (see \cite[Chapter 27]{con}).\\
Using the group structure of definite lattices and the cyclotomic structure of some coinvariant lattices we are able to prove the following:
\begin{cor}\label{cor:kum_maxp}
Let $X$ be a \hk manifold of Kummer $n$-type and let $\varphi\in Aut_s(X)$ be an automorphism such that $\nu(\varphi)$ has prime order $p$. Then $p\leq 5$. Moreover if $p=5$ we have $S_{\varphi}(X)\cong A_4(-1)$.
\begin{proof}
The first trivial remark is that an odd order automorphism is discriminant preserving, hence $S_\varphi(X)$ is negative definite and $\varphi$ acts trivially on $A_{S_\varphi(X)}$. By \Ref{oss}{G_tors} $S_\varphi(X)$ can be given the structure of a $\mathbb{D}_p$-lattice as in \Ref{oss}{int_to_cyclo}, which implies $rank(S_\varphi(X))=(p-1)m\leq 4$, therefore $p\leq5$. If $p=5$ then $S_\varphi(X)$ is a rank $1$ $\mathbb{D}_5$-lattice, thus by \Ref{ex}{cyclo_orderp} it is a multiple of $A_4$. By looking at conjugacy classes of isometries of $E_8$ one easily sees that $S_\varphi(X)=A_4(-1)$.
\end{proof}
\end{cor}
\begin{cor}
Let $X$ be a \hk manifold deformation equivalent to $Og_6$ and let $\varphi\in Aut_s(X)$ be an automorphism such that $\nu(\varphi)$ has prime order $p$. Then $p\leq 5$. If $p=5$ we have $S_{\varphi}(X)\cong A_4(-1)$.
\begin{proof}
As in \Ref{cor}{kum_maxp} we have $p\leq 5$ due to the limitation $rank(S_\varphi(X))\leq 5$ and we can proceed in the same way for the case $p=5$.
\end{proof}
\end{cor}
\begin{cor}\label{cor:og10_maxp}
Let $X$ be a \hk manifold of deformation equivalent to $Og_{10}$ and let $\varphi\in Aut_s(X)$ be an automorphism such that $\nu(\varphi)$ has prime order $p$. Then $p\leq 19$. Moreover if $p\geq 13$ then $S_\varphi(X)\cong A_{p-1}(-1)$.
\begin{proof}
As in \Ref{cor}{kum_maxp} we have $p\leq 19$ due to the limitation $rank(S_\varphi(X))\leq 21$ and in cases $p=13,17$ or $p=19$ we would get a rank 1 negative definite $\mathbb{D}_p$-lattice $A_{p-1}(n)$. However if $n\neq -1$ these lattices do not embed primitively into $H^2(X,\mathbb{Z})$, therefore $n=-1$.
\end{proof}
\end{cor}

\section{The \ktipo case}\label{sec:k32_case}
In this section we provide a specialization of the results of \Ref{sec}{gen_lattices} to manifolds of \ktipo. This section provides the closest possible generalization of \Ref{cap}{k3_case}. We are in fact able to give classification results for prime order symplectic automorphisms and also to give a way to compute Coinvariant lattices of symplectic automorphisms. We wish to remark that in this case the map $\nu$ of \eqref{numap} is injective, therefore we will not distinguish between $\overline{G}$ and $G$. Moreover the discriminant group of $H^2(X,\mathbb{Z})=L$ is $\mathbb{Z}_{/(2)}$ if $X$ is of \ktipo, therefore it has no nontrivial isometries, \ie all possible groups of isometries $G$ are discriminant preserving. Let us suppose that $G$ is a finite group of automorphisms such that $G=G_0$. 
\begin{lem}\label{lem:algaction}
Let $X$ be a manifold of \ktipo and let $G\subset Aut_s(X)$ be a finite group. Then the following assertions are true:
\begin{enumerate}
\item $S_G(X)$ is nondegenerate and negative definite.
\item $S_G(X)$ contains no element with square -2.
\item $T(X)\subset T_G(X)$ and $S_G(X)\subset S(X)$.
\item $G$ acts trivially on $A_{S_G(X)}$.
\end{enumerate}
\begin{proof}
By \Ref{lem}{algaction_general} we need only to prove the second assertion. Assume on the contrary that we have an element $c\in S_G(X)$ such that $(c,c)=-2$. Then by \Ref{thm}{num_eff}  
 it is known that either $\pm c$ or $\pm 2c$ is represented by an effective divisor D on X. Let $D'=\sum_{g\in G}gD$ which is also an effective divisor on X, but $[D']\in\,S_G(X)\cap T_G(X)=\{0\}$. This implies $D'$ is linearly equivalent to 0, which is impossible.
\end{proof}
\end{lem}
Notice that this amounts to saying that $(S_G(X),G)$ is a negative definite Leech couple in the sense of \Ref{defn}{leech_group}.

Now we can use \Ref{thm}{graph} to give sufficient conditions for an isometry  $\psi$ of $L$ to be induced by a birational map $\psi'$ of some marked Hyperk\"{a}hler manifold $(X,f)$ such that $f\circ\psi'^*\circ f^{-1}=\psi$. Thus we obtain a generalization of \Ref{thm}{nik_k31}: 

\begin{thm}\label{thm:cohom_to_aut}
Let $G$ be a finite subgroup of $O(L)$. Suppose that the following hold:
\begin{enumerate}
\item $S_G(L)$ is nondegenerate and negative definite.
\item $S_G(L)$ contains no element with square $(-2)$.
\end{enumerate}
Then $G$ is induced by a subgroup of $Bir(X)$ for some manifold $(X,f)$ of $K3^{[2]}$-type.
\begin{proof}
By the surjectivity of the period map and by \Ref{lem}{algaction} we can consider a marked $K3^{[2]}$-type 4-fold $(X,f)$ such that $T(X)\stackrel{f}{\rightarrow} T_G(L)$ is an isomorphism and also $S(X)\stackrel{f}{\rightarrow} S_G(L)$ is.\\
Let $g\in G$, let us consider the marked varieties $(X,f)$ and $(X,g\circ f)$. They have the same period in $\Omega$ and hence by \Ref{thm}{graph} we have $f^{-1}\circ g\circ f=\Gamma_*$. Here $\Gamma=Z+\sum_{j}Y_j$ in $X\times X$, where $Z$ is the graph of a bimeromorphic map from $X$ to itself and $Y_j$'s are cycles with $codim(\pi_i(Y_j))\geq 1$.\\ 
 We will prove that all $Y_j$'s contained in $\Gamma$ have $codim(\pi_i(Y_j))>1$, thus implying $\Gamma_*=Z_*$ on $H^2_{\mathbb{Z}}$. We know those of codimension 1 are uniruled and effective, moreover it is known (see \Ref{prop}{birat_cone_unir}
 ) that uniruled divisors cut out the closure of the birational K\"{a}hler cone $\mathcal{BK}_X$, \ie $(\alpha,D)\geq0$ for all $\alpha\in\overline{\mathcal{BK}}_X$ and for all uniruled $D$. We wish to remark that the manifold $X$ we chose has $\overline{\mathcal{BK}}_X=\overline{\mathcal{C}}_X$ by \Ref{thm}{birat_cone} 
  (it contains no -2 divisors).\\ Let $\beta\in\mathcal{C}_X$ be a K\"{a}hler class and let $D\in Pic(X)$ be a uniruled divisor, we can write
\begin{equation}\nonumber
\beta=\alpha+\gamma,\,\,f(\alpha)\,\in\,T_G(L)\otimes\mathbb{R},\,\,f(\gamma)\,\in\,S_G(L)\otimes\mathbb{R}.
\end{equation}
Hence $0<(\beta,D)=(\gamma,D)$ and moreover we have $(f^{-1}\circ g\circ f(\beta),D)=(f^{-1}\circ g\circ f(\gamma),D)=(\gamma,f^{-1}\circ g^{-1}\circ f(D))\geq0$ because $f^{-1}\circ g\circ f(\beta)\in \overline{\mathcal{BK}}_X$ and $D$ is uniruled. Here is the contradiction: 
\begin{equation}\nonumber
0<(\beta,\sum_{h\in G}f^{-1}\circ h\circ f(D)),
\end{equation}
which implies $0\neq D'=\sum_{h\in G}hD\in f^{-1}(T_G(L)\cap S_G(L))=0$, hence there are no uniruled divisors inside $Pic(X)$. Moreover we obtain $\Gamma_*=Z_*$, \ie there exists a bimeromorphic map $\psi'$ of $X$ such that $\psi'^*=f^{-1}\circ g\circ f$ on $H^2(X)$. 
\end{proof}
\end{thm}


\begin{prop}\label{prop:sfiga}
Let $(S,G)$ be a couple consisting in a Leech-type lattice and its Leech automorphism group as in \Ref{defn}{leech_group}. Let moreover $S\subset N$, one of the 24 Niemeier lattices.\\ Suppose there exists a primitive embedding $S\rightarrow L$.\\
Then $G$ extends to a group of bimeromorphisms on some Hyperk\"{a}hler manifold $X$ of $K3^{[2]}$-type.
\begin{proof}
This is an immediate consequence of \Ref{thm}{cohom_to_aut}: $G$ acts trivially on $A_S$, therefore we can extend $G$ to a group of isometries of $L$ acting trivially on $S^{\perp_L}$. Thus we have $S_G(L)\cong S$. Moreover since $S$ is a Leech-type lattice contained in a negative definite lattice $N$ the other conditions of \Ref{thm}{cohom_to_aut} are satisfied.
\end{proof}
\end{prop}
We are now ready to prove the main result of this section:
\begin{thm}\label{thm:sporadic}
Let $X$ be a Hyperk\"{a}hler manifold of \ktipo and let $G$ be a finite group of symplectic automorphisms of $X$, then $G\subset Co_1$. 
\begin{proof}
Let $b=Rank(S_G(X))$, by \Ref{lem}{algaction} $S_G(X)$ has signature $(0,b)$. By \Ref{oss}{s_in_24} we have a lattice $T'$ of signature $(4,20-b)$ such that $A_{T'}=A_{S_G(X)}$ and $q_{T'}=-q_{A_{S_G(X)}}$. Therefore we can apply \Ref{lem}{nik_esiste} obtaining a lattice $T$ of signature $(0,24-b)$ and discriminant form $-q_{A_{S_G(X)}}$. Thus by \Ref{lem}{nik_immerge1} there exists a primitive embedding $S_G(X)\rightarrow N$, where $N$ is one of the lattices contained in \Ref{tab}{nieme}. Again by \Ref{lem}{algaction} we see that $(S_G(X),G)$ is a Leech couple, hence $G$ lies inside the Leech group of $N$. A direct computation using the 23 holy constructions shows that all these groups are contained in $Co_0$. Obviously the central involution of $Co_0$ has a co-invariant lattice of rank 24, hence we can restrict ourselves to $Co_1$. 
\end{proof}
\end{thm}
\subsection{Prime order symplectic automorphisms in the \ktipo case}
The aim of this subsection is to give a first application of \Ref{thm}{sporadic}, \ie the classification of prime order symplectic automorphisms on manifolds of \ktipo up to their fixed locus and their co-invariant lattice. Let us first give a bound to the possible prime orders: 
\begin{lem}\label{lem:max_p_order}
Let $\varphi$ be a symplectic automorphism of prime order $p$ on a \hk fourfold $X$ of \ktipo. Then $p\leq 11$.
\begin{proof}
By \Ref{thm}{sporadic} the order of a symplectic automorphism must divide the order of the group $Co_1$. That sorts out all primes apart for $2,3,5,7,11,13,23$. An automorphism of order 23 has a co-invariant lattice which is negative definite and of rank 22, therefore it cannot embed into $H^2(X,\mathbb{Z})$. This can be explicitly computed using an order 23 element of $M_{24}$ and letting it act on $\Lambda$ or on $N_{23}$, otherwise we can just rely on \Ref{ex}{cyclo_orderp} for a different method. The only Niemeier lattice with an automorphism of order 13 is $\Lambda$, where all elements of order 13 are conjugate. It is a well known fact that these automorphisms have no fixed points on $\Lambda$, as in \Ref{ex}{A122} or as in \Ref{ex}{craig_leech}.

\end{proof}
\end{lem}
Then we need to analyze what happens only for $p\leq 11$, however our result for $p=2$ can be proven separately without using \Ref{thm}{sporadic}:
\begin{prop}\label{prop:cohom_p2}
Let $X$ be a manifold of \ktipo and let $\varphi\subset Aut(X)$ be a symplectic involution. Then $X^{\varphi}$ consists of 28 isolated points and a $K3$ surface and $S_{\varphi}(X)\cong E_{8}(-2)$.
\begin{proof}
First of all by \Ref{thm}{fixed_p2} we have that $X^{\varphi}$ consists of 28 isolated points and a $K3$ surface and $S_{\varphi}(X)$ has rank 8. We now define an isometry $\overline{\varphi}$ of $L'$ as in \Ref{oss}{s_in_24} such that $S_{\varphi}(X)\cong S_{\overline{\varphi}}(X)$, hence $l(S_{\overline{\varphi}}(L'))\leq 8$, and so does its unimodular complement $T_{\overline{\varphi}}(L')$. This means that we can apply \Ref{lem}{nik_spezza} obtaining $T_{\overline{\varphi}}(L')=U\oplus T'$, which means that we can define an involution of $U^3\oplus E_8(-1)^2$ having $S_{\varphi}(X)$ as the anti-invariant lattice. By \Ref{lem}{algaction} this involution satisfies the conditions of  
\Ref{thm}{nik_k31} which implies that this involution on $U^3\oplus E_8(-1)^2$ is induced by a symplectic involution $\psi$ on some $K3$ surface $S$ and hence also $S_{\psi}(S)\,\cong\,S_{\varphi}(X)$.\\ Thus, by the work of Morrison on involutions \cite{mor}, we know $S_{\varphi}(X)=E_8(-2)$. 
\end{proof}
\end{prop}
Then we can proceed to prove our result.
\begin{thm}\label{thm:prime_autom_k32}
Let $\varphi$ be a symplectic automorphism of prime order $p$ on a \hk fourfold $X$ of \ktipo. Then the following holds:
\begin{table}[ht]\label{tab:prime_autom_tab}
\begin{tabular}{|c|c|c|}
\hline
$p$ & Fixed locus $X^{\varphi}$ & Lattice $S_{\varphi}(X)$\\
\hline
2 & $1$ K3 surface and 28 isolated points & $E_8(-2)$\\
\hline
3 & 27 isolated points & $K_{12}(-2)$, as in \Ref{ex}{cox-todd}\\
\hline
3 & 1 abelian surface & $W(-1)$ as in \Ref{ex}{wall}\\
\hline
5 & 14 isolated points & $S_{5.K3}$ as in \Ref{ex}{S_5K3}\\
\hline
7 & 9 isolated points & $S_{7.K3}$ as in \Ref{ex}{S_7K3}\\
\hline
11 & 5 isolated points & $S_{11.K3^{[2]}}$ as in \Ref{ex}{p11A124}\\
\hline
\end{tabular}
\end{table}
\begin{proof}
This result for $p=2$ is contained in \Ref{prop}{cohom_p2}.
By \Ref{thm}{sporadic} and its proof it is sufficient to look at all possible Leech couples $(S,\varphi)$ where $S$ is in a Niemeier lattice $N$ and $\varphi$ is a prime order isometry inside $Aut(N)/W(N)$. First of all let us work on the co-invariant lattice. For $p=3$ we obtain our result by \Ref{prop}{p3_nieme} and \Ref{thm}{fixedp3}. For $p=5$ this is proven by  \Ref{prop}{p5_nieme} and \Ref{thm}{fixedp5}. For $p=7$ this is precisely \Ref{prop}{p7_nieme}. Finally we already proved the result for $p=11$ in \Ref{ex}{p11A124} and \Ref{ex}{p11A212}. For the fixed locus we still need to prove the result for $p=3,7$ and $11$, however this is just \Ref{prop}{fixed_p3_sarti}, \Ref{prop}{fixed_p7_sarti} and \Ref{prop}{fixed_p11_sarti}.
\end{proof}
\end{thm}
As a consequence we have an improvement of \Ref{thm}{standard_morph}:
\begin{cor}\label{cor:prime_stand}
Let $(X,\varphi)$ be a couple consisting in a \hk manifold and a symplectic automorphism of order $i=2,3,5$. If $i=3$ suppose moreover that $X^\varphi$ consists of 27 isolated points. Then $(X,\varphi)$ is standard.
\begin{proof}
By \Ref{thm}{prime_autom_k32} these hypothesis are equivalent to those of \Ref{thm}{standard_morph}, therefore the claim holds.
\end{proof}
\end{cor}

\section{The $K3^{[n]}$-type case}\label{sec:k3n_case}
In this section we specialize to the case of manifolds of $K3^{[n]}$-type. The results in this section are similar to those of \Ref{sec}{k32_case}, however we cannot compute the fixed locus of prime order symplectic automorphisms since the computations of \Ref{cap}{fixed} are not possible in this more general setting. We keep the notation $L_n\cong H^2(X,\mathbb{Z})$ for $X$ a manifold of \kntipo. Let us recall that also in this case the map \eqref{numap} is injective.
Let us recall that in \Ref{defn}{numerical_exc} we defined a class $\mathcal{NE}xc$ of numerically exceptional divisors which have an effective power by \Ref{thm}{num_eff}. Let us give a further definition:
\begin{defn}
Let $X$ be a manifold of \kntipo and let $ME(X)$\index{Exceptional classes, Multiples of exceptional divisors, $ME(X)$} be the set of elements $v$ in $S(X)$ such that $mv\in\mathcal{NE}xc(X)$ for some $m$.
\end{defn}
\begin{lem}\label{lem:algaction_k3n}
Let $X$ be a manifold of \kntipo. Let moreover $G\subset Aut_s(X)$ be a finite group. Then the following hold:
\begin{itemize}
\item $S_G(X)$ is nondegenerate and negative definite.
\item $S_G(X)\cap ME(X)=\emptyset$.
\item $T(X)\subset T_G(X)$ and $S_G(X)\subset S(X)$.
\item Suppose moreover that $G$ is discriminant preserving, then $G$ acts trivially on $A_{S_G(X)}$.
\end{itemize}
\begin{proof}
Apart for the second assertion everything has been proven in \Ref{lem}{algaction_general}. The elements of $ME(X)$ have a multiple which is effective, therefore we can reason as in \Ref{lem}{algaction} to conclude that they cannot be inside $S_G(X)$.
\end{proof}
\end{lem}

\begin{thm}\label{thm:bir_to_aut_k3n}
Let $n\geq 2$. Let $L_n$ be as above and let $G$ be a finite subgroup of $O(L_n)$. Suppose that the following hold:
\begin{enumerate}
\item $S_G(L_n)$ is nondegenerate and negative definite.
\item $S_G(L_n)\cap h(ME(Y))=\emptyset$ for all marked manifolds $(Y,h)$ of \kntipo.
\end{enumerate}
Then $G$ is induced by a subgroup of $Bir_s(X)$ for some manifold $(X,f)$ of \kntipo.
\begin{proof}
We can proceed as in \Ref{thm}{cohom_to_aut} choosing $(X,f)$ as a manifold with $f(T(X))=T_G(L_n)$. Again we have $\overline{\mathcal{BK}}_X=\overline{\mathcal{C}}_X$ by \Ref{thm}{birat_cone} and we prove in the same way that all elements of $G$ are induced by a (symplectic) birational morphism of $X$.
\end{proof}
\end{thm}
 
\begin{thm}
Let $X$ be a manifold of \kntipo and let $G\subset Aut_s(X)$ be a finite discriminant preserving group. Then there exists and embedding $S_G(X)\rightarrow N$, where $N$ is one of the 24 Niemeier lattices and $G$ extends to a group of isometries of $N$ acting trivially on $S_G(X)^{\perp_N}$. Moreover $G\subset Leech(N)\subset Co_1$.
\begin{proof}
The proof goes as in \Ref{thm}{sporadic}. Since we have no $-2$ vectors inside $S_G(X)$ (they would be in $ME(X)$), we obtain $G\subset Leech(N)\subset Co_1$.
\end{proof}
\end{thm}


\begin{cor}\label{cor:k3n_prime}
Let $X$ be a manifold of \kntipo and let $\varphi\in Aut_s(X)$ be of prime order $p\neq 2$. Then one of the following holds:
\begin{table}[ht]\label{tab:prime_autom_k3n_tab}
\begin{tabular}{|c|c|}
\hline
$p$ & Lattice $S_{\varphi}(X)$\\
\hline
3 &  $S_{3.exo}$, as defined in \Ref{ex}{p3E83} \\
\hline
3 &  $K_{12}(-2)$, where $K_{12}$ is as in \Ref{ex}{cox-todd}\\
\hline
3 &  $W(-1)$ as in \Ref{ex}{wall}\\
\hline
5 &  $S_{5.K3}$ as in \Ref{ex}{S_5K3}\\
\hline
5 &  $S_{5.exo}$ as in \Ref{ex}{slat_5rk4}\\
\hline
7 &  $S_{7.K3}$ as in \Ref{ex}{S_7K3}\\
\hline
11 & $S_{11.K3^{[2]}}$ as in \Ref{ex}{p11A124}\\
\hline
\end{tabular}
\end{table}
\begin{proof}
This is just a direct consequence of \Ref{prop}{p3_nieme}, \Ref{prop}{p5_nieme}, \Ref{prop}{p7_nieme} and \Ref{ex}{p11A124}. 
\end{proof}
\end{cor}
\begin{oss}\label{oss:n_repr}
A comparison between \Ref{cor}{k3n_prime} and \Ref{thm}{prime_autom_k32} shows that the cases in the corollary might not happen in some dimensions, let us look a little more into this.
Obviously all the cases corresponding to standard automorphisms exist in all possible dimensions, as the examples \ref{ex:3_stand}, \ref{ex:5_stand} and \ref{ex:7_stand} show.
To analyze all of the other cases we must embed the lattices $S_i$ contained in \Ref{tab}{prime_autom_k3n_tab} inside the Mukai lattice $L'$ and look at their orthogonal: if it represents the integer $2(n-1)$ with a primitive vector then there exists a primitive embedding  of $S_i$ inside $L_n$, \ie by \Ref{thm}{bir_to_aut_k3n} there exists a \hk manifold of \kntipo having a birational morphism $\varphi$ such that $S_{\varphi}\cong S_i$. We will compute what happens up to $n=9$, let us look at all cases one by one:
\begin{itemize}
\item[$S_i=W(-1)$] In this case we will look at a greater lattice: let $F$ be the orthogonal inside $\Lambda$ to the $\mathcal{S}$-lattice $2^93^6$. Then $W(-1)\subset F$ by $2^93^6\subset 2^{27}3^{36}$. Let us now embed $F$ into $L'$ and let $T=F^{\perp_L}$. A necessary condition for $W(-1)\rightarrow L_n$ is that $2(n-1)$ is represented by a primitive vector of $T$. By \Ref{ex}{slat_3rk4} $T\cong A_2\oplus A_2(3)$ and a direct computation shows that it represents the integers $2,6,8$ and $14$, therefore $W(-1)$ primitively embeds into $L,L_4,L_5$ and $L_8$ but might not embed into $L_3,L_6$ and $L_7$.
\item[$S_i=S_{3.exo}$] Let us fix an embedding $S_i\rightarrow L'$ and let $T=S_i^\perp$. Suppose that the integer $r$ is represented by a primitive element of $T$, then the lattice $T'=r^{\perp_T}$ exists, which in turn implies $l(A_{T'})\leq 7$, therefore by \Ref{lem}{nik_immerge} $r$ must be a multiple of 3, \ie $S_i\rightarrow L_n$ implies $n=3m+1$. 
\item[$S_i=S_{5.exo}$] Let us fix an embedding $B\rightarrow L'$ and let $T=B^\perp$. By \Ref{ex}{slat_5rk4} $T$ is in the same genus of the $\mathcal{S}$-lattice $2^53^{10}(-1)$, however there is only one lattice in this genus, which we recall is
\begin{equation*}
\left(
\begin{array}{cccc}
4 & 1 & 1 & -1\\
1 & 4 & -1 & 1\\
1 & -1 & 4 & 1\\
-1 & 1 & 1 & 4
\end{array}
\right).
\end{equation*}
A direct computation shows that its primitive vectors represent the integers $4,6,10,12,14$ and $16$ but not $8$, therefore it does not embed into $L_5$. 
\item[$S_i=S_{11.K3^{[2]}}$] Let us fix an embedding $S_{11.K3^{[2]}}\rightarrow L'$ and let $T=S_{11.K3^{[2]}}^\perp$. $T$ has determinant 121 and, by \cite{nip}, there is only one genus of such lattices, containing the following:
\begin{equation*}
\left(
\begin{array}{cccc}
4 & 2 & 1 & 0\\
2 & 4 & 1 & 1\\
1 & 1 & 4 & 2\\
0 & 1 & 2 & 4
\end{array}
\right),\,\,\,
\left(
\begin{array}{cccc}
2 & 1 & 1 & 0\\
1 & 2 & 1 & 1\\
1 & 1 & 8 & 4\\
0 & 1 & 4 & 8
\end{array}
\right),\,\,\,
\left(
\begin{array}{cccc}
2 & 0 & 1 & 0\\
0 & 2 & 0 & 1\\
1 & 0 & 6 & 0\\
0 & 1 & 0 & 6
\end{array}
\right).
\end{equation*}
A direct computation shows that the integers $2,4,6,8,10,12,14$ and $16$ are represented by these lattices, therefore $S_{11.K3^{[2]}}$ embeds into $L,L_3,L_4,L_5,L_6,L_7,L_8$ and $L_9$.
\end{itemize}
\end{oss}
\section{Examples revisited}\label{sec:exa2}
In this section we look back at the examples of \Ref{cap}{exoexa} and we use our classification results to compute the Picard lattice of several of our examples.
\subsection{Standard automorphisms}
Let us look at the automorphisms defined in the examples \ref{ex:fano_invol}, \ref{ex:fano_3stand}, \ref{ex:fano_3stand2}, \ref{ex:fano5} and \ref{ex:epw_5}. By \Ref{cor}{prime_stand} all of these automorphisms are standard. We can easily compute the Neron-Severi lattice $S(X)$ for the generic elements of the families in the above examples. Let $X$ be a generic Fano scheme of lines defined in \Ref{ex}{fano_invol} and let $\varphi$ be its symplectic involution. Then $S(X)$ is an overlattice of $(6)\oplus E_8(-2)$ by \Ref{thm}{prime_autom_k32} and \Ref{oss}{simple_sympl}. However there are no nontrivial overlattices of it, therefore we have $S(X)=(6)\oplus E_8(-2)$.\\ 
Let $X$ be a generic Fano scheme of lines satisfying the hypothesis of \Ref{ex}{fano_3stand} or of \Ref{ex}{fano_3stand2} and let $\varphi$ be its symplectic automorphism of order 3. Then $S(X)$ is an overlattice of $(6)\oplus K_{12}(-2)$ by \Ref{thm}{prime_autom_k32} and \Ref{oss}{simple_sympl}. However also in this case there are no nontrivial overlattices of it, therefore we have $S(X)=(6)\oplus K_{12}(-2)$.\\
Analogously let $X$ be a generic Double-EPW-sextic satisfying the hypothesis of \Ref{ex}{epw_3stand} and let $\varphi$ be its symplectic automorphism of order 3. Then $S(X)$ is an overlattice of $(2)\oplus K_{12}(-2)$ by \Ref{thm}{prime_autom_k32} and \Ref{prop}{EPW_autom}. However also in this case there are no nontrivial overlattices of it, therefore we have $S(X)=(2)\oplus K_{12}(-2)$.\\
Let $X$ be a generic Fano scheme of lines satisfying the hypothesis of \Ref{ex}{fano5} and let $\varphi$ be its symplectic automorphism of order 5. Then $S(X)$ is an overlattice of $(6)\oplus S_{5.K3}$ by \Ref{thm}{prime_autom_k32} and \Ref{oss}{simple_sympl}. However also in this case there are no nontrivial overlattices of it, therefore we have $S(X)=(6)\oplus S_{5.K3}$.\\
Finally let $X$ be a generic Double-EPW-sextic satisfying the hypothesis of \Ref{ex}{epw_5} and let $\varphi$ its symplectic automorphism of order 5. Then $S(X)$ is an overlattice of $(2)\oplus S_{5.K3}$ by \Ref{thm}{prime_autom_k32} and \Ref{prop}{EPW_autom}. However also in this case there are no nontrivial overlattices of it, therefore we have $S(X)=(2)\oplus S_{5.K3}$.\\
 
\subsection{Fano scheme of lines on Fermat's cubic}
Let us look back at the Fano scheme of lines $F$ defined in \Ref{ex}{fanoA6}, and let $G=(\mathbb{Z}_{/(3)})^4.A_6$ be the group of symplectic automorphisms of $F$ contained in that example. By \Ref{thm}{sporadic} we can evaluate the action of $G$ on $H^2(X,\mathbb{Z})$ by looking at its action on the Niemeier lattices $N$. Since $G\subset Leech(N)$ only if $N=\Lambda$ and since there is only one conjugacy class of $G$ inside $Co_1$ we have that $S_G(F)\cong (2^93^6)^{\perp_\Lambda}$, where $2^93^6$ is the $\mathcal{S}$-lattice of \Ref{ex}{slat_3rk4}. A direct computation using primitive embeddings of $S_G(F)$ into the Mukai lattice $L'$ as in \Ref{oss}{n_repr} shows that $T_G(F)=(6)\oplus A_2(3)$. Therefore we have $T(F)\cong A_2(3)$ and $S(F)\cong (-2)\oplus A_2(-3)\oplus U\oplus E_8(-1)^2$.

\subsection{Exotic automorphism of order 3}
In our classification of \Ref{thm}{prime_autom_k32} we proved that there exists an automorphism $\varphi$ of order 3 which fixes an abelian surface. Moreover it is defined on manifolds belonging to a 3-dimensional subset $\mathcal{M}_{3.exo}$ of the Moduli space of manifolds of \ktipo. A projective example of these automorphisms is given by \Ref{ex}{fano3}. We also proved that it has $S_\varphi(X)\cong W(-1)$ for all $X\in \mathcal{M}_{3.exo}$. Therefore $T_\varphi(X)$ has discriminant $2\cdot3^5$ and its discriminant group is $\mathbb{Z}_{/(6)}\times (\mathbb{Z}_{/(3)})^4$ which has 5 generators. Let us suppose that $X$ is projective: thus there exists a polarization $v\in T_\varphi(X)$. However $v^2$ must be a multiple of $3$, otherwise $v^\perp\subset T_\varphi(X)$ would have a discriminant group with 5 generators and rank 4, which is impossible. This implies that automorphisms of order 3 which fix a surface cannot be found on Double EPW-sextics, varieties of sums of powers or subspaces of the Grassmannian where the automorphism preserves their natural polarization (\ie it is induced by an automorphism of $\mathbb{P}^5$ or of the Grassmannian respectively). This implies that the double EPW-sextics in \Ref{ex}{epw_3exo} 
cannot be resolved while preserving their automorphisms.

\subsection{Automorphisms of order 11}
In this subsection we analyze deformation classes of manifolds of \ktipo with a symplectic automorphism of order 11 and we look at their possible invariant polarizations.
Let $X$ be a manifold of \ktipo with a symplectic automorphism $\psi$ of order 11 and let $\omega$ be a $\psi$-invariant \kahl class. First of all let us remind that \Ref{thm}{prime_autom_k32} implies that non trivial deformations of $(X,\psi)$ are  of maximal dimension 1, moreover the twistor family $TW_{\omega}(X)$ is naturally endowed with a symplectic automorphism of order 11 as in \Ref{oss}{twistor_deform}.   
Therefore $TW_{\omega}(X)$ is already a family of the maximal dimension for such manifolds $(X,\psi)$. Moreover we have that the twistor family $TW_{\omega}(X)$ is actually a family over the base $\mathbb{P}(T_{\psi}(X)\otimes\mathbb{R})$ since $T_{\psi}(X)=<\omega,\sigma_X,\overline{\sigma}_X>\cap H^2_{\mathbb{Z}}(X)$.\\ Thus what we really need to analyze are the possible lattices $T_{\psi}(X)$ up to isometry. We have already proved that there exists only one isometry class of lattices $S_{\psi}(X)$. However there might be several isomorphism classes of lattices $T_{\psi}(X)$. In fact \Ref{thm}{sporadic} and \Ref{prop}{sfiga} can be used only to compute its genus.\\ A direct computation shows that there are two such lattices, namely the following:
\begin{equation}\label{TX_1}
T^1_{11}=\left(\begin{array}{ccc}  2&1&0\\ 1&6&0\\ 0&0&22\end{array}\right),
\end{equation}
\begin{equation}\label{TX_2}
T^2_{11}=\left(\begin{array}{ccc} 6&-2&-2\\ -2& 8&-3\\-2&-3&8\end{array}\right).
\end{equation}
Therefore there are 2 distinct families of \hk manifolds endowed with a symplectic automorphism of order 11, let us call $TW(X_1)$ the first and $TW(X_2)$ the second.
Another direct computation shows that $T_\psi(X)$ has a primitive element of square $2$ only if we are in case \eqref{TX_1}, therefore the Double-EPW-sextic of \Ref{ex}{epw_11} belongs to $TW(X_1)$. Moreover there is an element of square $6$ and divisibility $2$ in $T_\psi(X)$ only in case \eqref{TX_2}, therefore the Fano scheme of lines of \Ref{ex}{fano11} belongs to $TW(X_2)$. 
\backmatter
\chapter*{Acknowledgements}
First and foremost I would like to thank my advisor, Prof. K. G. O'Grady for his guidance during the course of my PhD, for suggesting me this problem and for his help regarding the present work. Then I would like to thank Prof. E. Sernesi, my tutor at RomaTRE university, for his help in the early stages of my PhD and for his counsel concerning deformation theory. I am also very grateful to A. Rapagnetta for his patience and for his comments on an early version of this work. I would like to thank Prof. S. Boissi\'{e}re, C. Camere and Prof. A. Sarti for their comments and for useful discussions concerning \Ref{cap}{fixed} and \Ref{sec}{k32_case}. I am grateful to Prof. D. Huybrechts for his comments on symplectic involutions and to Prof. V. V. Nikulin for a motivating discussion. I would also like to thank Prof. M. Sch\"{u}tt and X. Roulleau for their interest in \Ref{ex}{fano11} and for related interesting discussions. Finally I would like to thank also Prof. B. Van Geemen for communicating me a proof of \Ref{prop}{vgir_mail} and M. Lelli Chiesa and G. Sacc\'{a} for useful discussions concerning \hk geometry.\\
I am also indebted to the University of RomaTRE and the department of Mathematics for funding my Phd, to the I.N.d.A.M for funding most of the schools, workshops and conferences I have attended, thus allowing me to meet a wide part of the mathematical community and to exchange ideas with them. For the same reason I would like to thank also the Field's Institute of Toronto, The Newton's Institute of Cambridge, The Banff center, the I.C.T.P. of Trieste, Luminy's CIRM and Trento's CIRM.\\ 
I would like also to add special thanks to the referees and the members of the commission of this thesis defence, Prof. A. Lopez, Prof. M. Lehn and Prof. B. Van Geemen, for their helpful comments.

\printindex
\end{document}